\documentclass[10pt]{article}

\usepackage{amsmath}
\usepackage{amsthm}
\usepackage{amssymb}
\usepackage{amscd}
\usepackage{amsfonts}  
\usepackage{makeidx}
\usepackage{hyperref}
\usepackage{makeidx}     
\usepackage{graphicx}    
\usepackage{multicol}    

\DeclareMathAlphabet\EuR{U}{eur}{m}{n}
\SetMathAlphabet\EuR{bold}{U}{eur}{b}{n}

\makeindex             

\newcommand{\indexnotation}[1]{\label{#1}}

\newcommand{\scal}[2]{\left\langle #1,#2 \right\rangle}

\newcommand{\squarematrix}[4]                
{                                            
\left( \begin{array}{cc} #1 & #2 \\ #3 &
#4
\end{array} \right)
}

\newcommand{\caleam}{{\mathcal E}\!{\mathcal A}\!{\mathcal M}}

\newcommand{\calfin}{{\mathcal F}\!{\mathcal I}\!{\mathcal N}}

\newcommand{\cala}{{\cal A}}
\newcommand{\calb}{{\cal B}}
\newcommand{\calc}{{\cal C}}
\newcommand{\cald}{{\cal D}}

\newcommand{\calf}{{\cal F}}
\newcommand{\calg}{{\cal G}}
\newcommand{\calh}{{\cal H}}

\newcommand{\caln}{{\cal N}}

\newcommand{\calr}{{\cal R}}

\newcommand{\calu}{{\cal U}}

\newcommand{\calz}{{\cal Z}}

\newcommand{\bbC}{{\mathbb C}}

\newcommand{\bbH}{{\mathbb H}}

\newcommand{\bbQ}{{\mathbb Q}} 
\newcommand{\bbR}{{\mathbb R}}

\newcommand{\bbZ}{{\mathbb Z}}

\newcommand{\bfP}{{\mathbf P}}

\newcommand{\bfT}{{\mathbf T}}

\newcommand{\Alt}{\operatorname{Alt}}

\newcommand{\aut}{\operatorname{aut}}

\newcommand{\coker}{\operatorname{coker}}
\newcommand{\colim}{\operatorname{colim}}

\newcommand{\con}{\operatorname{con}}

\newcommand{\defi}{\operatorname{def}}
\newcommand{\dom}{\operatorname{dom}}

\newcommand{\frk}{\operatorname{f-rk}}

\newcommand{\id}{\operatorname{id}}
\newcommand{\im}{\operatorname{im}}

\newcommand{\Isom}{\operatorname{Isom}}

\newcommand{\Rad}{\operatorname{Rad}}
\newcommand{\res}{\operatorname{res}}
\newcommand{\rk}{\operatorname{rk}}
\newcommand{\vol}{\operatorname{vol}}
\newcommand{\sign}{\operatorname{sign}}
\newcommand{\sing}{\operatorname{sing}}

\newcommand{\Tor}{\operatorname{Tor}}
\newcommand{\tors}{\operatorname{tors}}
\newcommand{\tr}{\operatorname{tr}}

\newcommand{\virt}{\operatorname{virt}}

\newcommand{\Wh}{\operatorname{Wh}}

\newtheorem{theorem}{Theorem}[section]
\newtheorem{lemma}[theorem]{Lemma}

\newtheorem{definition}[theorem]{Definition}

\newtheorem{example}[theorem]{Example}
\newtheorem{remark}[theorem]{Remark}

\newtheorem{conjecture}[theorem]{Conjecture}

\newtheorem{slogan}[theorem]{Slogan}
{\catcode`@=11\global\let\c@equation=\c@theorem}

\renewcommand{\labelenumi}{(\roman{enumi})}

\newcommand{\entry}[2]{#1, ~ \pageref{#2}}   
                                             
\newcommand{\EGF}[2]{E_{#2}(#1)}               

\newcommand{\comsquare}[8]                   
{\begin{CD}
#1 @>#2>> #3\\
@V{#4}VV @VV{#5}V\\
#6 @>>#7> #8
\end{CD}
}

\newcommand{\action}{\curvearrowright} 

\begin{document}

\typeout{----------------------------  ltwoalg.tex  ----------------------------}

\title{$L^2$-Invariants from the Algebraic Point of View}
\author{
Wolfgang L\"uck\thanks{\noindent email:
lueck@math.uni-muenster.de\protect\\
www: ~http://www.math.uni-muenster.de/u/lueck/\protect\\
FAX: 49 251 8338370\protect}\\
Fachbereich Mathematik\\ Universit\"at M\"unster\\
Einsteinstr.~62\\ 48149 M\"unster\\Germany}
\maketitle


\typeout{-----------------------  Abstract  ------------------------}

\begin{abstract}
We give a survey on $L^2$-invariants such as $L^2$-Betti numbers and
$L^2$-torsion taking an algebraic point of view. We discuss their
basic definitions, properties and applications to problems arising in
topology, geometry, group theory and $K$-theory. 
\end{abstract}

\\[1mm]


\typeout{--------------------   Section 0: Introduction --------------------------}
\setcounter{section}{-1}
\section{Introduction}
\label{sec: Introduction}

The purpose of this survey article is to present an algebraic approach
to $L^2$-invariants such as $L^2$-Betti numbers and $L^2$-torsion.
Originally these were defined analytically  in terms of heat
kernels. After it was discovered that they have simplicial and
homological algebraic counterparts, there have been many applications to various
problems in topology, geometry, group theory and algebraic $K$-theory,
which on the first glance do not involve any $L^2$-notions. Therefore
it seems to be useful to give a quick and friendly introduction to
these notions in particular for mathematicians who have more algebraic than
analytic background. This does not at all mean that the analytic
aspects are less important, but for certain applications it is not
necessary to know the analytic approach and it is possible and easier
to focus on the algebraic aspects. Moreover, questions about
$L^2$-invariants of heat kernels such as the Atiyah Conjecture or the
zero-in-the-spectrum-Conjecture turn
out to be strongly related to algebraic questions about modules over
group rings.

The hope of the author is that more people take notice of 
$L^2$-invariants and $L^2$-methods, and may be able to apply them to their
favourite problems, which not necessarily come a priori from an
$L^2$-setting. Typical examples of such instances will be discussed
in this survey article. There are many open questions and conjectures
which have the potential to stimulate further activities.

The author has tried to write this article in a way which makes it
possible to quickly pick out specific topics of interest and read them
locally without having to study too much of the previous text. 

These notes are based on a series of lectures which were presented by the
author at the LMS Durham Symposium on Geometry and Cohomology in Group Theory
in July 2003. The author wants to thank the organizers Martin Bridson,
Peter Kropholler and Ian Leary and the London Mathematical Society for
this wonderful symposium and Michael Weiermann for proof reading the manuscript.

In the sequel ring will always mean associative ring with unit
and $R$-module will mean left $R$-module unless explicitly
stated differently. The letter $G$ denotes a discrete group.
Actions of $G$ on spaces are always from the left.

\tableofcontents


\typeout{--------------------   The Group von Neumann Algebra --------------------------}

\section{Group von Neumann Algebras}
\label{sec: Group von Neumann Algebras}

The integral group ring $\bbZ G$ plays an important role in topology and geometry,
since for a $G$-space its singular chain complex or for a $G$-$CW$-complex its 
cellular chain complex are $\bbZ G$-chain complexes. 
However, this ring is rather complicated and
does not have some of the useful properties which  
other rings such as  fields or semisimple rings have.
Therefore it is very hard to analyse modules over $\bbZ G$. Often in algebra one studies
a complicated ring by investigating certain localizations or completions of it which do have
nice properties. They still contain and focus on useful information about 
the original ring, which now becomes accessible.
Examples are the quotient field of an integral domain, 
the $p$-adic completion of the integers
or the algebraic closure of a field. In this section 
we present a kind of completion of the complex
group ring $\bbC G$ given by the group von Neumann algebra 
and discuss its ring theoretic properties.


\subsection{The Definition of the Group von Neumann Algebra}
\label{sec: The Definition of the Group von Neumann Algebra}

Denote by $l^2(G)$%
\indexnotation{l^2(G)}
the Hilbert space $l^2(G)$ consisting of  formal sums
$\sum_{g \in G} \lambda_g \cdot g$
for complex numbers $\lambda_g$ such that
$\sum_{g \in G} |\lambda_g|^2 < \infty$. The scalar
product is defined by
\begin{eqnarray*}
\scal{\sum_{g \in G} \lambda_g \cdot g}
{\sum_{g \in G} \mu_g \cdot g}
& := &
\sum_{g \in G} \lambda_g \cdot \overline{\mu_g}.
\end{eqnarray*}
This is the same as the Hilbert space completion of the complex group ring
$\bbC G$  with respect to the pre-Hilbert space structure
for which $G$ is an orthonormal basis.
Notice that left multiplication with elements in $G$ induces
an isometric $G$-action on $l^2(G)$. Given a Hilbert space $H$,
denote by $\calb(H)$%
\indexnotation{calb(H)}
the $C^*$-algebra of bounded (linear) operators from $H$ to itself,
where the norm is the operator norm and the involution is given by
taking adjoints.

\begin{definition}[Group von Neumann algebra]
\label{def: group von Neumann algebra}
The \emph{group von Neumann algebra}%
\index{group von Neumann algebra}
\index{von Neumann algebra!of a group}
$\caln(G)$%
\indexnotation{caln(G)} of the group $G$ is defined
as the algebra of $G$-equivariant bounded operators
from $l^2(G)$ to $l^2(G)$
\begin{eqnarray*}
\caln(G)
& := &
\calb(l^2(G))^G.
\end{eqnarray*}
\end{definition}

In the sequel we will view the complex group ring $\bbC G$ as a subring
of $\caln(G)$ by the embedding of $\bbC$-algebras
$\rho_r \colon \bbC G \to \caln(G)$ which sends $g \in G$ to the $G$-equivariant operator
$r_{g^{-1}} \colon l^2(G) \to l^2(G)$ given by right multiplication with $g^{-1}$.

\begin{remark}[The general definition of von Neumann algebras]
\label{rem: definition of a von Neumann algebra}
\em
In general a \emph{von Neumann algebra}%
\index{von Neumann algebra}
$\cala$%
\indexnotation{cala}
 is a sub-$\ast$-algebra of $\calb(H)$
for some Hilbert space $H$,
which is closed in the weak topology and contains $\id\colon H \to H$. Often in the literature
the group von Neumann algebra
$\caln(G)$ is defined as the closure in the  weak topology of the complex group ring
$\bbC G$ considered as $\ast$-subalgebra of $\calb(l^2(G))$. This definition and
Definition~\ref{def: group von Neumann algebra} agree 
(see \cite[Theorem 6.7.2 on page 434]{Kadison-Ringrose(1986)}).
\em
\end{remark}

\begin{example}[The von Neumann algebra of a finite group] 
\label{exa: group von Neumann algebra of a finite group} \em
If $G$ is finite, then nothing happens, namely $\bbC G = l^2(G) = \caln(G)$.
\em
\end{example}

\begin{example}[The von Neumann algebra of $\bbZ^n$] 
\label{exa: group von neumann algebra of Z^n} \em
In general there is no concrete model for $\caln(G)$.
However, for $G = \bbZ^n$, there is the following illuminating model
for the group von Neumann algebra $\caln(\bbZ^n)$.
Let $L^2(T^n)$%
\indexnotation{L^2(T^n)}
be the Hilbert space of equivalence classes of $L^2$-integrable
complex-valued functions on the $n$-dimensional torus $T^n$,
where two such functions are
called equivalent if they differ only on a subset of measure zero.
Define the ring $L^{\infty}(T^n)$%
\indexnotation{L^{infty}(T^n)}
by equivalence classes of essentially bounded measurable functions
$f\colon T^n \to \bbC$,
where essentially bounded%
\index{function!essentially bounded}
means that there is a constant
$C > 0$ such that the set $\{x \in T^n\mid |f(x)| \ge C\}$
has measure zero. An element $(k_1, \ldots , k_n)$
in $\bbZ^n$ acts isometrically on $L^2(T^n)$ by pointwise multiplication with
the function $T^n \to \bbC$, which maps
$(z_1, z_2, \ldots, z_n)$ to $z_1^{k_1} \cdot \ldots \cdot z_n^{k_n}$.
Fourier transform yields an
isometric $\bbZ^n$-equivariant isomorphism
$l^2(\bbZ^n) \xrightarrow{\cong} L^2(T^n)$.
Hence $\caln(\bbZ^n) =  \calb(L^2(T^n))^{\bbZ^n}$.
We obtain an isomorphism (of $C^*$-algebras)
$$L^{\infty}(T^n) \xrightarrow{\cong} \caln(\bbZ^n)$$
by sending $f \in L^{\infty}(T^n)$ to
the $\bbZ^n$-equivariant operator
$$M_f\colon L^2(T^n) \to L^2(T^n),
\quad g \mapsto g \cdot f,$$
where $g\cdot f(x)$ is defined by $g(x)\cdot f(x)$. \em
\end{example}

Let $i \colon H \to G$ be an injective group homomorphism.
It induces a ring homomorphism $\bbC i \colon \bbC H \to \bbC G$,
which extends to a ring homomorphism
\begin{eqnarray}
&\caln(i)%
\indexnotation{induced ring homomorphism caln(i)}
\colon \caln(H) \to  \caln(G) &
\label{caln(i)}
\end{eqnarray}
as follows. Let $g \colon l^2(H) \to l^2(H)$ be a $H$-equivariant bounded 
operator. Then $\bbC G \otimes_{\bbC H} l^2(H) \subseteq l^2(G)$ is a dense 
$G$-invariant subspace and  
$$\id_{\bbC G}\otimes_{\bbC H} g \colon \bbC G \otimes_{\bbC H} l^2(H) \to 
\bbC G \otimes_{\bbC H} l^2(H)$$  
is a $G$-equivariant linear map, which is bounded with respect to the norm coming from
$l^2(G)$. Hence it induces a $G$-equivariant bounded operator
$l^2(G) \to l^2(G)$,
which is by definition the image of $g \in \caln(H)$ under $\caln(i)$.

In the sequel we will ignore the functional analytic aspects of
$\caln(G)$ and will only consider its algebraic properties as a ring.


\subsection{Ring Theoretic Properties of the Group von Neumann Algebra}
\label{sec: Ring Theoretic Properties of the Group von Neumann Algebra}

On the first glance the von Neumann algebra $\caln(G)$ looks not very nice as a ring.
It is an \emph{integral domain,}%
\index{integral domain}
\index{ring!integral domain}
i.e. has no non-trivial zero-divisors if and only if $G$ is trivial.
It is Noetherian if and only if $G$ is finite (see \cite[Exercise 9.11]{Lueck(2002)}).
It is for instance easy to see
that $\caln(\bbZ^n) \cong L^{\infty}(T^n)$ does 
contain non-trivial zero-divisors and is not Noetherian. The main advantage of
$\caln(G)$ is that it contains many more idempotents than $\bbC G$. This has the effect that
$\caln(G)$ has the following ring theoretic property. A ring $R$ is called
\emph{semihereditary}%
\index{semihereditary ring}
\index{ring!semihereditary}
if every finitely generated submodule of a projective module is  again projective.
This implies that the category of finitely presented $R$-modules is an abelian category.

\begin{theorem}[Von Neumann algebras are semihereditary]
\index{Theorem!Von Neumann Algebras are semihereditary}
\label{the: von Neumann algebras are semihereditary}
Any von Neumann algebra $\cala$ is semihereditary.
\end{theorem}
\proof This follows from
the facts that any von Neumann algebra is a Baer $*$-ring and hence in 
particular a Rickart $C^*$-algebra 
\cite[Definition 1, Definition 2 and Proposition 9 in 
Chapter 1.4]{Berberian(1972)} and that a $C^*$-algebra is 
semihereditary if and only if it is
Rickart \cite[Corollary 3.7 on page 270]{Ara-Goldstein(1993)}.  \qed

\begin{remark}[Group von Neumann algebras are semihereditary]
\em \label{rem: elementary proof of semihereditary}
It is quite useful to study the following elementary proof 
of Theorem~\ref{the: von Neumann algebras are semihereditary} 
in the special case of a group von Neumann algebra $\caln(G)$.
One easily checks that it suffices to show for a finitely generated
submodule $M \subseteq \caln(G)^n$ that $M$ is projective.
Let $f \colon \caln(G)^m \to \caln(G)^n$ be an $\caln(G)$-linear map.
Choose a  matrix $A \in M(m,n;\caln(G))$ such that $f$ is given by right multiplication with 
$A$. Because of $\caln(G) = \calb(l^2(G))^G$ we can define a $G$-equivariant bounded operator
$$\nu(f)  \colon l^2(G)^m \to l^2(G)^n, \quad 
(u_1, \ldots, u_m) ~ \mapsto ~ 
\left(\sum_{i = 1}^m \overline{a_{i,1}^*(\overline{u_i})}, \ldots , 
\sum_{i = 1}^m \overline{a_{i,n}^*(\overline{u_i})}\right),$$
where by definition
$\overline{\sum_{g \in G} \lambda_g \cdot g} :=  \sum_{g \in G} \overline{\lambda_g} \cdot g$
and $a_{i,j}^*$ denotes the adjoint of $a_{i,j}$.
With these conventions $\nu(\id) = \id$, 
$\nu(r\cdot f + s \cdot g) = r \cdot \nu(f) + s \cdot \nu(g)$ and 
$\nu(g \circ f) = \nu(g) \circ \nu(f)$ for $r,s \in \bbC$ and 
$\caln(G)$-linear maps $f$ and $g$. Moreover we have $\nu(f)^* = \nu(f^*)$ 
for an $\caln(G)$-map 
$f \colon \caln(G)^m \to \caln(G)^n$, where $f^* \colon \caln(G)^n \to \caln(G)^m$ 
is given by right multiplication with the matrix $(a_{j,i}^*)$, 
if $f$ is given by  right multiplication with the matrix $(a_{i,j})$, 
and $\nu(f)^*$ is the adjoint of the operator $\nu(f)$.

Every equivariant bounded operator
$l^2(G)^m \to l^2(G)^n$ can be written as $\nu(f)$ for a unique $f$. Moreover, 
the sequence $\caln(G)^m \xrightarrow{f} \caln(G)^n \xrightarrow{g} \caln(G)^p$ of $\caln(G)$-modules 
is exact if and only if the sequence of bounded $G$-equivariant operators 
$l^2(G)^m \xrightarrow{\nu(f)}  l^2(G)^n \xrightarrow{\nu(g)} l^2(G)^p$ is exact.
More details and explanations for the last two statements can be found
in \cite[Section 6.2]{Lueck(2002)}.

Consider the finitely generated $\caln(G)$-submodule $M \subseteq \caln(G)^n$.
Choose an $\caln(G)$-linear map $f \colon \caln(G)^m \to \caln(G)^n$
with image $M$.  The kernel of $\nu(f)$ is a closed $G$-invariant linear subspace of $l^2(G)^m$.
Hence there is an $\caln(G)$-map $p \colon \caln(G)^m \to \caln(G)^m$ such that
$\nu(p)$ is a $G$-equivariant projection, whose image is $\ker(\nu(f))$. 
Now $\nu(p) \circ \nu(p) = \nu(p)$ implies $p \circ p  = p$ and
 $\im(\nu(p)) = \ker(\nu(f))$ implies
$\im(p) = \ker(f)$.  Hence $\ker(f)$ is a direct summand in
$\caln(G)^m$ and $\im(f) = M$ is projective. 

The point is that in order to get the desired projection 
$p$ one passes to the interpretation by Hilbert spaces
and uses orthogonal projections there. We have enlarged the group ring
$\bbC G$ to the group von Neumann algebra $\caln(G)$, which 
does contain these orthogonal projections in contrast to $\bbC G$.
\em
\end{remark}


\subsection{Dimension Theory over the Group von Neumann Algebra}
\label{sec: Dimension Theory over the Group von Neumann Algebra}

An important feature of the group von Neumann algebra is its trace.

\begin{definition}[Von Neumann trace] \label{def: trace of the group von Neumann algebra}%
The \emph{von Neumann trace}%
\index{trace!von Neumann trace}
\index{von Neumann trace}
on $\caln(G)$ is defined by
$$\tr_{\caln(G)}\colon%
\indexnotation{tr_{caln(G)}}
 \caln(G) \to \bbC, \quad
f \mapsto \scal{f(e)}{e}_{l^2(G)},$$
where $e \in G \subseteq l^2(G)$ is the unit element.
\end{definition}

It enables us to define a dimension for finitely generated projective $\caln(G)$-modules.

\begin{definition}[Von Neumann dimension for finitely generated projective $\caln(G)$-modules] 
\label{def: Dimension of a finitely generated projective module}
Let $P$ be a finitely generated projective $\caln(G)$-module.
Choose a matrix $A = (a_{i,j}) \in M(n,n;\caln(G))$ with $A^2 = A$ such that the image of the $\caln(G)$-linear map
$r_A \colon \caln(G)^n \to \caln(G)^n$ given by right multiplication with $A$ is $\caln(G)$-isomorphic to $P$.
Define the \emph{von Neumann dimension}%
\index{von Neumann dimension!of a finitely generated projective $\caln(G)$-module}
\index{dimension!von Neumann dimension of a finitely generated projective $\caln(G)$-module}
of $P$ by
\begin{eqnarray*}
\dim_{\caln(G)}(P)
 & := & \sum_{i=1}^n \tr_{\caln(G)}(a_{i,i})
\quad \in [0,\infty).
\end{eqnarray*}
\end{definition}

We omit the standard proof that $\dim_{\caln(G)}(P)$ depends only on the isomorphism class of
$P$ but not on the choice of the matrix $A$.  Obviously 
$$\dim_{\caln(G)}(P \oplus Q) = \dim_{\caln(G)}(P) + \dim_{\caln(G)}(Q).$$
It is not hard to show
that $\dim_{\caln(G)}$ is \emph{faithful}, i.e. $\dim_{\caln(G)}(P) = 0 \Leftrightarrow P = 0$
holds for any finitely generated projective $\caln(G)$-module $P$.

Recall that the dual $M^{\ast}$ of a left or right 
$R$-module $M$ is the right or left  $R$-module
$\hom_{R}(M,R)$ respectively, where the $R$-multiplication is given by
$(fr)(x) = f(x)r$ or $(rf)(x) = rf(x)$ respectively for
$f \in M^{\ast}$, $x \in M$ and $r \in R$.

\begin{definition}[Closure of a submodule] \label{def: closure, K(M), P(M)}
Let $M$ be an $R$-submodule of $N$. Define
the \emph{closure}%
\index{closure of a submodule}
of $M$ in $N$
to be the $R$-submodule of $N$
$$\overline{M}%
\indexnotation{overline{M}}
~: = ~ \{ x \in N ~ \mid ~ f(x) = 0 \mbox{ for all }
f \in N^{\ast}  \mbox{ with } M \subseteq \ker(f)\}.$$
For an $R$-module $M$
define the $R$-submodule $\bfT M$%
\indexnotation{bfT M} and
the quotient $R$-module $\bfP M$%
\indexnotation{bfP M} by
\begin{eqnarray*}
\bfT M & := &\{x \in M ~ \mid ~ f(x) = 0
\mbox{ for all } f \in M^{\ast}\};
\\
\bfP M  & := & M/\bfT M.
\end{eqnarray*}
\end{definition}

Notice that $\bfT M$ is the closure of the trivial submodule in $M$.
It can also be described as the kernel of the canonical
map $i(M) \colon  M \to (M^{\ast})^{\ast}$, which sends
$x \in M$ to the map
$M^{\ast} \to R, \quad f \mapsto f(x)$.
Notice that $\bfT\bfP M = 0$, $\bfP\bfP M = \bfP M$, $M^* = (\bfP M)^*$ and that
$\bfP M = 0$ is equivalent to
$M^{\ast} = 0$. 

The next result is the key ingredient in the definition of $L^2$-Betti numbers for $G$-spaces.
Its proof can be found in \cite[Theorem 0.6] {Lueck(1998a)}, 
\cite[Theorem 6.7]{Lueck(2002)}.

\begin{theorem}{\bf (Dimension function for arbitrary $\caln(G)$-modules).}
\index{Theorem!Dimension function for arbitrary $\caln(G)$-modules}
\label{the: prop. ext. dim.}

\begin{enumerate}

\item \label{the: prop. ext. dim.: closure direct summand}
If $K \subseteq M$ is a submodule of the finitely generated
$\caln(G)$-module $M$, then $M/\overline{K}$ is finitely generated
projective and $\overline{K}$ is a direct summand in $M$;

\item
\label{the: prop. ext. dim.: decomposition in PM oplus TM}
If $M$ is a finitely generated $\caln(G)$-module, then
$\bfP M$ is finitely generated projective, there is an exact sequence
$0 \to \caln(G)^n \to \caln(G)^n \to \bfT M \to 0$ and
$$M \cong \bfP M \oplus \bfT M;$$

\item \label{the: prop. ext. dim.: properties of dim}
There exists precisely one dimension function
$$\dim_{\caln(G)} \colon \{\caln(G)\text{-modules}\} ~ \to ~ 
[0,\infty]%
\indexnotation{[0,infty]}
 := \{r \in \bbR \mid r \ge 0\} \amalg \{\infty\}$$
which satisfies:

\begin{enumerate}

\item
\label{the: prop. ext. dim.: properties of dim: extension property}
Extension Property\par\noindent
If $M$ is a finitely generated projective $\caln(G)$-module, then
$\dim_{\caln(G)}(M)$ agrees with the expression introduced in 
Definition \ref{def: Dimension of a finitely generated projective module};

\item
\label{the: prop. ext. dim.: properties of dim: additivity}
Additivity \par\noindent
If $0 \to M_0 \to M_1 \to M_2 \to 0$ is an
exact sequence of $\caln(G)$-modules, then
$$\dim_{\caln(G)}(M_1) ~ = ~ \dim_{\caln(G)}(M_0) + \dim_{\caln(G)}(M_2),$$
where for $r,s \in [0,\infty]$ we define $r + s$  by the ordinary sum of two
real numbers if both $r$ and $s$ are not $\infty$, and by $\infty$ otherwise;

\item
\label{the: prop. ext. dim.: properties of dim: cofinality}
Cofinality \par\noindent

Let $\{M_i\mid i \in I\}$ be a cofinal system of submodules of $M$,
i.e. $M = \bigcup_{i \in I} M_i$ and for two indices $i$ and $j$
there is an index $k$ in $I$ satisfying $M_i,M_j \subseteq M_k$.
Then
$$\dim_{\caln(G)}(M) ~ = ~ \sup\{\dim_{\caln(G)}(M_i) \mid i \in I\};$$

\item
\label{the: prop. ext. dim.: properties of dim: continuity}
Continuity \par\noindent
If $K \subseteq M$ is a submodule of the finitely generated
$\caln(G)$-module $M$, then
$$\dim_{\caln(G)}(K) ~ = ~ \dim_{\caln(G)}(\overline{K}).$$
\end{enumerate}

\end{enumerate}
\end{theorem}

\begin{definition}[Von Neumann dimension for arbitrary $\caln(G)$-modules] 
\label{def: Von Neumann dimension for arbitrary caln(G)-modules}
In the sequel we mean for an (arbitrary) $\caln(G)$-module $M$ by $\dim_{\caln(G)}(M)$%
\indexnotation{dim_{caln(G)}(M)}
the value of the dimension function appearing in
Theorem~\ref{the: prop. ext. dim.} and call it the \emph{von Neumann dimension of $M$}.%
\index{von Neumann dimension!of an $\caln(G)$-module}
\index{dimension!von Neumann dimension of an $\caln(G)$-module}
\end{definition}

\begin{remark}[Uniqueness of the dimension function] 
\label{rem: only possible candidate for dim_{caln(G)}} \em
There is only one possible definition for the dimension function
appearing in Theorem~\ref{the: prop. ext. dim.}, namely one must have
\begin{multline*}
\dim_{\caln(G)}(M) ~ := ~
\sup\{\dim_{\caln(G)}(P) \mid P \subseteq M \text{ finitely generated }
\\
\text{projective submodule}\}  \hspace{3mm} \in [0,\infty].
\end{multline*}
Namely, consider the directed system of finitely generated $\caln(G)$-submodules 
$\{M_i \mid i \in I\}$ of $M$ which is directed by inclusion.
By Cofinality
$$\dim_{\caln(G)}(M) ~ = ~ \sup\{\dim_{\caln(G)}(M_i) \mid i \in I\}.$$
From Additivity and 
Theorem~\ref{the: prop. ext. dim.}
\ref{the: prop. ext. dim.: decomposition in PM oplus TM} we conclude
$$\dim_{\caln(G)}(M_i) ~ = ~ \dim_{\caln(G)}(\bfP M_i)$$
and that $\bfP M_i$ is finitely generated projective. This shows uniqueness of
$\dim_{\caln(G)}$.  The hard part in the proof of Theorem~\ref{the: prop. ext. dim.}
\ref{the: prop. ext. dim.: properties of dim} is to show that the 
definition above does have all the desired properties.

We also see what $\dim_{\caln(G)}(M) = 0$ means. It is equivalent to the condition
that $M$ contains no non-trivial projective $\caln(G)$-submodule, or,
equivalently, no non-trivial finitely generated projective $\caln(G)$-submodule.
\em
\end{remark}

\begin{example}[The von Neumann dimension for finite groups] \em \label{exa: dim for finite G}
If $G$ is finite, then $\caln(G) = \bbC G$ and 
$\tr_{\caln(G)}\left(\sum_{g \in G} \lambda_g \cdot g\right)$ 
is the coefficient $\lambda_e$ of the unit element
$e \in G$. For an $\caln(G)$-module $M$ its von Neumann dimension
$\dim_{\caln(G)}(V)$ is $\frac{1}{|G|}$-times the complex dimension
of the underlying complex vector space $M$. 
\em
\end{example}

The next example implies that $\dim_{\caln(G)}(P)$ 
for a finitely generated projective $\caln(G)$-module
can be any non-negative real number.

\begin{example}[The von Neumann dimension for $\bbZ^n$] 
\label{exa: possible dimension of Hilbert modules} \em
Consider $G = \bbZ^n$. Recall that $\caln(\bbZ^n) = L^{\infty}(T^n)$.
Under this identification we get for the von Neumann trace
$$\tr_{\caln(\bbZ^n)} \colon \caln(\bbZ^n) \to \bbC, \quad f \mapsto \int_{T^n} f d\mu,$$
where $\mu$ is the standard Lebesgue measure on $T^n$.

Let $X \subseteq T^n$ be any measurable set and $\chi_X \in L^{\infty}(T^n)$
be its characteristic function. Denote by
$M_{\chi_X}\colon  L^2(T^n) \to  L^2(T^n)$
the $\bbZ^n$-equivariant unitary projection given by
multiplication with $\chi_X$. Its image $P$ is a finitely generated projective 
$\caln(\bbZ^n)$-module, whose von Neumann dimension 
$\dim_{\caln(\bbZ^n)}(P)$ is the volume $\mu(X)$ of $X$.
\em
\end{example}

In view of the results above the following slogan makes sense.
\begin{slogan} \em \label{slo: caln(G) like Z}
The group von Neumann algebra $\caln(G)$ behaves like the ring of integers $\bbZ$ provided
one ignores the properties integral domain and Noetherian.
\em
\end{slogan}
Namely, Theorem~\ref{the: prop. ext. dim.}
\ref{the: prop. ext. dim.: decomposition in PM oplus TM}
corresponds to the statement that a finitely generated $\bbZ$-module $M$ 
decomposes into $M = M/\tors(M) \oplus \tors(M)$ and that there exists an exact sequence
of $\bbZ$-modules $0 \to \bbZ^n \to \bbZ^n \to \tors(M) \to 0$,
where $\tors(M)$ is the $\bbZ$-module consisting of 
torsion elements. One obtains the obvious analog
of  Theorem~\ref{the: prop. ext. dim.}
\ref{the: prop. ext. dim.: properties of dim} if one considers
$$ \{\bbZ \text{-modules}\} ~ \to ~ [0,\infty], 
\quad M ~ \mapsto \dim_{\bbQ}(\bbQ \otimes_{\bbZ} M).$$
One basic difference between the case $\bbZ$ and $\caln(G)$  is
that there exist projective $\caln(G)$-modules with 
finite dimension which are not finitely generated,
which is not true over $\bbZ$. For instance
take the direct sum $P = \bigoplus_{i=1}^{\infty} P_i $ of $\caln(\bbZ^n)$-modules $P_i$ 
appearing in Example \ref{exa: possible dimension of Hilbert modules}  with 
$\dim_{\caln(\bbZ^n)}(P_i) = 2^{-i}$. Then $P$ is 
projective but not finitely generated and satisfies
$\dim_{\caln(\bbZ^n)}(P) = 1$.

The proof of the following two results is given in \cite[Theorem 6.13
and Theorem 6.39]{Lueck(2002)}.

\begin{theorem}[Dimension and colimits]
\index{Theorem!Dimension and colimits}
\label{the: dim. and col.}
Let $\{M_i \mid i \in I\}$ be a directed system of $\caln(G)$-modules over the directed set
$I$. For $i \le j$ let $\phi_{i,j}\colon  M_i \to M_j$
be the associated morphism of $\caln(G)$-modules.
For $i \in I$ let $\psi_{i}\colon  M_i \to \colim_{i \in I} M_i$
be the canonical morphism of $\caln(G)$-modules. Then:
\begin{enumerate}

\item
\label{the: dim. and col.: 1}
We get for the dimension of the $\caln(G)$-module given by
the colimit
$$\dim_{\caln(G)}\left(\colim_{i \in I} M_i\right) ~ = ~
\sup\left\{\dim_{\caln(G)}(\im(\psi_i)) \mid i \in I\right\};$$

\item \label{the: dim. and col.: 2}
Suppose
for each $i \in I$ that there exists  $i_0 \in I$ with 
$i \le i_0$ such that $\dim_{\caln(G)}(\im(\phi_{i,i_0})) < \infty$
holds. Then
\begin{eqnarray*}
\lefteqn{\dim_{\caln(G)}\left(\colim_{i \in I} M_i\right)} & &
\\
&  = &
\sup\left\{\inf\left\{
\dim_{\caln(G)}(\im(\phi_{i,j}\colon  M_i \to M_j)) \mid j \in I, i \le j
\right\}\mid i \in I\right\}.
\end{eqnarray*}
\end{enumerate}
\end{theorem}

\begin{theorem}[Induction and dimension] \label{the: Induction and Dimension}
\index{Theorem!Induction and dimension}
Let $i \colon  H \to G$ be an injective group homomorphism. Then
\begin{enumerate}

\item \label{the: Induction and Dimension: flat}
Induction with $\caln(i) \colon \caln(H) \to \caln(G)$ is a faithfully flat%
\index{faithfully flat functor} functor $M \mapsto i_*M%
\indexnotation{i_*M}
 := \caln(G) \otimes_{\caln(i)} M$
from the category of $\caln(H)$-modules to the category
of $\caln(G)$-modules, i.e. a  sequence
of $\caln(H)$-modules
$M_0 \to M_1 \to M_2$ is exact
at $M_1$ if and only if the induced sequence of $\caln(G)$-modules
$i_{\ast}M_0 \to i_{\ast}M_1
\to i_{\ast}M_2$ is exact at $i_{\ast}M_1$;

\item \label{the: Induction and Dimension: dim}
For any $\caln(H)$-module $M$ we have:
$$\dim_{\caln(H)}(M) ~ = ~ \dim_{\caln(G)}(i_{\ast}M).$$

\end{enumerate}
\end{theorem}

\begin{example}[The von Neumann dimension and {$\bbC [\bbZ^n]$}-modules]
\label{exa: dimension for G = Z^n} \em
Consider the case $G = \bbZ^n$. Then $\bbC[\bbZ^n]$ is a commutative integral domain
and hence has a quotient field $\bbC[\bbZ^n]_{(0)}$. Let $\dim_{\bbC[\bbZ^n]_{(0)}}$ denote
the usual dimension for vector spaces over $\bbC[\bbZ^n]_{(0)}$. Let $M$ be a
$\bbC[\bbZ^n]$-module. Then
\begin{eqnarray}
\hspace{-7mm} \dim_{\caln(\bbZ^n)}\left(\caln(\bbZ^n) \otimes_{\bbC[\bbZ^n]}M\right)
& = & 
\dim_{\bbC[\bbZ^n]_{(0)}}\left(\bbC[\bbZ^n]_{(0)} \otimes_{\bbC[\bbZ^n]} M \right).
\label{dim_{N(Z^n)} and dim_{C Z^n_{(0)}}}
\end{eqnarray}

This follows from the following considerations.
Let $\{M_i \mid i \in I\}$ be the directed system of finitely generated submodules of $M$.
Then $M = \colim_{i \in I} M_i$. Since the tensor product has a right adjoint, it is
compatible with colimits. This implies together with 
Theorem \ref{the: dim. and col.}
\begin{eqnarray*}
\dim_{\caln(\bbZ^n)}\left(\caln(\bbZ^n) \otimes_{\bbC[\bbZ^n]}M\right)
& = & \sup\left\{
\dim_{\caln(\bbZ^n)}\left(\caln(\bbZ^n) \otimes_{\bbC[\bbZ^n]}M_i\right)
\right\};
\\
\dim_{\bbC[\bbZ^n]_{(0)}}\left(\bbC[\bbZ^n]_{(0)} \otimes_{\bbC[\bbZ^n]}M\right)
& = & \sup\left\{
\dim_{\bbC[\bbZ^n]_{(0)}}\left(\bbC[\bbZ^n]_{(0)} \otimes_{\bbC[\bbZ^n]}M_i\right)
\right\}.
\end{eqnarray*}
Hence it suffices to prove the claim for a finitely generated $\bbC[\bbZ^n]$-module $N$.
The case $n = 1$ is easy. Then $\bbC[\bbZ]$ is a principal integral domain and we can write
$$N ~ = ~ \bbC[\bbZ]^r ~ \oplus ~ \bigoplus_{i=1}^k \bbC[\bbZ]/(u_i)$$
for non-trivial elements $u_i \in \bbC[\bbZ]$ and some non-negative
integers $k$ and $r$. One easily checks
that there is an exact $\caln(\bbZ)$-sequence 
$$0 \to \caln(\bbZ) \xrightarrow{r_{u_i}} \caln(\bbZ) \to \caln(\bbZ) \otimes_{\bbC[\bbZ]} 
\bbC[\bbZ]/(u_i) \to 0$$
using the identification $\caln(\bbZ) = L^{\infty}(S^1)$
from Example~\ref{exa: group von neumann algebra of Z^n}  to show injectivity
of the map $r_{u_i}$ given by multiplication with $u_i$. This implies
$$
\dim_{\caln(\bbZ)}\left(\caln(\bbZ) \otimes_{\bbC[\bbZ]}N\right)
~ = ~ r ~ = ~ 
\dim_{\bbC[\bbZ]_{(0)}}\left(\bbC[\bbZ]_{(0)} \otimes_{\bbC[\bbZ]}N\right).
$$
In the general case $n \ge 1$ one knows that there  exists a finite free 
$\bbC [\bbZ^n]$-resolution  of $N$.
Now the claim follows from \cite[Lemma 1.34]{Lueck(2002)}.

This example is the commutative version of a general setup for arbitrary groups,
which will be discussed in 
Subsection~\ref{subsec: The Ring Theoretic Version of the Atiyah Conjecture}.
\em
\end{example}

A center-valued dimension function for finitely generated projective modules
will be introduced in Definition~\ref{def: dim^u}. It can be used
to classify finitely generated projective $\caln(G)$-modules
(see Theorem~\ref{the: K_0 of fin. vN-alg.}) and shows that
the representation theory of finite dimensional representations over a finite group
extends to infinite groups if one works with $\caln(G)$ (see 
Remark~\ref{Group von Neumann algebras and representation theory}).


\typeout{--------------------   Section 2 --------------------------}

\section{Definition and Basic Properties of $L^2$-Betti Numbers}
\label{sec: Definition and Basic Properties of L^2-Betti Numbers}

In this section we define $L^2$-Betti  numbers for
arbitrary $G$-spaces and study their basic properties.
Our general algebraic definition is very general and is very flexible.
This allows to apply standard techniques such as spectral sequences and 
Mayer-Vietoris arguments directly. The original analytic definition for free proper
smooth $G$-manifolds with $G$-invariant Riemannian metrics is due to Atiyah and will
be briefly discussed in Subsection~\ref{subsec: Comparison with Other Definitions}.


\subsection{The Definition of $L^2$-Betti Numbers}

\label{sec: The Definition of L^2-Betti Numbers}

\begin{definition}[$L^2$-Betti numbers of $G$-spaces] 
\label{def: L^2-Betti numbers of G-spaces}
Let $X$ be a (left) G-space. Equip $\caln(G)$
with the obvious $\caln(G)$-$\bbZ G$-bimodule structure.
The \emph{singular homology} $H_p^G(X;\caln(G))$%
\indexnotation{H_p^G(X;caln(G))} \emph{of $X$ with coefficients in} $\caln(G)$%
\index{singular homology with coefficients}
 is the homology of the
$\caln(G)$-chain complex $\caln(G) \otimes_{\bbZ G} C_*^{\sing}(X)$, where
$C_*^{\sing}(X)$ is the singular chain complex of $X$ with the induced
 $\bbZ G$-structure. Define the
\emph{$p$-th $L^2$-Betti number of $X$}%
\index{L2-Betti number@$L^2$-Betti number!for arbitrary $G$-spaces}
by
$$b_p^{(2)}(X;\caln(G))%
\indexnotation{b_p^{(2)}(X;caln(G)) for arbitrary X}
 ~ := ~ \dim_{\caln(G)}\left(H_p^G(X;\caln(G))\right)
\quad \in [0,\infty],$$
where $\dim_{\caln(G)}$ is the dimension function of 
Definition~\ref{def: Von Neumann dimension for arbitrary caln(G)-modules}.

If $G$ and its action on $X$ are clear from the context,
we often omit $\caln(G)$ in the notation above. For instance, for a connected $CW$-complex
$X$ we denote by $b_p^{(2)}(\widetilde{X})$ the 
$L^2$-Betti number $b_p^{(2)}(\widetilde{X};\caln(\pi_1(X)))$ of its universal
covering $\widetilde{X}$%
\indexnotation{widetilde{X}}
 with respect to the obvious $\pi_1(X)$-action.
\end{definition}

Notice that we have \emph{no} assumptions on the $G$-action or on the  topology on $X$, 
we do \emph{not} need to require
that the operation is free, proper, simplicial or
cocompact. Thus we can apply this definition to the 
\emph{classifying space for free proper $G$-actions}%
\index{classifying space!for free proper $G$-actions}
$EG$,%
\indexnotation{EG}
which is a free $G$-$CW$-complex which is contractible (after
forgetting the group action). Recall that $EG$ is unique up to
$G$-homotopy. Its quotient $BG = G\backslash EG$%
\indexnotation{BG}
is a connected $CW$-complex, which is up to homotopy uniquely
determined by the property that $\pi_n(BG) = \{1\}$ for $n \ge 2$ and
$\pi_1(BG) \cong G$ holds, and called \emph{classifying space of $G$}.%
\index{classifying space!of a group}
Moreover, $G \to EG \to BG$ is the universal $G$-principal bundle.%
\index{universal!G-principal bundle@$G$-principal bundle}

\begin{definition}[$L^2$-Betti numbers of groups] 
\label{def: L^2-Betti numbers of groups}
Define for any (discrete) group $G$ its \emph{$p$-th $L^2$-Betti number}%
\index{L2-Betti number@$L^2$-Betti number!of a group} by
$$b^{(2)}_p(G)%
\indexnotation{b_p^(2)(G) for arbitrary G}
 ~ := ~ b^{(2)}_p(EG,\caln(G)).$$
\end{definition}

\begin{remark}[Comparison with the approach by Cheeger and Gromov] 
\em \label{rem: Comparison with the approach to Cheeger and Gromov}
A detailed comparison of our approach with the one by Cheeger and
Gromov \cite[section 2]{Cheeger-Gromov(1986)} can be found in
\cite[Remark 6.76]{Lueck(2002)}. Cheeger and Gromov \cite[Section 2]{Cheeger-Gromov(1986)}
define $L^2$-cohomology and $L^2$-Betti numbers of a
$G$-space $X$ by considering the category
whose objects are $G$-maps $f\colon  Y \to X$ for
a simplicial complex $Y$ with cocompact free simplicial
$G$-action and then using inverse limits to
extend the classical notions for
finite free $G$-$CW$-complexes such as $Y$
to $X$. Their approach is technically more complicated 
because for instance they work with cohomology instead of homology and therefore
have to deal with inverse limits instead of directed limits.
Our approach is closer to standard notions, the only non-standard
part is the verification of the properties of the extended
dimension function (Theorem \ref{the: prop. ext. dim.}). 
\em
\end{remark}

\begin{remark}[$L^2$-Betti numbers for von Neumann algebras]
\label{rem: L^2-Betti numbers for von Neumann algebras}
\em
The algebraic approach to $L^2$-Betti numbers of groups as
$$b_p^{(2)}(G) = \dim_{\caln(G)}\left(\Tor_p^{\bbC G}(\bbC,\caln(G))\right)$$ 
based on the dimension 
function for arbitrary modules and homological algebra 
plays a role in the definition of $L^2$-Betti numbers for certain
von Neumann algebras by Connes-Shlyakhtenko \cite{Connes-Shlyakhtenko(2003)}.
The point of their construction is to introduce invariants which depend
on the group von Neumann algebra $\caln(G)$ only.
If one could show that their invariants applied to $\caln(G)$ agree with 
the $L^2$-Betti numbers of $G$, one would get a positive answer to the open problem,
whether the von Neumann algebras of two finitely generated free groups 
$F_1$ and $F_2$ are isomorphic
as von Neumann algebras if and only if the groups $F_1$ and $F_2$ are isomorphic.
\em
\end{remark}

\begin{definition}[$G$-$CW$-complex]
\label{def: G-CW-complex}
A \emph{$G$-$CW$-complex}%
\index{G-CW-complex@$G$-$CW$-complex}
$X$ is a $G$-space together with a $G$-invariant filtration
$$\emptyset = X_{-1} \subseteq X_0 \subseteq X_1 \subseteq \ldots \subseteq
X_n \subseteq \ldots \subseteq \bigcup_{n \ge 0} X_n = X$$
such that $X$ carries the colimit topology
with respect to this filtration
(i.e. a set $C \subseteq X$ is closed if and only if $C \cap X_n$ is closed
in $X_n$ for all $n \ge 0$) and $X_n$ is obtained
from $X_{n-1}$ for each $n \ge 0$ by attaching
equivariant $n$-dimensional cells, i.e. there exists
a $G$-pushout
$$\comsquare
{\coprod_{i \in I_n} G/H_i \times S^{n-1}}
{\coprod_{i \in I_n} q_i}
{X_{n-1}}
{}
{}
{\coprod_{i \in I_n} G/H_i \times D^{n}}
{\coprod_{i \in I_n} Q_i}
{X_n}$$
\end{definition}

The space $X_n$ is called the \emph{$n$-skeleton}%
\index{skeleton}
of $X$. A $G$-$CW$-complex $X$ is \emph{proper}%
\index{G-CW-complex@$G$-$CW$-complex!proper}
if and only if all its isotropy groups are finite.  A $G$-space
is called \emph{cocompact}%
\index{cocompact}
if $G\backslash X$ is compact. A $G$-$CW$-complex $X$ is \emph{finite}%
\index{G-CW-complex@$G$-$CW$-complex!finite}
if $X$ has only finitely many equivariant cells. A $G$-$CW$-complex
is finite if and only if it is cocompact. A $G$-$CW$-complex $X$ is
\emph{of finite type}%
\index{G-CW-complex@$G$-$CW$-complex!of finite type}
if each $n$-skeleton is finite. It is called
\emph{of dimension $\le n$}%
\index{G-CW-complex@$G$-$CW$-complex!of dimension $\le n$}
if $X = X_n$ and \emph{finite dimensional}%
\index{G-CW-complex@$G$-$CW$-complex!finite dimensional}
if it is of dimension $\le n$ for some integer $n$.
A free $G$-$CW$-complex $X$ is the same as
a regular covering $X \to Y$ of a $CW$-complex $Y$ with
$G$ as group of deck transformations.

Notice that Definition~\ref{def: G-CW-complex}
also makes sense in the case where $G$ is a topological group.
Every proper smooth cocompact $G$-manifold is a proper $G$-$CW$-complex
by means of an equivariant triangulation.

For a $G$-$CW$-complex one can use the cellular 
$\bbZ G$-chain complex instead of the singular chain
complex in the definition of $L^2$-Betti numbers 
by the next result. Its proof can be found in \cite[Lemma 4.2]{Lueck(1998a)}.
For more information about $G$-$CW$-complexes we refer for instance 
to \cite[Sections II.1 and II.2]{Dieck(1987)},
\cite[Sections 1 and 2]{Lueck(1989)}, \cite[Subsection 1.2.1]{Lueck(2002)}.

\begin{lemma} \label{lem: singular = cellular}
Let $X$ be a $G$-$CW$-complex. Let $C_*^c(X)$ be its cellular $\bbZ G$-chain complex.
Then there is a $\bbZ G$-chain homotopy equivalence
$C_*^{\sing}(X) \to C_*^c(X)$ and we get
$$b_p^{(2)}(X;\caln(G)) ~ = ~ \dim_{\caln(G)}\left(H_p\left(\caln(G) \otimes_{\bbZ G} C_*^c(X)\right)\right).$$
\end{lemma}

The definition of $b_p^{(2)}(X;\caln(G))$  and the above lemma extend in the obvious way to pairs
$(X,A)$.


\subsection{Basic Properties of $L^2$-Betti Numbers}
\label{sec: Basic Properties of L^2-Betti Numbers}

The basic properties of $L^2$-Betti numbers are summarized in the following
theorem. Its proof can be found in \cite[Theorem 1.35 and Theorem
6.54]{Lueck(2002)} except for 
assertion~\ref{the: properties of gen. b_p^{(2)}: Restriction with epimorphisms with finite kernel}
which follows from \cite[Lemma 13.45]{Lueck(2002)}.

\begin{theorem}[$L^2$-Betti numbers for arbitrary spaces]
\index{Theorem!L2-Betti numbers for arbitrary spaces@$L^2$-Betti numbers for arbitrary spaces}
\label{the: properties of gen. b_p^{(2)}}%
\ \\
\begin{enumerate}
\item
\label{the: properties of gen. b_p^{(2)}: homology invariance}
Homology invariance
\\[1mm]
We have for a $G$-map  $f\colon  X \to Y$:
\begin{enumerate}
\item
\label{the: properties of gen. b_p^{(2)}: homotopy invariance: n-connected}
Suppose for $n \ge 1$
that for each subgroup $H \subseteq G$ the induced
map $f^H\colon  X^H \to Y^H$ is $\bbC$-homologically $n$-connected, i.e.
the  map
$$H_p^{\sing}(f^H;\bbC)\colon  H_p^{\sing}(X^H;\bbC) \to H_p^{\sing}(Y^H;\bbC)$$ 
induced by $f^H$ on singular homology with complex coefficients
is bijective for $p < n$ and surjective for $p = n$.
Then 
\begin{eqnarray*}
b_p^{(2)}(X) & =   & b_p^{(2)}(Y)
\hspace{10mm} \mbox{ for } p < n;
\\
b_p^{(2)}(X) & \ge & b_p^{(2)}(Y)
\hspace{10mm} \mbox{ for } p = n;
\end{eqnarray*}

\item
\label{the: properties of gen. b_p^{(2)}: homotopy invariance: H_*-iso}
Suppose that for each subgroup $H \subseteq G$ the induced
map $f^H\colon  X^H \to Y^H$
is a $\bbC$-homology equivalence, i.e.
$H_p^{\sing}(f^H;\bbC)$ is bijective for $p \ge 0$.
Then 
$$b_p^{(2)}(X) ~ = ~ b_p^{(2)}(Y) \hspace{10mm} \mbox{ for } p \ge 0;$$
\end{enumerate}

\item \label{the: properties of gen. b_p^{(2)}: comparison with Borel}
Comparison with the Borel construction
\\[1mm]
Let $X$ be a $G$-$CW$-complex. Suppose that for all $x \in X$ the isotropy group $G_x$
is finite or satisfies $b_p^{(2)}(G_x) = 0$ for all $p \ge 0$. Then
\begin{eqnarray*}
b_p^{(2)}(X;\caln(G)) & = &
b_p^{(2)}(EG \times X;\caln(G)) \hspace{8mm}
\mbox{ for }p \ge 0,
\end{eqnarray*}
where $G$ acts diagonally on $EG \times X$;

\item \label{the: properties of gen. b_p^{(2)}: inv. und. non. equiv. hom. equi.}
Invariance under non-equivariant $\bbC$-homology equivalences\\[1mm]
Suppose that
$f\colon  X \to Y$ is a $G$-equivariant map of $G$-$CW$-complexes
such that the induced map $H_p^{\sing}(f;\bbC)$ on singular
homology with complex coefficients is bijective for all $p$.
Suppose that for all $x \in X$ the isotropy group $G_x$
is finite or satisfies $b_p^{(2)}(G_x) = 0$ for all $p \ge 0$, and
analogously for all $y \in Y$. Then we have for all $p \ge 0$
$$b_p^{(2)}(X;\caln(G)) ~ = ~ b_p^{(2)}(Y;\caln(G));$$

\item
\label{the: properties of gen. b_p^{(2)}: sees only the proper part}
Independence of equivariant cells with infinite isotropy
\\[1mm]
Let $X$ be a $G$-$CW$-complex. Let $X[\infty]$%
\indexnotation{X[infty]}
be the $G$-$CW$-subcomplex
consisting of those points whose isotropy subgroups are infinite. Then
we get for all $p \ge 0$
$$b_p^{(2)}(X;\caln(G)) ~ = ~ b_p^{(2)}(X,X[\infty];\caln(G));$$

\item \label{the: properties of gen. b_p^{(2)}: Kuenneth formula}
K\"unneth formula
\\[1mm]
Let $X$ be a $G$-space and $Y$ be an $H$-space.
Then $X \times Y$ is a $G \times
H$-space and we get for all $n \ge 0$
\begin{eqnarray*}
b_n^{(2)}(X \times Y)
& = &
\sum_{p + q = n} b_p^{(2)}(X) \cdot b_q^{(2)}(Y),
\end{eqnarray*}
where we use the convention that $0 \cdot \infty = 0$, $r \cdot \infty
= \infty$ for $r \in (0,\infty]$ and $r + \infty = \infty$ for
$r \in [0,\infty]$;

\item \label{the: properties of gen. b_p^{(2)}: induction}
Induction
\\[1mm]
Let $i\colon  H \to G$ be an inclusion of groups
and let $X$ be an $H$-space. Let
$\caln(i)\colon  \caln(H) \to \caln(G)$ be the induced ring homomorphism
(see \eqref{caln(i)}). Then:
\begin{eqnarray*}
H_p^G(G\times_H X;\caln(G)) & = &
\caln(G) \otimes_{\caln(i)} H_p^H(X;\caln(H));
\\
b_p^{(2)}(G\times_H X;\caln(G)) & = &
b_p^{(2)}(X;\caln(H));
\end{eqnarray*}

\item \label{the: properties of gen. b_p^{(2)}: restriction}
Restriction to subgroups of finite index
\\[1mm]
Let $H \subseteq G$ be a subgroup of finite index $[G:H]$. 
Let $X$ be a $G$-space and let $\res_G^H X$ be the $H$-space
obtained from $X$ by restriction. Then
$$b_p^{(2)}(\res_G^H X;\caln(H)) ~ = ~
[G:H] \cdot b_p^{(2)}(X;\caln(G));$$

\item \label{the: properties of gen. b_p^{(2)}: Restriction with epimorphisms with finite kernel}
Restriction with epimorphisms with finite kernel\\[1mm]
Let $p \colon G \to Q$ be an epimorphism of groups with finite kernel $K$.
Let $X$ be a $Q$-space. Let $p^*X$ be the $G$-space obtained from $X$ using $p$. Then
\begin{eqnarray*}
b_p^{(2)}(p^*X;\caln(G)) & = &
\frac{1}{|K|} \cdot b_p^{(2)}(X;\caln(Q));
\end{eqnarray*}

\item \label{the: properties of gen. b_p^{(2)}: b_0^{(2)}}
Zero-th homology and $L^2$-Betti number
\\[1mm]
Let $X$ be a path-connected $G$-space. Then:
\begin{enumerate}
\item \label{the: properties of gen. b_p^{(2)}: b_0^{(2)}: H_0^{(2)}}
There is an $\caln(G)$-isomorphism
$H_0^G(X;\caln(G)) \xrightarrow{\cong}
~ \caln(G)\otimes_{\bbC G} \bbC$;

\item \label{the: properties of gen. b_p^{(2)}: b_0^{(2)}: value}
$b^{(2)}_0(X;\caln(G)) = |G|^{-1}$, where
$|G|^{-1}$ is defined to be zero if the order $|G|$ of
$G$ is infinite;

\end{enumerate}

\item \label{the: properties of gen. b_p^{(2)}: Euler-Poincar'e formula}
Euler-Poincar\'e formula
\\[1mm]
Let $X$ be a free finite $G$-$CW$-complex.
Let $\chi(G\backslash X)$ be the Euler characteristic of the finite
$CW$-complex
$G\backslash X$, i.e.
\begin{eqnarray*}
\chi(G\backslash X)
& := &
 \sum_{p \ge 0} (-1)^p \cdot |I_p(G\backslash X)|
\hspace{10mm} \in \bbZ ,
\end{eqnarray*}
where
$|I_p(G\backslash X)|$ is the number of $p$-cells of $G\backslash X$. Then
\begin{eqnarray*}
\chi(G\backslash X)
& = &
\sum_{p\ge 0} (-1)^p \cdot b^{(2)}_p(X);
\end{eqnarray*}

\item \label{the: properties of gen. b_p^{(2)}:
Morse inequalities}
Morse inequalities
\\[1mm]
Let $X$ be a free $G$-$CW$-complex of finite type.  Then
we get for $n \ge 0$
\begin{eqnarray*}
\sum_{p=0}^n (-1)^{n-p} \cdot b_p^{(2)}(X)
& \le &
\sum_{p=0}^n (-1)^{n-p} \cdot |I_p(G\backslash X)|;
\end{eqnarray*}

\item \label{the: properties of gen. b_p^{(2)}: Poincar'e duality}
Poincar\'e duality
\\[1mm]
Let $M$ be a cocompact free proper $G$-manifold of dimension $n$
which is orientable. Then
\begin{eqnarray*}
b^{(2)}_p(M)
& = &
b^{(2)}_{n-p}(M,\partial M);
\end{eqnarray*}

\item \label{the: properties of gen. b_p^{(2)}: wedges}
Wedges\\[1mm]
Let $X_1$, $X_2$, \ldots, $X_r$ be
connected (pointed) $CW$-complexes of finite type
and $X = \bigvee_{i=1}^r X_i$ be their wedge.
Then
\begin{eqnarray*}
b_1^{(2)}(\widetilde{X}) - b_0^{(2)}(\widetilde{X})
& = &
r -1 + \sum_{j=1}^r
\left( b_1^{(2)}(\widetilde{X_j}) - b_0^{(2)}(\widetilde{X_j})\right);
\\
b_p^{(2)}(\widetilde{X})
& = &
\sum_{j=1}^r b_p^{(2)}(\widetilde{X_j})
\hspace{10mm} \mbox{ for } 2 \le p;
\end{eqnarray*}

\item \label{the: properties of gen. b_p^{(2)}: connected sums}
Connected sums\\[1mm]
Let $M_1$, $M_2$, \ldots, $M_r$ be
compact connected $m$-dimensional manifolds for $m \ge 3$.
Let $M$ be their connected sum $M_1 \# \ldots \# M_r$. Then
\begin{eqnarray*}
b_1^{(2)}(\widetilde{M}) - b_0^{(2)}(\widetilde{M})
& = &
r -1 + \sum_{j=1}^r
\left( b_1^{(2)}(\widetilde{M_j}) - b_0^{(2)}(\widetilde{M_j})\right);
\\
b_p^{(2)}(\widetilde{M})
& = &
\sum_{j=1}^r b_p^{(2)}(\widetilde{M_j})
\hspace{10mm} \mbox{ for } 2 \le p \le m-2.
\end{eqnarray*}

\end{enumerate}
\end{theorem}

\begin{example} \label{exa: b_p^{(2)} cdot || = b_p for finite G}
\em If $G$ is finite, then $b_p^{(2)}(X;\caln(G))$ reduces to the classical Betti number
$b_p(X)$ multiplied with the factor $|G|^{-1}$. 
\end{example}

\begin{remark}[Reading off $L^2$-Betti numbers from $H_p(X;\bbC)$]
\label{rem: comments on properties of L^2-Betti numbers} \em
If $f \colon X \to Y$ is a $G$-map of free $G$-$CW$-complexes which induces
isomorphisms $H_p^{\sing}(f;\bbC)$ for all $p \ge 0$, then
Theorem~\ref{the: properties of gen. b_p^{(2)}}
\ref{the: properties of gen. b_p^{(2)}: homology invariance} implies
$$b_p^{(2)}(X;\caln(G)) ~ = ~ b_p^{(2)}(Y;\caln(G)).$$

This does not necessarily mean that one can read off
$b_p^{(2)}(X;\caln(G))$ from the singular homology $H_p(X;\bbC)$
regarded as a $\bbC G$-module in general. In general there is for a free $G$-$CW$-complex $X$
a spectral sequence
converging to $H_{p+q}^G(X;\caln(G))$, whose $E^2$-term is
$$E^2_{p,q}  ~ = ~ \Tor^{\bbC G}_p(H_q(X;\bbC),\caln(G)).$$
There is no reason why the equality of the dimension of the $E^2$-term for two free
$G$-$CW$-complexes $X$ and $Y$ implies that the dimension of $H_{p+q}^G(X;\caln(G))$
and $H_{p+q}^G(Y;\caln(G))$ agree. 
However, this is the case if the spectral sequence collapses
from the dimension point of view. For instance, if we make the assumption
$\dim_{\caln(G)}\left(\Tor^{\bbC G}_p(M,\caln(G))\right) ~ = ~ 0$ for all $\bbC G$-modules
$M$ and $p \ge 2$, Additivity and Cofinality of $\dim_{\caln(G)}$ 
(see Theorem \ref{the: prop. ext. dim.}) imply
\begin{multline*}
b_p^{(2)}(X;\caln(G)) ~ = ~ 
\\
\dim_{\caln(G)}\left(\caln(G) \otimes_{\bbC G}   H_p(X;\bbC)\right) + 
\dim_{\caln(G)}\left(\Tor^{\bbC G}_1(H_{p-1}(X;\bbC),\caln(G))\right).
\end{multline*}
The assumption above is satisfied if $G$ is amenable 
(see Theorem~\ref{the: dimension of higher Tor-s vanish in the amenable case})
or $G$ has cohomological dimension $\le 1$ over $\bbC$, for instance, if $G$ is virtually
free.
\em
\end{remark}

\begin{remark}[$L^2$-Betti numbers ignore infinite isotropy]  \label{ingnoring infinite isotropy}
\em 
Theorem~\ref{the: properties of gen. b_p^{(2)}}
\ref{the: properties of gen. b_p^{(2)}: sees only the proper part}
says that the $L^2$-Betti numbers do not see the part of
a $G$-space $X$ whose isotropy groups are infinite. In particular
$b_p^{(2)}(X;\caln(G)) = 0$ if $X$ is a $G$-$CW$-complex whose isotropy groups are all
infinite. This follows from the fact that for a subgroup $H \subseteq G$
$$\dim_{\caln(G)}\left(\caln(G) \otimes_{\bbC G} \bbC [G/H]\right) 
~ = ~
\left\{\begin{array}{lll} \frac{1}{|H|} & & \text {if } |H| < \infty;\\
0 & & \text {if } |H| = \infty.
\end{array}
\right.
$$
\em
\end{remark}

\begin{remark}[$L^2$-Betti numbers often vanish] 
\label{rem: L^2-Betti numbers often vanish}
\em
An important phenomenon is that
the $L^2$-Betti numbers of universal coverings of spaces and of 
groups tend to vanish more often than the classical
Betti numbers. This allows to draw interesting conclusions as we will see later.
\em
\end{remark}

\subsection{Comparison with Other Definitions}

\label{subsec: Comparison with Other Definitions}

In this subsection we give a short overview of the previous definitions
of $L^2$-Betti numbers. Originally they were defined in terms of heat kernels.
Their analytic aspects are important, but we  will only focus on their
algebraic aspects in this survey article. So a reader may skip the brief explanations below.

The notion of $L^2$-Betti numbers is due to Atiyah \cite{Atiyah(1976)}.
He defined for a smooth Riemannian manifold
with a free proper cocompact $G$-action by isometries
its \emph{analytic $p$-th $L^2$-Betti number}%
\index{L2-Betti number@$L^2$-Betti number!analytic}
by the following expression in terms of the 
\emph{heat kernel}
$e^{-t\Delta_p}(x,y)$ of the $p$-th Laplacian $\Delta_p$
\begin{eqnarray}
b_p^{(2)}(M)
& = &
\lim_{t \to \infty} \int_{\calf} \tr_{\bbC}(e^{-t\Delta_p}(x,x)) ~ dvol_x ,
\label{analytic L^2-Betti number and heat kernel}
\end{eqnarray}
where $\calf$ is a fundamental domain for the $G$-action and 
$\tr_{\bbC}$ denotes the trace of an endomorphism of a finite-dimensional vector space.
The $L^2$-Betti numbers are invariants of the large times asymptotic of the heat
kernel.

A \emph{finitely generated Hilbert $\caln(G)$-module}%
\index{Hilbert $\caln(G)$-module!finitely generated}
is a Hilbert space $V$ together with a
linear $G$-action by isometries such that there exists a linear isometric
$G$-embedding into $l^2(G)^n$ for some $n \ge 0$. 
One can assign to it its \emph{von Neumann dimension}%
\index{von Neumann dimension!of a finitely generated Hilbert $\caln(G)$-module}
by
$$\dim_{\caln(G)}(V) ~ := ~ \tr_{\caln(G)}(A) \quad\in [0,\infty),$$
where $A$ is any idempotent matrix $A \in M(n,n;\caln(G))$ such that the image of 
the $G$-equivariant operator $l^2(G)^n \to l^2(G)^n$ induced by $A$ is
isometrically linearly $G$-isomorphic to $V$.

The expression in \eqref{analytic L^2-Betti number and heat kernel} 
can be interpreted as the von Neumann
dimension of the \emph{space $\calh^p_{(2)}(M)$ of square-integrable harmonic $p$-forms} on
$M$, which is a finitely generated Hilbert $\caln(G)$-module
(see \cite[Proposition 4.16 on page 63]{Atiyah(1976)})
\begin{eqnarray}
\lim_{t \to \infty} \int_{\calf} \tr_{\bbC}(e^{-t\Delta_p}(x,x)) ~ dvol_x 
& = & \dim_{\caln(G)}\left(\calh^p_{(2)}(M)\right).
\label{b_p^{(2)} = dim(calh^p)}
\end{eqnarray}

Given a cocompact free $G$-$CW$-complex $X$,  one obtains a chain complex
of finitely generated Hilbert $\caln(G)$-modules 
$C_*^{(2)}(X)%
\indexnotation{C_*^(2)(X)}
 ~ := ~ C_*^c(X) \otimes_{\bbZ G} l^2(G)$.
Its \emph{reduced $p$-th $L^2$-homology} is the finitely generated Hilbert $\caln(G)$-module
\begin{eqnarray}
H_p^{(2)}(X;l^2(G))%
\indexnotation{H_p^{(2)}(X;l^2(G))}
 & = & \ker(c_p^{(2)})/\overline{\im(c^{(2)}_{p+1}))}.
\label{reduced homology H_p^{(2)}(X;l^2(G))}
\end{eqnarray}
Notice that we divide out the closure of the image of the $(p+1)$-th differential
$c_{p+1}^{(2)}$ of $C_*^{(2)}(X)$ in order to ensure that we obtain a Hilbert space.
Then by a result of Dodziuk \cite{Dodziuk(1977)} there is an isometric bijective $G$-operator
\begin{eqnarray}
\calh^p_{(2)}(M) & \xrightarrow{\cong} & H_p^{(2)}(K;l^2(G)),
\label{Dodziuk's result}
\end{eqnarray}
where $K$ is an equivariant triangulation of $M$. Finally one can show
\cite[Theorem 6.1]{Lueck(1997a)}
\begin{eqnarray}
b_p^{(2)}(K;\caln(G)) & = & \dim_{\caln(G)}\left(H_p^{(2)}(K;l^2(G))\right),
\label{our def = dodziuk's def}
\end{eqnarray}
where $b_p^{(2)}(K;\caln(G))$ is the $p$-th $L^2$-Betti number in the sense of
Definition \ref{def: L^2-Betti numbers of G-spaces}.

All in all we see that our Definition~\ref{def: L^2-Betti numbers of G-spaces} 
of $L^2$-Betti numbers for arbitrary $G$-spaces extends the heat kernel definition
of \eqref{analytic L^2-Betti number and heat kernel} for  
smooth Riemannian manifolds with a free proper cocompact $G$-action by isometries.
More details of all these  definitions and of their identifications
can be found in \cite[Chapter 1]{Lueck(2002)}.


\subsection{$L^2$-Euler Characteristic}
\label{subsec: L2-Euler Characteristic}

In this section we introduce the notion of
$L^2$-Euler characteristic.

If $X$ is a $G$-$CW$-complex, denote by
$I(X)$%
\indexnotation{I(X)}
the set of its equivariant cells.
For a cell $c \in I(X)$ let $(G_c)$ be
the conjugacy class of subgroups of $G$ given by its orbit type
and let $\dim(c)$ be its dimension.
Denote by $|G_c|^{-1}$ the inverse of the order of any representative
of $(G_c)$, where $|G_c|^{-1}$ is to be understood to be zero
if the order is infinite.

\begin{definition}[$L^2$-Euler characteristic]
\label{def: i(X,caln(G)), h^{(2)}(X,caln(G)) and chi^{(2)}(X;caln(G))}
Let $G$ be a group and let $X$ be a $G$-space. Define
$$\begin{array}{lll}
h^{(2)}(X;\caln(G))%
\indexnotation{h(2)(X;caln(G))}
& := & \sum_{p \ge 0} b_p^{(2)}(X;\caln(G)) ~ \in [0,\infty];
\\
\chi^{(2)}(X;\caln(G))%
\indexnotation{chi^{(2)}(X;caln(G))}
& := & \sum_{p \ge 0} (-1)^p \cdot b_p^{(2)}(X;\caln(G))
~ \in \bbR, \text{ if } h^{(2)}(X;\caln(G)) < \infty;
\\
m(X;G)%
\indexnotation{m(X;caln(G))}
& := & \sum_{c \in I(X)} |G_c|^{-1} ~ \in [0,\infty],
\quad\text{ if } X
\text{ is a } G\text{-}CW\text{-complex};
\\
h^{(2)}(G)
\indexnotation{h2G}
& := & h^{(2)}(EG;\caln(G)) ~ \in [0,\infty];
\\
\chi^{(2)}(G)%
\indexnotation{chi^{(2)}(G)} 
& := & \chi^{(2)}(EG;\caln(G))
~ \in \bbR, \quad\text{ if } h^{(2)}(G) < \infty.
\end{array}$$
We call $\chi^{(2)}(X;\caln(G))$ and $\chi^{(2)}(G)$
the \emph{$L^2$-Euler characteristic}%
\index{Euler characteristic!$L^2$-Euler characteristic}
\index{L2-Euler characteristic@$L^2$-Euler characteristic}
of  $X$ and $G$.
\end{definition}

The condition $h^{(2)}(X;\caln(G)) < \infty$ ensures
that the sum which appears in the definition of
$\chi^{(2)}(X;\caln(G))$ converges absolutely and that the following
results are true. The reader should compare the next theorem
with \cite[Theorem 0.3 on page 191]{Cheeger-Gromov(1986)}.
It essentially follows from 
Theorem~\ref{the: properties of gen. b_p^{(2)}}.
Details of its proof can be found in \cite[Theorem 6.80]{Lueck(2002)}.

\begin{theorem}[$L^2$-Euler characteristic]
\index{Theorem!L2-Euler characteristic@$L^2$-Euler characteristic}
\label{the: basic properties of the l2-Euler characteristic}
~\
\begin{enumerate}
\item
\label{the: basic properties of the l2-Euler characteristic: m and h}
Generalized Euler-Poincar\'e formula
\\[1mm]
Let $X$ be a $G$-$CW$-complex with $m(X;G) < \infty$. Then
\begin{eqnarray*}
h^{(2)}(X;\caln(G)) & < & \infty;
\\
\sum_{c \in I(X)} (-1)^{\dim(c)} \cdot |G_c|^{-1} & = &
\chi^{(2)}(X;\caln(G));
\end{eqnarray*}

\item
\label{the: basic properties of the l2-Euler characteristic: sum formula}
Sum formula
\\[1mm]
Consider the following $G$-pushout
$$\comsquare{X_0}{i_1}{X_1}{i_2}{j_1}{X_2}{j_2}{X}$$
such that $i_1$ is a $G$-cofibration.
Suppose that $h^{(2)}(X_i;\caln(G)) < \infty$
for $i=0,1,2$. Then
\begin{eqnarray*}
h^{(2)}(X;\caln(G)) & < & \infty;
\\
\chi^{(2)}(X;\caln(G)) & = &
\chi^{(2)}(X_1;\caln(G)) + \chi^{(2)}(X_2;\caln(G))-
\chi^{(2)}(X_0;\caln(G));
\end{eqnarray*}

\item
\label{the: basic properties of the l2-Euler characteristic: EG times X}
Comparison with the Borel construction
\\[1mm]
Let $X$ be a $G$-$CW$-complex. If for all $c \in I(X)$
the group $G_c$ is finite or $b_p^{(2)}(G_c) = 0$ for all $p \ge 0$, then
\begin{eqnarray*}
b_p^{(2)}(X;\caln(G)) & = &
b_p^{(2)}(EG \times X;\caln(G)) \hspace{8mm}
\text{ for }p \ge 0;
\\
h^{(2)}(X;\caln(G)) & = &
h^{(2)}(EG \times X;\caln(G));
\\
\chi^{(2)}(X;\caln(G)) & = &
\chi^{(2)}(EG \times X;\caln(G)), 
\text{if } h^{(2)}(X;\caln(G)) < \infty;
\\
\sum_{c \in I(X)} (-1)^{\dim(c)} \cdot |G_c|^{-1}& = &
\chi^{(2)}(EG \times X;\caln(G)), \hspace{1mm}
\text{if } m(X;G) < \infty;
\end{eqnarray*}

\item \label{the: basic properties of the l2-Euler characteristic: X to Y}
Invariance under non-equivariant $\bbC$-homology equivalences
\\[1mm]
Suppose that
$f\colon  X \to Y$ is a $G$-equivariant map of $G$-$CW$-complexes
with $m(X;G) < \infty$ and $m(Y;G) < \infty$,
such that the induced map $H_p(f;\bbC)$ on
homology with complex coefficients is bijective for all $p \ge 0$.
Suppose that for all  $c \in I(X)$
the group $G_c$ is finite or $b_p^{(2)}(G_c) = 0$ for all $p\ge 0$, and
analogously for all $d \in I(Y)$. Then
\begin{eqnarray*}
\chi^{(2)}(X;\caln(G))
& = &
\sum_{c \in I(X)} (-1)^{\dim(c)} \cdot |G_c|^{-1}
\\
& = &
\sum_{d \in I(Y)} (-1)^{\dim(d)} \cdot |G_d|^{-1}
\\
& = &
\chi^{(2)}(Y;\caln(G));
\end{eqnarray*}

\item
\label{the: basic properties of the l2-Euler characteristic: Kuenneth}
K\"unneth formula
\\[1mm]
Let $X$ be a $G$-$CW$-complex and $Y$ be an $H$-$CW$-complex.
Then we get for the $G \times H$-$CW$-complex $X \times Y$
\begin{eqnarray*}
m(X\times Y;G \times H)
& = &
m(X;G) \cdot m(Y;H);
\\
h^{(2)}(X\times Y;\caln(G \times H))
& = &
h^{(2)}(X;\caln(G)) \cdot h^{(2)}(Y;\caln(H));
\\
\chi^{(2)}(X\times Y;\caln(G \times H)) & = &
\chi^{(2)}(X;\caln(G)) \cdot \chi^{(2)}(Y;\caln(H)),
\\
& &
\hspace{5mm} \text{if } h^{(2)}(X;\caln(G)), h^{(2)}(Y;\caln(H)) < \infty,
\end{eqnarray*}
where we use the convention that $0 \cdot \infty = 0$ and
$r \cdot \infty= \infty$ for $r \in (0,\infty]$;

\item
\label{the: basic properties of the l2-Euler characteristic: induction}
Induction
\\[1mm]
Let $H\subseteq G$ be a subgroup and let $X$ be an $H$-space. Then
\begin{eqnarray*}
m(G\times_H X;G)
& = &
m(X;H);
\\
h^{(2)}(G\times_H X;\caln(G))
& = &
h^{(2)}(X;\caln(H));
\\
\chi^{(2)}(G\times_H X;\caln(G))
& = &
\chi^{(2)}(X;\caln(H)),
\hspace{10mm} \text{if } h^{(2)}(X;\caln(H)) < \infty;
\end{eqnarray*}

\item
\label{the: basic properties of the l2-Euler characteristic: restriction}
Restriction to subgroups of finite index
\\[1mm]
Let $H \subseteq G$ be a subgroup of finite index $[G:H]$.
Let $X$ be a $G$-space and let $\res_G^H X$ be the $H$-space
obtained from $X$ by restriction. Then
\begin{eqnarray*}
m(\res_G^H X;H) & = & [G:H] \cdot m(X;G);
\\
h^{(2)}(\res_G^H X;\caln(H)) & = & [G:H] \cdot h^{(2)}(X;\caln(G));
\\
\chi^{(2)}(\res_G^H X;\caln(H)) & = & [G:H] \cdot \chi^{(2)}(X;\caln(G)),
\hspace{1mm} \text{if } h^{(2)}(X;\caln(G)) < \infty,
\end{eqnarray*}
where $[G:H] \cdot \infty$ is understood to be $\infty$;

\item \label{the: basic properties of the l2-Euler characteristic: Restr. with epi. with fin. ker.}
Restriction with epimorphisms with finite kernel\\[1mm]
Let $p \colon G \to Q$ be an epimorphism of groups with finite kernel $K$.
Let $X$ be a $Q$-space. Let $p^*X$ be the $G$-space obtained from $X$ using $p$. Then
\begin{eqnarray*}
m(p^*X;G) & = & |K|^{-1} \cdot m(X;Q);
\\
h^{(2)}(p^*X;\caln(G)) & = & |K|^{-1}  \cdot h^{(2)}(X;\caln(Q));
\\
\chi^{(2)}(p^*X;\caln(G)) & = & |K|^{-1}  \cdot \chi^{(2)}(X;\caln(Q)),
\hspace{4mm} \text{if } h^{(2)}(X;\caln(Q)) < \infty.
\end{eqnarray*}

\end{enumerate}
\end{theorem}

\begin{remark}[$L^2$-Euler characteristic and virtual Euler characteristic]
 \label{rem: virtual and L^2-Euler characteristics} \em
The $L^2$-Euler characteristic generalizes the notion of the virtual
Euler characteristic.
Let $X$ be a $CW$-complex which is \emph{virtually homotopy finite}, i.e.
there is a $d$-sheeted covering $p\colon  \overline{X} \to X$
for some positive integer $d$ such that
$\overline{X}$ is homotopy equivalent to a finite $CW$-complex.
Define the \emph{virtual Euler characteristic}%
\index{Euler characteristic!virtual Euler characteristic}
following Wall \cite{Wall(1961)}
$$
\chi_{\virt}(X)
\indexnotation{chivirtuell} 
~ := ~ \frac{\chi(\overline{X})}{d}.
$$
One easily checks that this is independent of the choice of
$p\colon  \overline{X} \to X$ since the classical Euler characteristic is
multiplicative under finite coverings. Moreover, we conclude from
Theorem~\ref{the: basic properties of the l2-Euler characteristic}
\ref{the: basic properties of the l2-Euler characteristic: m and h}
and
\ref{the: basic properties of the l2-Euler characteristic: restriction}
that for virtually homotopy finite $X$
\begin{eqnarray*}
m(\widetilde{X};\pi_1(X)) & < & \infty;
\\
\chi^{(2)}(\widetilde{X};\caln(\pi_1(X))) & = &
\chi_{\virt}(X).
\end{eqnarray*}
\em
\end{remark}

\begin{remark}[$L^2$-Euler characteristic and orbifold Euler characteristic]
\label{rem: virtual and orbifold characteristics} \em
If $X$ is a finite $G$-$CW$-complex, then 
$\sum_{c \in I(X)} (-1)^{\dim(c)} \cdot |G_c|^{-1}$ is also called 
\emph{orbifold Euler characteristic}%
\index{Euler characteristic!orbifold Euler characteristic}
and agrees with the $L^2$-Euler characteristic by 
Theorem~\ref{the: basic properties of the l2-Euler characteristic}
\ref{the: basic properties of the l2-Euler characteristic: m and h}.
\em
\end{remark}
 

\typeout{--------------------   Section 3 --------------------------}

\section{Computations of $L^2$-Betti Numbers}
\label{sec: Computations of L^2-Betti Numbers}

In this section we state some cases where the $L^2$-Betti numbers
$b_p^{(2)}(\widetilde{X})$ for certain compact manifolds or finite $CW$-complexes
$X$ can explicitly be computed. These computations give evidence for certain conjectures
such as the Atiyah Conjecture~\ref{con: Atiyah Conjecture} for $(G,d,\bbQ)$ and the Singer
Conjecture~\ref{con: Singer Conjecture} which we will discuss later. Sometimes we will also
make a few comments on their proofs in order to give some insight into the methods. 
Besides analytic methods, which will not be discussed, standard techniques
from topology and algebra such as spectral sequences and  Mayer-Vietoris sequences
will play a role.
With our algebraic setup and the nice properties of the dimension function
such as Additivity and Cofinality these tools are directly available, whereas
in the original settings, which we have briefly discussed in
Subsection~\ref{subsec: Comparison with Other Definitions}, these methods do not apply directly
and, if at all, only after some considerable technical efforts.


\subsection{Abelian Groups}
\label{subsec: abelian groups}

 Let $X$ be a $\bbZ^n$-space. Then we get from 
\eqref{dim_{N(Z^n)} and dim_{C Z^n_{(0)}}}
\begin{eqnarray}
b_p^{(2)}(X;\caln(\bbZ^n)) 
& = & \dim_{\bbC[\bbZ^n]_{(0)}}
\left(\bbC[\bbZ^n]_{(0)} \otimes_{\bbC[\bbZ^n]} H_p^{\sing}(X;\bbC)\right).
\label{L^2-Betti numbers of Z^n-spaces}
\end{eqnarray}
Notice that $b_p^{(2)}(X;\caln(\bbZ^n))$ is always an integer or $\infty$.


\subsection{Finite Coverings}
\label{subsec: Finite Coverings}

Let $p \colon X \to Y$ be a finite covering with $d$-sheets. Then we conclude from 
Theorem~\ref{the: properties of gen. b_p^{(2)}}
\ref{the: properties of gen. b_p^{(2)}: restriction}
\begin{eqnarray}
b_p^{(2)}(\widetilde{X}) & = & d\cdot b_p^{(2)}(\widetilde{Y}).
\label{L^2-Betti numbers of finite coverings}
\end{eqnarray}
This implies for every connected $CW$-complex $X$ which admits a selfcovering
$X \to X$ with $d$-sheets for $d \ge 2$ that $b_p^{(2)}(\widetilde{X)} = 0$ for all 
$p \in \bbZ$. In particular
\begin{eqnarray}
b_p^{(2)}(\widetilde{S^1}) ~ = ~ 0 \quad \text{ for all } p \in \bbZ.
\label{L^2-Betti numbers of widetilde{S^1}}
\end{eqnarray}


\subsection{Surfaces}
\label{subsec: Surfaces}Let $F^d_g$%
\indexnotation{F_g^d}
be the orientable closed  surface of genus $g$ with $d$ embedded
 2-disks removed. (As any non-orientable compact surface is
finitely covered by an orientable surface,
it suffices to handle the orientable case by 
\eqref{L^2-Betti numbers of finite coverings}.)
From the value of the zero-th $L^2$-Betti number,
the Euler-Poincar\'e formula and Poincar\'e  duality 
(see Theorem~\ref{the: properties of gen. b_p^{(2)}}
\ref{the: properties of gen. b_p^{(2)}: b_0^{(2)}},
\ref{the: properties of gen. b_p^{(2)}: Euler-Poincar'e formula} and 
\ref{the: properties of gen. b_p^{(2)}: Poincar'e duality})
and from the fact
that a compact surface with boundary is homotopy equivalent
to a bouquet of circles, we conclude
\begin{eqnarray*}
b_0^{(2)}(\widetilde{F^d_g})
& = &
\left\{ \begin{array}{lll} 1 & & \mbox{ if } g =0, d=0,1; \\ 0 & &
\mbox{ otherwise;}
\end{array} \right. 
\\
b_1^{(2)}(\widetilde{F^d_g})
& = &
\left\{ \begin{array}{lll} 0 & & \mbox{ if } g = 0, d=0,1; \\
d + 2\cdot (g-1) & &
 \mbox{ otherwise;}
\end{array} \right. 
\\
b_2^{(2)}(\widetilde{F^d_g})
& = & \left\{ \begin{array}{lll} 1 & & \mbox{ if } g =0, d=0; \\
0 & & \mbox{ otherwise.}
\end{array} \right. 
\end{eqnarray*}
Of course $b_p^{(2)}(\widetilde{F_g^d}) = 0$ for $ p\ge 3$.


\subsection{Three-Dimensional Manifolds}
\label{subsec: Three-Dimensional Manifolds}

In this subsection we state the values of the $L^2$-Betti numbers of 
compact orientable $3$-manifolds.

We begin with collecting some basic notations and
facts about $3$-manifolds. In the sequel $3$-manifold means connected compact orientable
$3$-manifold, possibly with boundary. A
$3$-manifold $M$ is \emph{prime}%
\index{prime $3$-manifold}
\index{manifold!prime $3$-manifold}
if for any decomposition of $M$ as a connected sum $M_1 \#
M_2$, $M_1$ or $M_2$ is homeomorphic to $S^3$. It is
\emph{irreducible}
\index{irreducible $3$-manifold}
\index{manifold!irreducible $3$-manifold}
if every embedded $2$-sphere bounds an
embedded $3$-disk. Every
prime $3$-manifold is either irreducible or is homeomorphic to
$S^1 \times S^2$ \cite[Lemma 3.13]{Hempel(1976)}. 
A $3$-manifold $M$ has a prime decomposition, i.e.
one can write $M$ as a connected sum
\begin{eqnarray*}
M  & = & M_1 \# M_2 \# \ldots \# M_r,
\end{eqnarray*}
where each $M_j$ is prime, and this prime decomposition is
unique up to renumbering and orientation preserving homeomorphism
\cite[Theorems 3.15, 3.21]{Hempel(1976)}.
Recall that a connected $CW$-complex is called \emph{aspherical}%
\index{aspherical}
if $\pi_n(X) = 0$ for $n \ge 2$, or, equivalently, if $\widetilde{X}$ is contractible.
Any aspherical $3$-manifold is homotopy equivalent to an irreducible
$3$-manifold with infinite fundamental group or to a $3$-disk.
By the Sphere Theorem \cite[Theorem 4.3]{Hempel(1976)}, an irreducible 
3-manifold is aspherical if and only if
it is a 3-disk or has infinite fundamental group.

Let us say that a prime 3-manifold is \emph{exceptional}%
\index{exceptional prime $3$-manifold}
\index{manifold!exceptional prime $3$-manifold}
if it is closed
and no finite covering of it is homotopy equivalent
to a Haken, Seifert or hyperbolic 3-manifold. No exceptional prime
3-manifolds are
known. Both Thurston's Geometrization Conjecture and
Waldhausen's Conjecture that any $3$-manifold
is finitely covered by a Haken manifold imply that there are none.

Details of the proof of the following theorem can be found in
\cite[Sections 5 and 6]{Lott-Lueck(1995)}. The proof is quite interesting since
it uses both topological and analytic tools and relies on 
Thurston's Geometrization.

\begin{theorem}[$L^2$-Betti numbers of $3$-manifolds]
\index{Theorem!L2-Betti numbers of $3$-manifolds@$L^2$-Betti numbers of $3$-manifolds}
\label{the: L^2-Betti numbers of exceptional 3-manifolds}
Let $M$ be the connected sum
\mbox{$M_1 \# \ldots \# M_r$} of (compact connected orientable)
prime $3$-manifolds $M_j$ which are non-exceptional. 
Assume that $\pi_1(M)$ is infinite.
Then the $L^2$-Betti numbers of the universal covering
$\widetilde{M}$ are given by
\begin{eqnarray}
b_0^{(2)}(\widetilde{M}) & = & 0; \nonumber
\\
b_1^{(2)}(\widetilde{M}) & = & (r-1) -
\sum_{j=1}^r \frac{1}{\mid \pi_1(M_j)\mid}\; +
\left|\{C \in \pi_0(\partial M) \mid C \cong S^2\}\right|  - \chi(M); \nonumber
\\
b_2^{(2)}(\widetilde{M}) & = & (r-1) -
\sum_{j=1}^r \frac{1}{\mid \pi_1(M_j)\mid} \;+
\left|\{C \in \pi_0(\partial M) \mid C \cong S^2\}\right|; \nonumber
\\
b_3^{(2)}(\widetilde{M}) & = & 0. \nonumber
\end{eqnarray}
\end{theorem}

Notice that in the situation of 
Theorem~\ref{the: L^2-Betti numbers of exceptional 3-manifolds}
the $p$-th $L^2$-Betti number $b_p(\widetilde{M})$ is a rational number.
It is an integer, if $\pi_1(M)$ is torsion-free, and vanishes, if
$M$ is aspherical.


\subsection{Symmetric Spaces}
\label{subsec: Symmetric Spaces}

Let $L$ be a connected semisimple Lie group  with
finite center such that its Lie algebra has no compact ideal.
Let $K \subseteq L$ be a maximal compact subgroup. Then
the manifold $M := L/K$ equipped with a left $L$-invariant Riemannian metric
is a  symmetric space of non-compact
type with $L = \Isom(M)^0$ and $K = \Isom(M)_x^0$,
where $\Isom(M)^0$ is the identity component of the group of isometries
$\Isom(M)$ and $\Isom(M)_x^0$ is the isotropy group 
of some point $x \in M$ under the $\Isom(M)^0$-action.
Every  symmetric space $M$ of non-compact type
can be written in this way. The space $M$ is diffeomorphic to $\bbR^n$. 
Define its \emph{fundamental rank}%
\index{fundamental rank of a symmetric space}
$$\frk(M)%
\indexnotation{frk(M)}
 := \rk_{\bbC}(L) - \rk_{\bbC}(K),$$
where $\rk_{\bbC}(L)$ and  $\rk_{\bbC}(K)$ denotes the so called complex rank of the Lie algebra of $L$
and $K$ respectively (see \cite[page 128f]{Knapp(1986)}). For a compact Lie group $K$ this is the same
as the dimension of a maximal torus. The proof of the next result is due to
Borel~\cite{Borel(1985)}.

\begin{theorem}[$L^2$-Betti numbers of symmetric spaces of non-compact type]
\index{Theorem!L2-Betti numbers of symmetric spaces of non-compact type@
$L^2$-Betti numbers of symmetric spaces of non-compact type}
\label{the: l^2-inv. of sym. sp.}
Let $M$ be a closed Riemannian manifold whose universal covering
$\widetilde{M}$ is a symmetric space of non-compact type. 

Then $b_p^{(2)}(\widetilde{M}) \not= 0$ if and only if $\frk(\widetilde{M}) =  0$
and  $2p = \dim(M)$. 
If $\frk(\widetilde{M}) = 0$, then $\dim(M)$ is even 
and for $2p = \dim(M)$ we get 
$$0 ~ < ~ b_p^{(2)}(\widetilde{M}) ~ = ~ (-1)^p \cdot \chi(M).$$
\end{theorem}

This applies in particular to a hyperbolic manifold and thus we get
the result of Dodziuk~\cite{Dodziuk(1979)}.

\begin{theorem} \label{the: proof of the Chern conjecture for hyperbolic space}
Let $M$ be a hyperbolic closed Riemannian manifold of dimension $n$.
Then
\begin{eqnarray*}
&b_p^{(2)}(\widetilde{M})
\left\{
\begin{array}{ll}
= 0 & \mbox{  if }  2p \not= n
\\
> 0 & \mbox{ if }  2p = n
\end{array}\right.. &
\end{eqnarray*}
If $n$ is even, then
\begin{eqnarray*}
(-1)^{n/2}\cdot \chi(M)
& > &
0.
\end{eqnarray*}
\end{theorem}

The strategy of the proof of
Theorem~\ref{the: proof of the Chern conjecture for hyperbolic space}
is the following.  Because of the Euler-Poincar\'e formula
(see Theorem~\ref{the: properties of gen. b_p^{(2)}}
\ref{the: properties of gen. b_p^{(2)}: Euler-Poincar'e formula}) it suffices to
show that $b_p^{(2)}(\widetilde{M}) = 0$ for $2p \not= n$ and
$b_p^{(2)}(\widetilde{M}) > 0$ for $2p = n$. 
Because of the Hodge-deRham Theorem (see \eqref{Dodziuk's result}) and the facts
that the von Neumann dimension is faithful and $\widetilde{M}$ is isometrically
diffeomorphic to the hyperbolic space $\bbH^n$, it remains to show
that the space of harmonic $L^2$-integrable forms $\calh^p_{(2)}(\bbH^n)$ is trivial
for $2p \not= n$ and non-trivial for $2p = n$. Notice that this question is independent
of $M$ or the $\pi_1(M)$-action. Using the rotational symmetry of
$\bbH^n$, this question is answered positively by Dodziuk~\cite{Dodziuk(1979)}.

More generally one has the following so called Proportionality Principle
(see \cite[Theorem 3.183]{Lueck(2002)}.)

\begin{theorem}[Proportionality Principle for $L^2$-Betti numbers]
\index{Theorem!Proportionality principle for $L^2$-Betti numbers}
\label{the: proportionality principle for $L^2$-Betti numbers}
Let $M$ be a simply connected Riemannian manifold. Then there
are constants $B_p^{(2)}(M)$ for $p\ge 0$ depending only on the Riemannian manifold $M$ with
the following property: For every discrete group $G$
with a cocompact free proper  action on $M$
by isometries the following holds
\begin{eqnarray*}
b_p^{(2)}(M;\caln(G))
& = &
B_p^{(2)}(M) \cdot \vol(G\backslash M).
\end{eqnarray*}
\end{theorem}


\subsection{Spaces with $S^1$-Action}
\label{subsec: Spaces with $S^1$-Action}

The next two theorems are taken from
\cite[Corollary~1.43 and Theorem~6.65]{Lueck(2002)}.

\begin{theorem}{\bf ($L^2$-Betti numbers and $S^1$-actions).}
\index{Theorem!L2-Betti numbers and S^1-actions@$L^2$-Betti numbers and $S^1$-actions}
\label{the: $S^1$-actions and general $L^2$-Betti numbers}
Let $X$ be a connected $S^1$-$CW$-complex.
Suppose that for one orbit $S^1/H$ (and hence for all orbits) the inclusion
into $X$ induces a map on $\pi_1$ with infinite image.
(In particular the $S^1$-action has no fixed points.) 

Then we get 
\begin{eqnarray*}
b^{(2)}_p(\widetilde{X}) & = & 0 \quad \text{ for } p \in \bbZ;\\
\chi(X) & = & 0.
\end{eqnarray*}
\end{theorem}
\proof
We give an outline of the idea of the proof in the case where $X$ is a cocompact
$S^1$-$CW$-complex, because it is a very illuminating example.
The proof in the general case is given in
\cite[Theorem~6.65]{Lueck(2002)}.
It is useful to show the following slightly more general statement
that for any finite $S^1$-$CW$-complex $Y$ and $S^1$-map $f \colon Y \to X$ we get
$b_p^{(2)}(f^*\widetilde{X};\caln(\pi_1(X)) = 0$ for all $p \ge 0$, where $f^*\widetilde{X} \to Y$ is
the pullback of the universal covering $\widetilde{X} \to X$ with $f$. We prove the latter
statement by induction over the dimension and the number of $S^1$-equivariant cells in top
dimension of $Y$. In the induction
step we can assume that $Y$ is an $S^1$-pushout
$$\comsquare{S^1/H \times S^{n-1}}{q}{Z}{}{j}{S^1/H \times D^n}{Q}{Y}$$
for $n = \dim(Y)$. It induces a pushout of free finite $\pi_1(X)$-$CW$-complexes
$$\comsquare{q^*j^*f^*\widetilde{X}}{}{j^*f^*\widetilde{X}}
{}{}
{Q^*f^*\widetilde{X}}{}{f^*\widetilde{X}}$$
The associated long exact Mayer-Vietoris sequence looks like
\begin{multline*}
\ldots H_p(q^*f^*\widetilde{X};\caln(\pi_1(X))) 
\\
\to 
H_p(Q^*f^*\widetilde{X};\caln(\pi_1(X))) \oplus H_p(j^*f^*\widetilde{X};\caln(\pi_1(X)))
\\
\to 
H_p(f^*\widetilde{X};\caln(\pi_1(X)))
\to
H_{p-1}(q^*j^*f^*\widetilde{X};\caln(\pi_1(X))) 
\\
\to 
H_{p-1}(Q^*f^*\widetilde{X};\caln(\pi_1(X))) \oplus H_{p-1}(j^*f^*\widetilde{X};\caln(\pi_1(X)))
\to \ldots
\end{multline*}
Because of the Additivity of the dimension 
(see Theorem~\ref{the: prop. ext. dim.}
\ref{the: prop. ext. dim.: properties of dim: additivity}) it suffices to prove
for all $p \in \bbZ$
\begin{eqnarray*}
\dim_{\caln(\pi_1(X))}\left(H_p(j^*f^*\widetilde{X};\caln(\pi_1(X))) \right)
& = & 0;
\\
\dim_{\caln(\pi_1(X))}\left(H_p(q^*f^*\widetilde{X};\caln(\pi_1(X))) \right)
& = & 0;
\\
\dim_{\caln(\pi_1(X))}\left(H_p(Q^*f^*\widetilde{X};\caln(\pi_1(X)))  \right)
& = & 0.
\end{eqnarray*}
The induction hypothesis applies to $f \circ j\colon Z \to X$ and
$f \circ j \circ q \colon S^1/H \times S^{n-1} \to X$. Hence it remains to
show
$$\dim_{\caln(\pi_1(X))}\left(H_p(Q^*f^*\widetilde{X};\caln(\pi_1(X))) \right) ~ = ~ 0.$$
By elementary covering theory $Q^*f^*\widetilde{X}$ is $\pi_1(X)$-homeomorphic
to $\pi_1(X) \times_j  \widetilde{S^1/H} \times D^n$ for the injective group homomorphism
$j \colon \pi_1(S^1/H) \to \pi_1(X)$ induced by $f \circ Q$. We conclude from
the K\"unneth formula and the compatibility of dimension and induction
(see
Theorem \ref{the: properties of gen. b_p^{(2)}}
\ref{the: properties of gen. b_p^{(2)}: Kuenneth formula} and
\ref{the: properties of gen. b_p^{(2)}: induction})
\begin{eqnarray*}
\dim_{\caln(\pi_1(X))}\left(H_p(Q^*f^*\widetilde{X};\caln(\pi_1(X)))  \right)
& = & b_p^{(2)}(\widetilde{S^1/H}).
\end{eqnarray*}
Since $S^1/H$ is homeomorphic to $S^1$, we get
$b_p^{(2)}(\widetilde{S^1/H}) = 0$ from 
\eqref{L^2-Betti numbers of widetilde{S^1}}. \qed

The next result is taken from \cite[Corollary 1.43]{Lueck(2002)}.

\begin{theorem}
\label{the: fixed point free S^1-actions on Eilenberg MacLane spaces on
aspherical closed manifolds and L^2-Betti numbers}
Let $M$ be an aspherical closed manifold with
non-trivial $S^1$-action. (Non-trivial means that $sx \not= x$ holds for at least
one element $s \in S^1$ and one element $x \in M$). Then the action has no fixed
points and the inclusion of any orbit into $X$ induces an
injection on the fundamental groups. All $L^2$-Betti numbers
$b_p^{(2)}(\widetilde{M})$ are trivial and $\chi(M) = 0$.
\end{theorem}


\subsection{Mapping Tori}
\label{subsec: Mapping Tori}

Let $f\colon X \to X$ be a selfmap. Its \emph{mapping torus}%
\index{mapping torus}
$T_f$%
\indexnotation{T_f}
is obtained from the cylinder
$X \times [0,1]$ by glueing the bottom
to the top by the identification $(x,1) = (f(x),0)$.
There is a canonical map $p\colon T_f \to S^1$
which sends $(x,t)$ to $\exp(2\pi i t)$. It induces a canonical
epimorphism $\pi_1(T_f) \to \bbZ = \pi_1(S^1)$ if $X$ is path-connected.

The following result is taken from 
\cite[Theorem 6.63]{Lueck(2002)}. 

\begin{theorem}[Vanishing of $L^2$-Betti numbers of mapping tori]\ \\
\index{Theorem!Vanishing of $L^2$-Betti numbers of mapping tori}
\label{the: vanishing of general $L^2$-Betti numbers of mapping tori}
Let $f\colon  X \to X$ be a cellular
selfmap of a connected $CW$-complex $X$
and let $\pi_1(T_f) \xrightarrow{\phi} G \xrightarrow{\psi}\bbZ$
be a factorization of the canonical epimorphism into epimorphisms
$\phi$ and $\psi$. Suppose for given $p \ge 0$ that
$b_p^{(2)}(G \times_{\phi \circ i} \widetilde{X};\caln(G)) < \infty$ and
$b_{p-1}^{(2)}(G \times_{\phi \circ i} \widetilde{X};\caln(G)) < \infty$ holds,
where $i\colon  \pi_1(X) \to \pi_1(T_f)$ is the map induced by the obvious inclusion
of $X$ into $T_f$.  Let $\overline{T_f}$ be the covering of $T_f$ associated
to $\phi$, which is a free $G$-$CW$-complex. Then we get
\begin{eqnarray*}
b_p^{(2)}(\overline{T_f};\caln(G))
& = &
0.
\end{eqnarray*}
\end{theorem}
\proof We give the proof in the special case where $X$ is a connected finite
$CW$-complex and $\phi = \id$, i.e. we show for a connected finite $CW$-complex $X$
that $b_p^{(2)}(\widetilde{T_f}) = 0$ for all $p \ge 0$. 
For each positive integer $d$ there is a finite $d$-sheeted covering
$\overline{T_f} \to T_f$ associated to the subgroup of index $d$ in $\pi_1(T_f)$
which is the preimage of $d\bbZ \subseteq \bbZ$ under the canonical
homomorphism $\pi_1(T_f) \to \bbZ$.  There is a homotopy equivalence
$T_{f^d} \to \overline{T_f}$.
We conclude from
\eqref{L^2-Betti numbers of finite coverings} and homotopy invariance of
$L^2$-Betti numbers
(see Theorem~\ref{the: properties of gen. b_p^{(2)}}
\ref{the: properties of gen. b_p^{(2)}: homology invariance})
$$b_p^{(2)}(\widetilde{T_f}) ~ = ~ \frac{b_p^{(2)}(\widetilde{T_{f^d}})}{d}.$$
There is a $CW$-complex structure on  $T_{f^d}$ with 
$\beta_p(X) + \beta_{p-1}(X)$ $p$-cells, if $\beta_p(X)$ is the number of $p$-cells in $X$.
We conclude from Additivity of the dimension function
(see Theorem~\ref{the: prop. ext. dim.}
\ref{the: prop. ext. dim.: properties of dim: additivity})
\begin{multline*} b_p^{(2)}(\widetilde{T_{f^d}})  
~ \le ~  
\dim_{\caln(\pi_1(T_{f^d}))}
\left(\caln(\pi_1(T_{f^d})) \otimes_{\bbZ\pi_1(T_{f^d})}  C_p(\widetilde{T_{f^d}})\right)
\\
~ = ~\beta_p(X) + \beta_{p-1}(X).
\end{multline*}
This implies for all positive integers $d$
$$0 ~ \le ~  b_p^{(2)}(\widetilde{T_f}) ~ \le ~ \frac{\beta_p(X) + \beta_{p-1}(X)}{d}.$$
Taking the limit for $d \to \infty$ implies $ b_p^{(2)}(\widetilde{T_f}) = 0$. \qed


\subsection{Fibrations}
\label{subsec: Fibrations}

The next result is proved in \cite[Lemma 6.6. and Theorem 6.67]{Lueck(2002)}.
The proof is based on standard spectral sequence arguments and 
the fact that the dimension function is defined for
arbitrary $\caln(G)$-modules.

\begin{theorem}[$L^2$-Betti numbers and fibrations]
\index{Theorem!L2-Betti numbers and fibrations@$L^2$-Betti number and fibrations}
\label{the: L^2-Betti numbers and fibrations}
\ \\
\begin{enumerate} 

\item \label{the: L^2-Betti numbers and fibrations: d to d+1}
Let $F \xrightarrow{i} E \xrightarrow{p} B$ be a fibration
of connected $CW$-complexes. Consider a factorization
$p_*\colon  \pi_1(E) \xrightarrow{\phi} G \xrightarrow{\psi} \pi_1(B)$
of the map induced by $p$ into epimorphisms
$\phi$ and $\psi$. Let $i_*\colon  \pi_1(F) \to \pi_1(E)$
be the homomorphism induced by the inclusion $i$.
Suppose for a given integer $d \ge 1$
that $b_p^{(2)}(G \times_{\phi \circ i_*}\widetilde{F};\caln(G))= 0$
for $p \le d-1$ and
$b_d^{(2)}(G \times_{\phi \circ i_*}\widetilde{F};\caln(G)) < \infty$
holds. Suppose that $\pi_1(B)$ contains an element
of infinite order or finite subgroups of arbitrarily large order. Then
$b_p^{(2)}(G \times_{\phi} \widetilde{E};\caln(G)) = 0$ for $p \le d$;

\item \label{the: L^2-Betti numbers and fibrations: vanishing of all L^2-Betti numbers}
Let $F \xrightarrow{i} E \to B$ be a fibration
of connected $CW$-complexes. Consider a group homomorphism
$\phi\colon  \pi_1(E) \to G$. Let $i_*\colon  \pi_1(F) \to \pi_1(E)$ be the
homomorphism induced by the inclusion $i$.
Suppose that for a given integer $d \ge 0$ the $L^2$-Betti number
$b_p^{(2)}(G \times_{\phi \circ i_*}\widetilde{F};\caln(G))$
vanishes for all $p \le d$. Then the $L^2$-Betti number
$b_p^{(2)}(G \times_{\phi} \widetilde{E};\caln(G))$
vanishes for all $p \le d$.
\end{enumerate}

\end{theorem}


\typeout{--------------------   Section 4 --------------------------}

\section{The Atiyah Conjecture}
\label{subsec: The Atiyah Conjecture}

In this section we discuss the Atiyah Conjecture

\begin{conjecture}[Atiyah Conjecture] \label{con: Atiyah Conjecture}
\index{Conjecture!Atiyah Conjecture}
Let $G$ be a discrete group with an upper bound 
on the orders of its finite subgroups. 
Consider $d \in \bbZ$, $d \ge 1$ such that the order of every finite subgroup of $G$ divides $d$.
Let $F$ be a field with $\bbQ \subseteq F \subseteq \bbC$.
The \emph{Atiyah Conjecture} for $(G,d,F)$ says 
that for any finitely presented $FG$-module $M$ we have 
$$d \cdot \dim_{\caln(G)}\left(\caln(G) \otimes_{FG} M\right) ~ \in ~ \bbZ.$$
\end{conjecture}


\subsection{Reformulations of the Atiyah Conjecture}

We present equivalent reformulations of the Atiyah 
Conjecture~\ref{con: Atiyah Conjecture}.

\begin{theorem}[Reformulations of the Atiyah Conjecture]
\label{the: Reformulations of the Atiyah Conjecture}
\index{Theorem!Reformulations of the Atiyah Conjecture}

Let $G$ be a discrete group.
Suppose that there exists $d \in \bbZ$, $d \ge 1$ such that the order 
of every finite subgroup of $G$ divides $d$. 
Let $F$ be a field with $\bbQ \subseteq F \subseteq \bbC$. 
Then the following assertions are equivalent:
\begin{enumerate}

\item \label{the: Reformulations of the  Atiyah Conjecture: original}
The  Atiyah Conjecture~\ref{con:  Atiyah Conjecture}
is true for $(G,d,F)$, i.e. for every finitely presented  $FG$-module $M$ we have 
$$d \cdot \dim_{\caln(G)}\left(\caln(G) \otimes_{FG} M\right) ~ \in ~ \bbZ;$$

\item \label{the: Reformulations of the  Atiyah Conjecture: arbitrary module}
For every  $FG$-module $M$ we have 
$$d \cdot \dim_{\caln(G)}\left(\caln(G) \otimes_{FG} M\right) ~ \in ~ \bbZ \amalg \{\infty\}.$$

\end{enumerate}
\end{theorem}
\proof See \cite[Lemma 10.7 and Remark 10.11]{Lueck(2002)}. \qed

We mention that the  Atiyah Conjecture~\ref{con:  Atiyah Conjecture}
is true for $(G,d,F)$ if and only if for any finitely generated
subgroup $H \subseteq G$ the Atiyah Conjecture~\ref{con:  Atiyah Conjecture}
is true for  $(H,d,F)$ (see \cite[Lemma 10.4]{Lueck(2002)}).

The next result explains that the  Atiyah 
Conjecture~\ref{con:  Atiyah Conjecture}
for $(G,d,\bbQ)$ for a finitely generated group $G$ is a statement about the possible values of
$L^2$-Betti numbers.

\begin{theorem}[Reformulations of the  Atiyah Conjecture for $F = \bbQ$]
\label{the: Reformulations of the  Atiyah Conjecture for F = Q}
\index{Theorem!Reformulations of the Atiyah Conjecture for $F = \bbQ$}

Let $G$ be a finitely generated group with an upper bound $d \in \bbZ, d \ge 1$ 
on the orders of its finite subgroups. 
Then the following assertions are equivalent:
\begin{enumerate}

\item \label{the: Reformulations of the  Atiyah Conjecture for F = Q: original}
The  Atiyah Conjecture~\ref{con:  Atiyah Conjecture}
is true for $(G,d,\bbQ)$;

\item \label{the: Reformulations of the  Atiyah Conjecture for F = Q: smooth manifold}
For every  free proper smooth cocompact $G$-manifold $M$ without boundary and $p \in \bbZ$ we have
$$d \cdot b_p^{(2)}(M;\caln(G)) ~ \in ~ \bbZ;$$

\item \label{the: Reformulations of the  Atiyah Conjecture for F = Q: cocompact G-CW}
For every  finite  free $G$-$CW$-complex $X$ and $p \in \bbZ$ we have
$$d \cdot b_p^{(2)}(X;\caln(G)) ~ \in ~ \bbZ;$$

\item \label{the: Reformulations of the  Atiyah Conjecture for F = Q: G-space}
For every  $G$-space $X$ and $p \in \bbZ$ we have
$$d \cdot b_p^{(2)}(X;\caln(G)) ~ \in ~ \bbZ \amalg \{\infty\}.$$

\end{enumerate}
\end{theorem}
\proof This follows from  \cite[Lemma 10.5]{Lueck(2002)} and
Theorem~\ref{the: Reformulations of the  Atiyah Conjecture}. \qed

We mention that all the explicit computations presented in
Section~\ref{sec: Computations of L^2-Betti Numbers} 
are compatible with the Atiyah Conjecture~\ref{con: Atiyah Conjecture}.


\subsection{The Ring Theoretic Version of the Atiyah Conjecture}
\label{subsec: The Ring Theoretic Version of the Atiyah Conjecture}

In this subsection we consider the following \emph{fundamental square of ring extensions}%
\index{fundamental square of ring extensions}
\begin{eqnarray}
\comsquare{\bbC G}{i}{\caln(G)}{j}{k}{\cald(G)}{l}{\calu(G)}
\label{fundamental square of ring extensions}
\end{eqnarray}
which we explain next. 

As before $\bbC G$ is the complex group ring and $\caln(G)$ is the group von Neumann algebra.

By $\calu(G)$%
\indexnotation{calu(G)} we denote the algebra of affiliated operators. Instead of
its functional analytic definition we describe it algebraically, namely, it is the Ore localization
of $\caln(G)$ with respect to the multiplicative subset of
non-trivial zero-divisors in $\caln(G)$. The proof that this multiplicative subset satisfies the Ore
condition and basic definitions and properties of Ore localization and of $\calu(G)$
can be found for instance in
\cite[Sections 8.1 and 8.2]{Lueck(2002)}. In particular
$\calu(G)$ is flat when regarded as an $\caln(G)$-module.  Moreover, 
the ring $\calu(G)$ is a 
\emph{von Neumann regular ring},%
\index{ring!von Neumann regular}
\index{von Neumann regular ring}
i.e. every finitely generated submodule of a projective module is a direct summand.
This is a stronger condition than being semihereditary.

Given a finitely generated projective $\calu(G)$-module $Q$,
there is a finitely generated projective $\caln(G)$-module $P$
such that $\calu(G) \otimes_{\caln(G)} P$ and $Q$ are $\calu(G)$-isomorphic.
If $P_0$ and $P_1$ are two finitely generated projective $\caln(G)$-modules,
then $P_0 \cong_{\caln(G)} P_1 \Leftrightarrow
\calu(G) \otimes_{\caln(G)} P_0 \cong_{\calu(G)}
\calu(G) \otimes_{\caln(G)}P_1$. This enables us to define a dimension function for 
$\dim_{\calu(G)}$ with properties  analogous  to $\dim_{\caln(G)}$
(see \cite[Section 8.3]{Lueck(2002)}, \cite{Reich(1999)} or 
\cite{Reich(2001)}).

\begin{theorem}{\bf (Dimension function for arbitrary $\calu(G)$-modules).}
\index{Theorem!Dimension function for arbitrary $\calu(G)$-modules}
\label{the: prop. ext. dim. for calu}

There exists precisely one dimension function
$$\dim_{\calu(G)}%
\indexnotation{dim_{calu(G)}(M)}
 \colon \{\calu(G)\text{-modules}\} ~ \to ~  [0,\infty]%
$$
which satisfies:

\begin{enumerate}

\item
\label{the: prop. ext. dim. for calu: extension property}
Extension Property\par\noindent
If $M$ is an $\caln(G)$-module, then
$$\dim_{\calu(G)}\left(\calu(G) \otimes_{\caln(G)} M\right) ~ = ~ \dim_{\caln(G)}(M);$$

\item
\label{the: prop. ext. dim. for calu: additivity}
Additivity \par\noindent
If $0 \to M_0 \to M_1 \to M_2 \to 0$ is an
exact sequence of $\calu(G)$-modules, then
$$\dim_{\calu(G)}(M_1) ~ = ~ \dim_{\calu(G)}(M_0) + \dim_{\calu(G)}(M_2);$$

\item
\label{the: prop. ext. dim. for calu: cofinality}
Cofinality \par\noindent
Let $\{M_i\mid i \in I\}$ be a cofinal system of submodules of $M$.
Then
$$\dim_{\calu(G)}(M) ~ = ~ \sup\{\dim_{\calu(G)}(M_i) \mid i \in I\};$$

\item
\label{the: prop. ext. dim. for calu: properties of dim: continuity}
Continuity \par\noindent
If $K \subseteq M$ is a submodule of the finitely generated
$\calu(G)$-module $M$, then
$$\dim_{\calu(G)}(K) ~ = ~ \dim_{\calu(G)}(\overline{K}).$$

\end{enumerate}

\end{theorem}

\begin{remark}[Comparing $\bbZ \subseteq \bbQ$ and $\caln(G) \subseteq \calu(G)$]
\label{rem: Comparing Z subseteq Q and calkn(G) subseteq calu(G)}
\em
Recall the Slogan~\ref{slo: caln(G) like Z} that the
 group von Neumann algebra $\caln(G)$ behaves like the ring of integers $\bbZ$, provided
one ignores the properties integral domain and Noetherian.  This is supported
by the construction and properties of $\calu(G)$. Obviously $\calu(G)$
 plays the same role for $\caln(G)$ as $\bbQ$ plays for $\bbZ$
as the definition of $\calu(G)$ as the Ore localization of $\caln(G)$ with respect to the multiplicative
subset of non-zero-divisors and 
Theorem~\ref{the: prop. ext. dim. for calu} show.
\em
\end{remark}

A subring $R\subseteq S$ is called \emph{division closed}%
\index{ring!division closed} 
if each element in $R$, which is invertible in $S$, is already invertible in $R$.
It is called  \emph{rationally closed}%
\index{ring!rationally closed} 
if each square matrix over $R$, which is invertible over $S$,
is already invertible over $R$.
Notice that the intersection of division closed subrings of $S$
is again division closed, and analogously for rationally closed subrings.
Hence the following definition
makes sense.

\begin{definition}[Division and rational closure]
\label{def: division and rational closure}
Let $S$ be a ring with subring $R \subseteq S$. The \emph{division closure}%
\index{division closure}
$\cald(R \subseteq S)$
\indexnotation{cald(R subset S)}
or  \emph{rational closure}%
\index{rational closure}
$\calr(R \subseteq S)$
\indexnotation{calr(R subset S)}
respectively is the smallest subring of $S$ which contains $R$ and is
division closed or rationally closed respectively.
\end{definition}

The ring $\cald(G)$ appearing in the fundamental square 
\eqref{fundamental square of ring extensions}
is the rational closure of $\bbC G$ in $\calu(G)$.

\begin{conjecture}[Ring theoretic version of the Atiyah Conjecture]
\label{con: Ring theoretic version of the Atiyah Conjecture}
\index{Conjecture!Ring theoretic version of the Atiyah Conjecture}
Let $G$ be a group for which there exists an upper bound on the orders of
its finite subgroups. Then:
\begin{enumerate}

\item[{\bf (R)}] The ring $\cald(G)$ is semisimple;

\item[{\bf (K)}] The composition
$$\bigoplus_{H \subseteq G, |H| < \infty}~ K_0(\bbC H) \xrightarrow{a} K_0(\bbC G) \xrightarrow{j} 
K_0(\cald(G))$$
is surjective, where $a$ is induced by the various inclusions $H \to G$.

\end{enumerate}
\end{conjecture}

\begin{lemma} \label{lem: Ring theoretic Atiyah implies Atiyah}
Let $G$ be a group. Suppose that there exists $d \in \bbZ$, $d \ge 1$ such that the order 
of every finite subgroup of $G$ divides $d$. 
If the group $G$ satisfies the
ring theoretic version of the Atiyah 
Conjecture~\ref{con: Ring theoretic version of the Atiyah Conjecture}, 
then the Atiyah Conjecture~\ref{con:  Atiyah Conjecture} for $(G,d,\bbC)$
is true.
\end{lemma}
 \proof Let $M$ be a finitely presented $\bbC G$-module. 
Then $\cald(G) \otimes_{\bbC G}M$ is a finitely generated projective
$\cald(G)$-module since $\cald(G)$ is semisimple by assumption. We obtain a well-defined homomorphism 
of abelian groups
$$D \colon K_0(\cald(G)) \to \bbR, \quad [P] \mapsto \dim_{\calu(G)}\left(\calu(G) \otimes_{\cald(G)} P\right).$$
Because of the fundamental square \eqref{fundamental square of ring extensions}
and Theorem~\ref{the: prop. ext. dim. for calu}
\ref{the: prop. ext. dim. for calu: extension property}
we have
$$\dim_{\caln(G)}(\caln(G) \otimes_{\bbC G} M) ~ = ~
D([ \cald(G) \otimes_{\bbC G} M ]).$$
Hence it suffices to show that $d \cdot \im(D)$ is contained in $\bbZ$.
Because of assumption ${\bf (K)}$ it suffices to check for
each finite subgroup $H \subseteq G$ and each finitely generated projective $\bbC H$-module $P$
$$d \cdot \dim_{\calu(G)}(\calu(G) \otimes_{\bbC G} \bbC G \otimes_{\bbC H} P) ~ \in ~ \bbZ.$$
Example~\ref{exa: dim for finite G} and
Theorem~\ref{the: Induction and Dimension} imply
\begin{eqnarray*}
\dim_{\calu(G)}(\calu(G) \otimes_{\bbC G} \bbC G \otimes_{\bbC H} P)
& = & 
\dim_{\caln(G)}(\caln(G) \otimes_{\bbC G} \bbC G \otimes_{\bbC H} P)
\\
& = & 
\dim_{\caln(G)}(\caln(G) \otimes_{\caln(H)} P)
\\
& = & 
\dim_{\caln(H)}(P)
\\
& = & 
\frac{\dim_{\bbC}(P)}{|H|}.
\end{eqnarray*}
Obviously $d \cdot \frac{\dim_{\bbC}(P)}{|H|}  ~ \in ~ \bbZ$. \qed


\subsection{The Atiyah Conjecture for Torsion-Free Groups}

\begin{remark}[The Atiyah Conjecture in the torsion-free case]
\label{rem: The Atiyah Conjecture in the torsion-free case} \em
Let $G$ be a torsion-free group. Then we can choose $d = 1$ in
the  Atiyah Conjecture~\ref{con:  Atiyah Conjecture}.
The  Atiyah Conjecture~\ref{con:  Atiyah Conjecture} for $(G,1,F)$ says that 
$\dim_{\caln(G)}(\caln(G) \otimes_{FG} M) ~\in~\bbZ$ holds for every 
finitely presented $FG$-module $M$  and Theorem 
\ref{the: Reformulations of the  Atiyah Conjecture} says
that then this holds automatically for all
$FG$-modules $M$ with $\dim_{\caln(G)}(\caln(G) \otimes_{FG} M) < \infty$. 
In the case, where $F = \bbQ$ and $G$ is a torsion-free finitely generated group $G$,
Theorem~\ref{the: Reformulations of the  Atiyah Conjecture for F = Q} implies that the
 Atiyah Conjecture~\ref{con:  Atiyah Conjecture}
 for $(G,1,F)$ is equivalent to the statement that 
$b_p^{(2)}(X;\caln(G)) ~ \in ~ \bbZ$ is true for all $G$-spaces $X$.\em
\end{remark}

\begin{remark}[The ring theoretic version of the Atiyah Conjecture in the torsion-free case]
\label{rem: The Ring theoretic version of the Atiyah Conjecture in the torsion-free case} \em
Let $G$ be a torsion-free group. Then the ring theoretic version of the Atiyah 
Conjecture~\ref{con: Ring theoretic version of the Atiyah Conjecture} reduces
to the statement that $\cald(G)$ is a skewfield. In this case
we can assign to every $\cald(G)$-module $N$ its dimension
$\dim_{\cald(G)}(N) \in \bbZ \amalg \{\infty\}$ in the usual way and we get for every
$\bbC G$-module $M$
\begin{multline*}
\dim_{\caln(G)}(\caln(G) \otimes_{\bbC G} M) ~ = ~
\dim_{\calu(G)}(\calu(G) \otimes_{\bbC G} M) ~ = ~
\dim_{\cald(G)}(\cald(G) \otimes_{\bbC G} M).
\end{multline*} 
\em
\end{remark}

\begin{example}[The case $G = \bbZ^n$] \label{exa: fundamental square for G = Z^n}
\em
In the case $G = \bbZ^n$ the fundamental square
of ring extensions~\eqref{fundamental square of ring extensions}
can be identified with
$$
\comsquare{\bbC[\bbZ^n]}{}{L^{\infty}(T^n)}
{}{}
{\bbC[\bbZ^n]_{(0)}}{}{MF(T^n)}
$$
where $MF(T^n)$ the ring of equivalence classes of measurable functions
$T^n \to \bbC$. We have already proved
$$\dim_{\caln(\bbZ^n)}(\caln(\bbZ^n) \otimes_{\bbC[\bbZ^n]} M) 
~ = ~
\dim_{\bbC[\bbZ^n]_{(0)}}(\bbC[\bbZ^n]_{(0)} \otimes_{\bbC[\bbZ^n]} M)
$$
in Example~\ref{exa: dimension for G = Z^n}.
\em
\end{example}


\subsection{The Atiyah Conjecture Implies the Kaplanski Conjecture}

The following conjecture is a well-known conjecture about group rings.

\begin{conjecture}[Kaplanski Conjecture]
\label{con: Kaplanski Conjecture}
\index{Conjecture!Kaplanski Conjecture}
Let $F$ be a field and let $G$ be a torsion-free group. 
Then $FG$ contains no non-trivial zero-divisors.
\end{conjecture}

\begin{theorem}[The Atiyah and the Kaplanski Conjecture]
\label{the: The Atiyah and the Kaplanski Conjecture}
\index{Theorem!The Atiyah and the Kaplanski Conjecture}
Let $G$ be a torsion-free group and let $F$ be a field with $\bbQ \subseteq F \subseteq \bbC$.
Then the Atiyah Conjecture~\ref{con:  Atiyah Conjecture}
for $(G,1,F)$  implies the Kaplanski Conjecture~\ref{con: Kaplanski Conjecture}
for $F$ and $G$.
\end{theorem}
\proof
Let $u \in FG$ be a zero-divisor. Then the kernel of the 
$\caln(G)$-map 
$r_u \colon \caln(G) \to \caln(G)$ given by right multiplication with $u$
is non-trivial. Since $\caln(G)$ is semihereditary, the image of $r_u$ is projective.
Hence both $\ker(r_u)$ and $\caln(G)/\ker(r_u)$ are finitely generated projective
$\caln(G)$-modules. Additivity of $\dim_{\caln(G)}$ implies 
$$0 < \dim_{\caln(G)}((\ker(r_u)) \le \dim_{\caln(G)}(\caln(G)) = 1.$$
We conclude from Remark~\ref{rem: The Atiyah Conjecture in the torsion-free case}
that $\dim_{\caln(G)}(\ker(r_u))$ is an integer. Additivity of
$\dim_{\caln(G)}$ implies 
$$\dim_{\caln(G)}\left(\caln(G)/\ker(r_u)\right) ~ = ~ 0.$$
We conclude $\caln(G)/\ker(r_u) = 0$ and hence $u = 0$. \qed


\subsection{The Status of the Atiyah Conjecture}
\label{subsec: The Status of the Atiyah Conjecture}

Let $l^{\infty}(G,\bbR)$ be the space of equivalence classes of bounded
functions from $G$ to $\bbR$ with the supremum norm.
Denote by $1$ the constant function with value $1$.

\begin{definition}[Amenable group] \label{def: amenable group}
A group $G$ is called \emph{amenable},
\index{group!amenable}
if there is a (left) $G$-invariant linear operator
$\mu\colon  l^{\infty}(G,\bbR) \to \bbR$
with $\mu(1) = 1$, which satisfies for all $f \in l^{\infty}(G,\bbR)$
$$
\inf\{f(g)\mid g \in G\} \le
\mu(f)  \le \sup\{f(g)\mid g\in G\}.
$$
The latter condition is equivalent to the condition that $\mu$ is
bounded and $\mu(f) \ge 0$ if $f(g) \ge 0$ for all
$g \in G $. 
\end{definition}

\begin{definition}[Elementary amenable group] \label{def: elementary amenable group}
The \emph{class of elementary amenable groups}%
\index{group!elementary amenable}
$\caleam$%
\indexnotation{caleam}
 is defined as the smallest
class of groups which has the following properties:
\begin{enumerate}

\item It contains all finite and all abelian groups;

\item It is closed under taking subgroups;

\item It is closed under taking quotient groups;

\item It is closed under extensions%
\index{closed under!extensions}, i.e. if
$1 \to H \to G \to K \to 1$ is an exact sequence of groups and
$H$ and $K$ belong to $\caleam$, then also $G \in\caleam$;

\item It is closed under \emph{directed unions},
\index{closed under!directed unions}
i.e. if $\{G_i \mid i \in I\}$ is a directed system of subgroups
such that $G = \bigcup_{i \in I} G_i$ and each $G_i$ belongs to $\caleam$, then
$G \in \caleam$. (Directed means that for two indices $i$ and $j$
there is a third index $k$ with $G_i,G_j \subseteq G_k$.)

\end{enumerate}
\end{definition}

The class of amenable groups satisfies all the conditions
appearing in Definition~\ref{def: elementary amenable group}.
Hence every elementary amenable group is amenable. The converse is not true.

\begin{definition}[Linnell's class of groups $\calc$]
 \label{def: Linnell's class calc}
Let $\calc$
\indexnotation{calc}
be the smallest class of groups, which contains all
free groups and is closed under directed unions and extensions with
elementary amenable quotients.
\end{definition}

The next result is due to Linnell \cite{Linnell(1993)}.

\begin{theorem}[Linnell's Theorem] \label{the: Linnells theorem}
\index{Theorem!Linnell's Theorem}
Let $G$ be a group in $\calc$.
Suppose that there exists $d \in \bbZ$, $d \ge 1$ such that the order 
of every finite subgroup of $G$ divides $d$. Then the ring theoretic version of the Atiyah 
Conjecture~\ref{con: Ring theoretic version of the Atiyah Conjecture} for $G$ 
and hence the Atiyah Conjecture~\ref{con: Atiyah Conjecture} for $(G,d,\bbC)$ are true.
\end{theorem}

The next definition and the next theorem are due to Schick
\cite{Schick(2000c)}.

\begin{definition} \label{def: class cald}
Let $\cald$%
\indexnotation{cald}
be the smallest non-empty class of groups such that
\begin{enumerate}

\item If $p\colon G \to A$ is an epimorphism of a torsion-free group $G$ onto an
elementary amenable group $A$ and if $p^{-1}(B) \in \cald$ for every finite group
$B \subseteq A$, then $G \in \cald$;

\item $\cald$ is closed under taking subgroups;

\item $\cald$ is closed under colimits and inverse limits over directed systems.

\end{enumerate}

\end{definition}

\begin{theorem} \label{the: Atiyah Conjecture for cald}
\begin{enumerate}

\item \label{the: Atiyah Conjecture for cald: Atiyah}
If the group $G$ belongs to $\cald$,
then $G$ is torsion-free and the 
 Atiyah Conjecture~\ref{con: Atiyah Conjecture}  for $(G,1,\bbQ)$ is true for $G$;

\item \label{the: Atiyah Conjecture for cald: properties of cald}
The class $\cald$ is closed under direct sums, direct products and free
products. Every residually torsion-free elementary amenable group
belongs to $\cald$.

\end{enumerate}
\end{theorem}

More information about the status of the Atiyah Conjecture~\ref{con: Atiyah Conjecture}
can be found for instance in \cite[Subsection 10.1.3]{Lueck(2002)}.


\subsection{Groups Without Bound on the Order of Its Finite Subgroups}
\label{subsec: Groups Without Bound on the Order of Its Finite Subgroups}

Given a group $G$, let $\calfin(G)$%
\indexnotation{calfin(G)} be the set of finite subgroups
of $G$. Denote by
\begin{eqnarray}
\frac{1}{|\calfin(G)|}\bbZ%
\indexnotation{frac{1}{midcalfin(G)mid}zz}
 & \subseteq & \bbQ
\label{1/|calfin G| Z}
\end{eqnarray}
the additive subgroup of $\bbR$ generated by the set of rational
numbers $\{\frac{1}{|H|} \mid H \in \calfin(G)\}$.

There is the following formulation of the Atiyah Conjecture
for arbitrary groups in the literature. 

\begin{conjecture}[Atiyah Conjecture for arbitrary groups $G$]
\label{con: Atiyah Conjecture for arbitrary groups G}
A group $G$ satisfies the \emph{Atiyah Conjecture}%
\index{Conjecture!Atiyah Conjecture for arbitrary groups}
if for every finitely presented $\bbC G$-module $M$ we have
\begin{eqnarray*}
\dim_{\caln(G)}(\caln(G) \otimes_{\bbC G} M) 
& \in & \frac{1}{|\calfin(G)|}\bbZ.
\end{eqnarray*}
\end{conjecture}

There do exist counterexamples to this conjecture.
The \emph{lamplighter group}%
\index{group!lamplighter group}
\index{lamplighter group}
$L$ is defined by the semidirect product
$$L := \left(\bigoplus_{n \in \bbZ} \bbZ/2\right)  \rtimes \bbZ$$
with respect to the shift automorphism
of  $\bigoplus_{n \in \bbZ} \bbZ/2$, which sends
$(x_n)_{n \in \bbZ}$ to $(x_{n-1})_{n \in \bbZ}$. Let $e_0 \in \bigoplus_{n \in \bbZ} \bbZ/2$
be the element whose entries are all zero except the entry at $0$. Denote by
$t \in \bbZ$ the standard generator of $\bbZ$ which we will also  view as an element of
$L$. Then $\{e_0t,t\}$ is a set of generators
for $L$. The associated \emph{Markov operator}%
\index{Markov operator}
$M\colon l^2(G) \to l^2(G)$ is given by right multiplication with
$\frac{1}{4} \cdot(e_0t + t + (e_0t)^{-1} + t^{-1})$.
It is related to the Laplace operator $\Delta_0\colon l^2(G) \to l^2(G)$ of the Cayley graph
of $G$ by $\Delta_0 = 4\cdot \id - 4 \cdot M$.
The following result is a special case of the main
 result in the paper of Grigorchuk and {\.Z}uk
\cite[Theorem 1 and Corollary 3]{Grigorchuk-Zuk(2001)}
(see also \cite{Grigorchuk-Linnell-Schick-Zuk(2000)}). An elementary
proof can be found in \cite{Dicks-Schick(2002)}.

\begin{theorem}[Counterexample to the Atiyah Conjecture for arbitrary groups]
\index{Theorem!Counterexample to the Atiyah Conjecture for arbitrary groups}
 \label{the: counterexample to Atiyah}
The von Neumann dimension of the kernel of the Markov operator $M$
of the lamplighter group $L$ associated to the set of generators
 $\{e_0t,t\}$ is $1/3$. In particular $L$ does not satisfy the
Atiyah Conjecture \ref{con: Atiyah Conjecture for arbitrary groups G}.
\end{theorem}

To the author's knowledge there is no example of a group $G$ for
which there is a finitely presented $\bbC G$-module $M$ such that
$\dim_{\caln(G)}(\caln(G) \otimes_{\bbZ G}M)$ is irrational.

Let $A= \bigoplus_{n \in \bbZ} \bbZ/2$. Because this group is locally
finite, it satisfies the Atiyah Conjecture for arbitrary 
groups~\ref{con: Atiyah Conjecture for arbitrary groups G}, i.e.
$\dim_{\caln(G)}(\caln(G) \otimes_{\bbC A} M) \in \bbZ[1/2]$ for every finitely
presented $\bbC A$-module $M$.
On the other hand, each non-negative real number $r$ can be realized
as $\dim_{\caln(G)}(\caln(G) \otimes_{\bbC A} M)$ for a finitely
generated $\bbC A$-module (see \cite[Example 10.13]{Lueck(2002)}).
Notice that there is no upper bound on the orders of finite subgroups
of $A$, so that this is no contradiction to
Theorem~\ref{the: Reformulations of the Atiyah Conjecture}.


\typeout{--------------------   Section 5 --------------------------}

\section{Flatness Properties of the Group von Neumann Algebra}
\label{sec: Flatness Poperties of caln(G) over CG}

The proof of  next result can be found in  \cite[Theorem 5.1] {Lueck(1998b)} or
\cite[Theorem 6.37]{Lueck(2002)}.

\begin{theorem}{\bf (Dimension-flatness of $\caln(G)$ over $\bbC G$ for
  amenable $G$).}
\index{Theorem!Dimension-flatness of $\caln(G)$ over $\bbC G$ for
  amenable $G$}
\label{the: dimension of higher Tor-s vanish in the amenable case}
Let $G$ be amenable and $M$ be a $\bbC G$-module. Then
\begin{eqnarray*}
\dim_{\caln(G)}\left(\Tor^{\bbC G}_p(\caln(G),M)\right)
& = & 0
\hspace{10mm} \mbox{ for } p \ge 1,
\end{eqnarray*}
where we consider $\caln(G)$ as an
$\caln(G)$-$\bbC G$-bimodule.
\end{theorem}

It implies using an easy spectral sequence argument

\begin{theorem}[$L^2$-Betti numbers and homology in the amenable case]
\index{Theorem!$L^2$-Betti numbers and homology in the amenable case}
\label{the: L^2-Betti numbers and homology in the amenable case}
Let $G$ be an amenable group and $X$ be a $G$-space. Then

\begin{enumerate} 

\item \label{the: L^2-Betti numbers and homology in the amenable case: L^2-Betti numbers}
$b_p^{(2)}(X;\caln(G)) ~ = ~
\dim_{\caln(G)}\left(\caln(G) \otimes_{\bbC G}
H_p^{\sing}(X;\bbC)\right);$

\item \label{the: L^2-Betti numbers and homology in the amenable case: L^2-Euler characteristic}
Suppose that $X$ is a $G$-$CW$-complex with $m(X;G) < \infty$. Then
\begin{eqnarray*}
\chi^{(2)}(X) & = & \sum_{c \in I(X)} (-1)^{\dim(c)} \cdot |G_c|^{-1} 
\\ 
& = & \sum_{p \ge 0} (-1)^p \cdot \dim_{\caln(G)}
\left(\caln(G) \otimes_{\bbC G} H_p(X;\bbC)\right).
\end{eqnarray*}
\end{enumerate}
\end{theorem}

Further applications of 
Theorem~\ref{the: dimension of higher Tor-s vanish in the amenable case} will be discussed
in Section~\ref{sec: Applications to Group Theory} and 
Section~\ref{sec: $G$- and $K$-Theory}.

\begin{conjecture}{\bf (Amenability and dimension-flatness of $\caln(G)$ over $\bbC G$).}
 \label{con: characterization of amenable by dim(Tor)}%
\index{Conjecture!Amenability and dimension-flatness of $\caln(G)$ over $\bbC G$}
A group $G$ is amenable if and only if for every $\bbC G$-module
$M$
\begin{eqnarray*}
\dim_{\caln(G)}\left(\Tor^{\bbC G}_p(\caln(G),M)\right)
& = & 0
\hspace{10mm} \mbox{ for } p \ge 1
\end{eqnarray*}
holds.
\end{conjecture}

\begin{remark}[Evidence for Conjecture~\ref{con: characterization of amenable by dim(Tor)}]
\label{rem: Evidence for Conjecture on dimension-flatness}
\em
Theorem~\ref{the: dimension of higher Tor-s vanish in the amenable case} proves the
``only if''-statement of Conjecture~\ref{con: characterization of amenable by dim(Tor)}.
Some evidence for the ``if''-statement of 
Conjecture~\ref{con: characterization of amenable by dim(Tor)}
comes from the following fact. Notice that a group which contains a non-abelian free group
as a subgroup, cannot be amenable.

Suppose that $G$ contains
a free group $\bbZ \ast \bbZ$ of rank $2$ as a subgroup. Notice that
$S^1 \vee S^1$ is a model for $B(\bbZ \ast \bbZ)$. Its cellular
$\bbC [\bbZ \ast \bbZ]$-chain complex yields an exact sequence
$0 \to \bbC [\bbZ \ast \bbZ]^2 \to \bbC [\bbZ \ast \bbZ] \to \bbC \to 0$,
where $\bbC$ is equipped with the trivial $\bbZ \ast \bbZ$-action.
One easily checks 
$b_1^{(2)}(\widetilde{S^1 \vee S^1}) = - \chi(S^1\vee S^1) = 1$.
This implies
\begin{eqnarray*}
\dim_{\caln(\bbZ \ast \bbZ)}\left(
\Tor_1^{\bbC[\bbZ \ast \bbZ]}(\caln(\bbZ \ast \bbZ),\bbC)\right) & = & 1.
\end{eqnarray*}
We conclude from 
Theorem~\ref{the: Induction and Dimension} \ref{the: Induction and Dimension: flat}
\begin{eqnarray*}
\caln(G) \otimes_{\caln(\bbZ \ast \bbZ)}
\Tor_1^{\bbC[\bbZ \ast \bbZ]}(\caln(\bbZ \ast \bbZ),\bbC) &  = &
\Tor_1^{\bbC G}(\caln(G),\bbC G \otimes_{\bbC[\bbZ \ast \bbZ]}\bbC).
\end{eqnarray*}
Theorem~\ref{the: Induction and Dimension} \ref{the: Induction and Dimension: dim} implies
$$\dim_{\caln(G)}\left(\Tor_1^{\bbC G}(\caln(G),\bbC G \otimes_{\bbC[\bbZ \ast \bbZ]}\bbC)\right)
~ = ~ 1.$$
\em
\end{remark}

One may ask for which groups the von Neumann algebra
$\caln(G)$ is flat as a $\bbC G$-module. This is true if
$G$ is virtually cyclic%
\index{group!virtually cyclic}, i.e.
$G$ is finite or contains $\bbZ$ as a normal subgroup of finite index. 
There is some evidence for the following conjecture
(see \cite[Remark 5.15]{Lueck(1998b)}).

\begin{conjecture}[Flatness of $\caln(G)$ over $\bbC G$]
\label{con: N(G) flat over c G}%
\index{Conjecture!Flatness of $\caln(G)$ over $\bbC G$}
The group von Neumann algebra $\caln(G)$ is flat over
$\bbC G$ if and only if $G$ is virtually cyclic.
\end{conjecture}


\typeout{--------------------   Section 6 --------------------------}

\section{Applications to Group Theory}
\label{sec: Applications to Group Theory}

Recall the Definition~\ref{def: L^2-Betti numbers of G-spaces} of
the $L^2$-Betti numbers of a group $G$ by $b_p^{(2)}(G) : = b_p^{(2)}(EG;\caln(G))$.
In this section
we present tools for and examples of computations of the $L^2$-Betti numbers and
discuss applications to group theory. 
We will explain
in Remark~\ref{rem: Rigidity of reduced $L^2$-cohomology of groups} 
that for a torsion-free group with a model of finite type for $BG$ 
the knowledge of $b_p^{(2)}(G;\caln(G))$ is the same as the
knowledge of the reduced $L^2$-homology $H_p^{(2)}(EG,l^2(G))$, or, equivalently, of
$\bfP H_p^G(EG;\caln(G))$ if $G$ satisfies
satisfies the  Atiyah Conjecture~\ref{con: Atiyah Conjecture} for $(G,1,\bbQ)$.


\subsection{$L^2$-Betti Numbers of Groups}
\label{subsec: L2-Betti Numbers of Groups}

Theorem~\ref{the: properties of gen. b_p^{(2)}} implies:

\begin{theorem}[$L^2$-Betti numbers and Betti numbers of groups]
\label{the: L^2-Betti numbers and Betti numbers of groups}
\index{Theorem!$L^2$-Betti numbers and Betti numbers of groups}
In the sequel we use the conventions $0 \cdot \infty = 0$, $r \cdot \infty = \infty$
for $r \in (0,\infty]$  and $r + \infty = \infty$ for $r \in [0,\infty]$ and put
$|G|^{-1} = 0$ for $|G| = \infty$.
Let $G_1, G_2, \ldots$ be a sequence of non-trivial groups. \begin{enumerate} 

\item \label{the: L^2-Betti numbers and Betti numbers of groups: free amalgamated products}
Free amalgamated products\\[1mm]
For
$r \in \{2, 3, \ldots\} \amalg \{\infty\}$ we get

\begin{eqnarray*}
b_0^{(2)}(\ast_{i=1}^r G_i) & = & 0;
\\
b_1^{(2)}(\ast_{i=1}^r G_i) & = &
\left\{ \begin{array}{lll}
r - 1 + \sum_{i=1}^r \left(b_1^{(2)}(G_i) -
\frac{1}{|G_i|}\right) & & \text{ , if } r < \infty;
\\
\infty & & \text{ , if } r = \infty;
\end{array}\right.
\\
b_p^{(2)}(\ast_{i=1}^r G_i) & = & 
\sum_{i=1}^r b_p^{(2)}(G_i) \hspace{8mm}\mbox{ for } p \ge 2;
\\
b_p(\ast_{i=1}^r G_i) & = &
\sum_{i=1}^r b_p(G_i) \hspace{10mm}\mbox{ for } p \ge 1;
\end{eqnarray*}

\item \label{the: L^2-Betti numbers and Betti numbers of groups: Kuenneth formula}
K\"unneth formula
\begin{eqnarray*}
b_p^{(2)}(G_1 \times G_2) & = &
\sum_{i=0}^p b_i^{(2)}(G_1) \cdot
b_{p-i}^{(2)}(G_2);
\\
b_p(G_1 \times G_2) & = &
\sum_{i=0}^p b_i(G_1) \cdot
b_{p-i}(G_2);
\end{eqnarray*}

\item \label{the: L^2-Betti numbers and Betti numbers of groups: Restriction to subr of fin ind.} 
Restriction to subgroups of finite index\\
For a subgroup $H \subseteq G$ of finite index $[G:H]$ we get
$$b_p^{(2)}(H) ~ = ~ [G:H] \cdot b_p^{(2)}(G);$$

\item \label{the: L^2-Betti numbers and Betti numbers of groups: normal finite sub.} 
Extensions with finite kernel\\[1mm]
Let $1 \to H \to G \to Q \to 1$ be an extension of groups with finite $H$. Then
$$b_p^{(2)}(Q) ~ = ~ |H| \cdot b_p^{(2)}(G);$$

\item \label{the: L^2-Betti numbers and Betti numbers of groups: b_0^{(2)}}
Zero-th $L^2$-Betti number\\[1mm]
We have $b_0^{(2)}(G) = 0$ for $|G| = \infty$ and
$b_0^{(2)}(G) = |G|^{-1}$ for $|G| < \infty$.

\end{enumerate}

\end{theorem}

\begin{example}[Independence of $L^2$-Betti numbers and Betti numbers]
\label{exa: Independence of $L^2$-Betti numbers and Betti numbers}
\em
Given an integer $l \ge 1$ and a sequence $r_1$, $r_2$,
$\ldots$, $r_l$ of non-negative rational numbers, we can construct
a group $G$ such that BG is of finite type and
\begin{eqnarray*}
b_p^{(2)}(G) & = &
\left\{\begin{array}{lll} r_p & & \mbox{ for } 1 \le p \le l;\\
0 & & \mbox{ for } l+1 \le p;
\end{array}
\right.
\\
b_p(G) & = & 0 \hspace{11mm} \mbox{ for } p \ge 1,
\end{eqnarray*}
holds as follows.

For integers $m \ge 0$, $n \ge 1$ and $i \ge 1$ define
$$G_i(m,n) = \bbZ/n \times \left(\ast_{k=1}^{2m+2} \bbZ/2 \right)\times
\left(\prod_{j=1}^{i-1} \ast_{l=1}^4\bbZ/2 \right)$$
One easily checks using 
Theorem~\ref{the: L^2-Betti numbers and Betti numbers of groups}.
\begin{eqnarray*}
b_i^{(2)}(G_i(m,n)) & = & \frac{m}{n};
\\
b_p^{(2)}(G_i(m,n)) & = & 0 \hspace{10mm} \mbox{ for } p \not= i;
\\
b_p(G_i(m,n)) & = & 0 \hspace{10mm} \mbox{ for } p \ge 1.
\end{eqnarray*}
Define the desired group $G$ as follows. For $l = 1$ put $G = G_1(m,n)$ if
$r_1 = m/n$. It remains to treat the case $l \ge 2$.
Choose integers $n\ge 1$ and $k \ge l$ with $r_1 = \frac{k-2}{n}$.
Fix for $i=2,3, \ldots , k$
integers $m_i \ge 0$ and $n_i \ge 1$
such that $\frac{m_i}{n \cdot n_i} = r_i$ holds for $1 \le i \le l$
and $m_i = 0$ holds for $i > l$. Put
$$G ~ = ~ \bbZ/n \times \ast_{i=2}^k G_i(m_i,n_i).$$
One easily checks using 
Theorem~\ref{the: L^2-Betti numbers and Betti numbers of groups}
that $G$ has the prescribed $L^2$-Betti numbers and Betti numbers and
a model for $BG$ of finite type.

On the other hand  we can construct
for any sequence $n_1$, $n_2$, $\ldots$ of non-negative integers
a $CW$-complex $X$ of finite type such that
$b_p(X) = n_p$ and $b_p^{(2)}(\widetilde{X}) = 0$ holds for $p\ge 1$,
namely take
$$X = B(\bbZ/2 \ast \bbZ/2) \times \bigvee_{p = 1}^{\infty}
\left(\bigvee_{i=1}^{n_p}S^p\right).$$

This example shows by considering the $(l+1)$-skeleton
that for a finite connected $CW$-complex $X$ the only general relation
between the $L^2$-Betti numbers $b_p^{(2)}(\widetilde{X})$ of its universal covering 
$\widetilde{X}$ and
the Betti numbers $b_p(X)$ of $X$ is given by the Euler-Poincar\'e
formula (see Theorem~\ref{the: properties of gen. b_p^{(2)}}
\ref{the: properties of gen. b_p^{(2)}: Euler-Poincar'e formula})
$$\sum_{p \ge 0} (-1)^p \cdot b_p^{(2)}(\widetilde{X}) 
~ = ~ 
\chi(X) 
~ = ~
\sum_{p \ge 0} (-1)^p \cdot b_p(X).$$

\em
\end{example}


\subsection{Vanishing of $L^2$-Betti Numbers of Groups}
\label{subsec: Vanishing of L2-Betti Numbers of Groups}

Let $d$ be a non-negative integer
or $d = \infty$. In this subsection we want to investigate the class of groups
\begin{eqnarray}
\calb_{d}%
\indexnotation{calb_{d}}
 & := &
\{G \mid b_p^{(2)}(G) = 0
\mbox{ for } 0 \le p \le d\}. \label{clas calb_d}
\end{eqnarray} 
Notice that $\calb_0$ is the class of infinite groups
by Theorem~\ref{the: L^2-Betti numbers and Betti numbers of groups}
\ref{the: L^2-Betti numbers and Betti numbers of groups: b_0^{(2)}}.

\begin{theorem}
\label{the: properties of bnull}
Let $d$ be a non-negative integer
or $d = \infty$. Then:

\begin{enumerate}
\item \label{the: properties of bnull: amenable groups}
The class $\calb_{\infty}$
contains all infinite amenable groups;

\item \label{the: properties of bnull: normal subgroups}
If $G$ contains a normal subgroup $H$
with $H \in \calb_d$, then
$G \in \calb_d$;

\item \label{the: properties of bnull: directed unions}
If $G$ is the union of a directed system of
subgroups $\{G_i \mid i \in I\}$ such that each
$G_i$ belongs to $\calb_d$, then $G \in \calb_d$;

\item \label{the: properties of bnull: amalgamated products}
Suppose that there are groups
$G_1$ and $G_2$ and group homomorphisms
$\phi_i\colon  G_0 \to G_i$ for
$i=1,2$ such that $\phi_1$ and $\phi_2$ are injective, $G_0$
belongs to $\calb_{d-1}$, $G_1$ and $G_2$
belong to $\calb_d$ and
$G$ is the amalgamated product $G_1 \ast_{G_0} G_2$
with respect to $\phi_1$ and $\phi_2$. Then $G$ belongs
to $\calb_d$;

\item \label{the: properties of bnull: extensions and endomorphisms of
quotient}
Let $1 \to H \to G \to K \to 1$
be an exact sequence of groups such that $b_p^{(2)}(H)$
is finite for all $p \le d$. Suppose that $K$
is  infinite amenable  or suppose that $BK$
has finite $d$-skeleton  and there is an injective endomorphism
$j\colon  K \to K$ whose image has finite index, but is
not equal to $K$. Then
$G \in \calb_{d}$;

\item \label{the: properties of bnull: extensions and non-torsion quotient}
Let $1 \to H \to G \to K \to 1$
be an exact sequence of groups such that
$H \in \calb_{d-1}$, $b_d^{(2)}(H) < \infty$ and
$K$ contains an element of infinite order or  finite subgroups
of arbitrary large order. Then $G \in \calb_d$;

\item \label{the: properties of bnull: Gaboriau's result}
Let $1 \to H \to G \to K \to 1$
be an exact sequence of infinite countable groups such that
$b_1^{(2)}(H) < \infty$. Then $G \in \calb_1$.

\end{enumerate}
\end{theorem}
\proof
\ref{the: properties of bnull: amenable groups}
We get $b_p^{(2)}(G) = 0$ for $p = 0$ from
Theorem~\ref{the: L^2-Betti numbers and Betti numbers of groups}
\ref{the: L^2-Betti numbers and Betti numbers of groups: b_0^{(2)}}.
The case $p \ge 1$ follows from
Theorem~\ref{the: L^2-Betti numbers and homology in the amenable case}
\ref{the: L^2-Betti numbers and homology in the amenable case: L^2-Betti numbers}
since $H_p^{\sing}(EG;\bbC) =  0$ for $p \ge 1$.
\\[1mm]
\ref{the: properties of bnull: normal subgroups}
Apply Theorem~\ref{the: L^2-Betti numbers and fibrations}
\ref{the: L^2-Betti numbers and fibrations: vanishing of all L^2-Betti numbers}
to the fibration $BH \to BG \to B(G/H)$.
\\[1mm]
\ref{the: properties of bnull: directed unions}
The proof is based on a  colimit argument. See \cite[Theorem 7.2 (3)]{Lueck(2002)}.
\\[1mm]
\ref{the: properties of bnull: amalgamated products}
The proof is based on a Mayer-Vietoris argument. See \cite[Theorem 7.2 (4)]{Lueck(2002)}.
\\[1mm]
\ref{the: properties of bnull: extensions and endomorphisms of
quotient}
See \cite[Theorem 7.2 (5)]{Lueck(2002)}.
\\[1mm]
\ref{the: properties of bnull: extensions and non-torsion quotient}
This follows from Theorem~\ref{the: L^2-Betti numbers and fibrations}
\ref{the: L^2-Betti numbers and fibrations: d to d+1}
applied to the fibration $BH \to BG \to B(G/H)$.
\\[1mm]
\ref{the: properties of bnull: Gaboriau's result}
This is proved by
Gaboriau \cite[Theorem 6.8]{Gaboriau(2001)}.
\qed

More information about the vanishing of the first $L^2$-Betti number 
can be found for instance in
\cite{Bekka-Valette(1997)}. Obviously the following is true

\begin{lemma} \label{lem: g in callb_infty implies vanishing of chi}
If $G$ belongs to $\calb_{\infty}$, then $\chi^{(2)}(G) = 0$.
\end{lemma}

\begin{remark}[The Theorem of Cheeger and Gromov]
\label{rem: Cheeger and Gromov}
\em
We rediscover from 
Theorem~\ref{the: properties of bnull}
the result of Cheeger and Gromov  \cite{Cheeger-Gromov(1986)}
that all the $L^2$-Betti numbers of an infinite amenable  group $G$ vanish. 
A detailed comparison of our approach and the approach by Cheeger and 
Gromov to $L^2$-Betti numbers can be found in \cite[Remark 6.76]{Lueck(2002)}.
\em
\end{remark}

\begin{remark}[Advantage of the general definition of $L^2$-Betti numbers]
\label{rem: Advantage of general definition of $L^2$-Betti numbers}
\em
Recall that we have given criterions for $G \in \calb_{\infty}$
in Theorem \ref{the: properties of bnull}. Now it becomes clear why it is worth while
to extend the classical notion of the Euler characteristic $\chi(G) := \chi(BG)$
 for groups $G$ with finite $BG$ to arbitrary groups.
For instance it may very well happen for a group
$G$ with finite $BG$ that $G$ contains a normal group
$H$ which is not even finitely generated and has in particular
no finite model for $BH$ and which belongs to $\calb_{\infty}$
(for instance, $H$ is amenable). Then the classical Euler characteristic is not defined
any more for $H$,  but we can still conclude that the classical Euler characteristic
of $G$ vanishes by Remark~\ref{rem: virtual and L^2-Euler characteristics},
Theorem~\ref{the: properties of bnull}
and Lemma~\ref{lem: g in callb_infty implies vanishing of chi}.
\em
\end{remark}


\subsection{$L^2$-Betti Numbers of Some Specific Groups}
\label{subsec: L2-Betti Numbers of Some Specific Groups}

\begin{example}[Thompson's group] \label{exa: Thompson's group} \em
Next we explain the following
observation about \emph{Thompson's group $F$}.%
\index{group!Thompson's group}
\indexnotation{Thompson's group}
It is the group of orientation preserving dyadic PL-automorphisms
of $[0,1]$, where dyadic means that all slopes are integral powers
of $2$ and the break points are contained in $\bbZ[1/2]$.
It has the presentation
$$F = \langle x_0, x_1, x_2, \ldots ~ \mid ~
x_i^{-1}x_nx_i = x_{n+1} \mbox{ for } i < n\rangle.$$
This group  has some very interesting properties.
Its classifying space $BF$ is of finite type \cite{Brown-Geoghegan(1984)}
but is not homotopy equivalent to a finite dimensional $CW$-complex
since $F$ contains $\bbZ^n$ as a subgroup for all $n \ge 0$
\cite[Proposition 1.8]{Brown-Geoghegan(1984)}.
It is not elementary amenable and does not contain
a subgroup which is free on two generators
\cite{Brin-Squier(1985)},
\cite{Cannon-Floyd-Parry(1996)}.
Hence it is a very interesting question whether $F$ is amenable or not. We conclude
from Theorem~\ref{the: properties of bnull}
\ref{the: properties of bnull: amenable groups} that
a necessary condition for $F$ to be amenable is
that $b^{(2)}_p(F)$ vanishes for all $p \ge 0$.
By \cite[Theorem 7.10]{Lueck(2002)} this condition is satisfied. \em
\end{example}

\begin{example}[Artin groups]\label{exa: Artin groups}
\em Davis and Leary \cite{Davis-Leary(2001z)} compute
for every Artin group $A$ the reduced $L^2$-cohomology and thus the
$L^2$-Betti numbers of the universal covering $\widetilde{S_A}$
of its \emph{Salvetti complex} $S_A$.
The Salvetti complex $S_A$ is a CW-complex which is 
conjectured to be a model for the classifying space $BA$ of $A$.
This conjecture is known to be true in many cases and 
implies that the $L^2$-Betti numbers of $A$ are given by the $L^2$-Betti numbers of
$\widetilde{S_A}$. 
\em
\end{example}

\begin{example}[Right angled Coxeter groups] \label{exa: Right angled Coxeter groups} \em
The $L^2$-homology and the $L^2$-Betti numbers of right angled Coxeter groups are treated 
by Davis and Okun \cite{Davis-Okun(2001)}. More details will be given in
Remark~\ref{rem: Davis-Okun}. 
\em
\end{example}

\begin{example}[Fundamental groups of surfaces and $3$-manifolds] 
\label{exa: Fundamental groups of surfaces and 3-manifolds} \em
Let $G$ be the fundamental group of a compact orientable surface $F_g^d$
of genus $g$ with $d$ boundary components.
Suppose that $G$ is non-trivial which is equivalent to the condition that
$d \ge 1$ or $g \ge 1$. Then $F_g^d$ is a model for $BG$ and we have computed
$b_p^{(2)}(G) = b_p^{(2)}(\widetilde{F_g^d})$ in Subsection~\ref{subsec: Surfaces}.

Let $G$ be the fundamental group of a compact orientable $3$-manifold $M$.
The case $|G| < \infty$ is clear, since then the universal covering is 
homotopy equivalent to a sphere or contractible. So let us assume $|G| = \infty$.
Under the condition that $M$ in non-exceptional, we have computed $b_p^{(2)}(\widetilde{M})$
in Theorem~\ref{the: L^2-Betti numbers of exceptional 3-manifolds}. 
If $M$ is prime, then either $M = S^1 \times S^2$ and $G = \bbZ$ and $b_p^{(2)}(G) = 0$ for all
$p \ge 0$ or $M$ is irreducible, in which case $M$ is aspherical and 
$b_p^{(2)}(G) = b_p^{(2)}(\widetilde{M})$.

Suppose that $M$ is not prime. Then still $b_1^{(2)}(G) = b_1^{(2)}(\widetilde{M})$
by Theorem~\ref{the: properties of gen. b_p^{(2)}}
\ref{the: properties of gen. b_p^{(2)}: homotopy invariance: n-connected}
since the classifying map $M \to BG$ is $2$-connected. Suppose the prime decomposition of $M$ looks like
$M = \#_{i=1}^r M_i$. Then $G = \ast_{i=1}^r G_i$ for $G_i = \pi_1(M_i)$.
We know $b_p^{(2)}(G_i)$  for each $i$ if each $M_i$ is non-exceptional and we get 
$b_p^{(2)}(G) ~ = ~ \sum_{i=1}^r b_p^{(2)}(G_i)$  for $p \ge 2$ from
Theorem~\ref{the: L^2-Betti numbers and Betti numbers of groups}
\ref{the: L^2-Betti numbers and Betti numbers of groups: free amalgamated products}.
\em
\end{example}

\begin{example}[One relator groups]
\label{exa: One relator groups} 
\em
Let 
$G ~ = ~ \langle g_1, g_2, \ldots g_s \mid R \rangle $ 
be a torsion-free one relator group for $s \in \{2,3  \ldots \} \amalg \{\infty\}$ 
and one non-trivial relation $R$.
Then
$$b_p^{(2)}(G) ~ = \left\{
\begin{array}{lll}
0 & & \text{ if } p \not= 1;
\\
s-2 & &  \text{ if } p = 1 \text { and } s < \infty;
\\
\infty & &  \text{ if } p = 1 \text { and } s = \infty.
\end{array}\right.
$$
We only treat the case $s < \infty$, the general case is 
obtained from it by taking the free amalgamated product with a free group.
Because the $2$-dimensional $CW$-complex $X$ associated to the given presentation is
a model for $BG$ (see \cite[chapter III \S\S 9 -11]{Lyndon-Schupp(1977)})
and satisfies $\chi(X) = s-2$, it suffices to prove $b_2^{(2)}(G) = 0$.
We sketch the argument of 
Dicks and Linnell for this claim.
Howie~\cite{Howie(1982)}
has shown that such a group $G$ is
locally indicable and hence left-orderable. A result of 
Linnell~\cite[Theorem 2]{Linnell(1992)}
for left-orderable groups says that an element $\alpha \in
\bbC G$ with $\alpha \not= 0$ is a non-zero-divisor in $\calu(G)$.
This implies that the second differential $c_2^{\calu(G)}$ in the chain complex
$ \calu(G) \otimes_{\bbC G} C_*(EG)$ is injective. Since $\calu(G)$ is flat over $\caln(G)$,
we get from Theorem~\ref{the: prop. ext. dim. for calu}
\begin{multline*}
b_2^{(2)}(G) 
~ = ~  
\dim_{\caln(G)}\left(H_p^G(EG;\caln(G))\right) 
~ = ~ 
\dim_{\calu(G)}\left(H_p^G(EG;\calu(G))\right) 
\\
~ = ~ 
\dim_{\calu(G)}\left(\ker\left(c_2^{\calu(G)}\right)\right) 
~ = ~ 0.
\end{multline*}
Linnell has an extensions of this argument to 
non-torsion-free one-relator groups $G$ with $s \ge 2$ generators.
(The case $s = 1$ is obvious.)
Such a group contains a cyclic subgroup
$\bbZ/k$ such that any finite subgroup is subconjugated to $\bbZ/k$ and then
$$b_p^{(2)}(G) ~ = \left\{
\begin{array}{lll}
0 & & \text{ if } p \not= 1;
\\
s-1-\frac{1}{k} & &  \text{ if } p = 1 \text { and } s < \infty;
\\
\infty & &  \text{ if } p = 1 \text { and } s = \infty.
\end{array}\right.
$$
\em
\end{example}

\begin{example}[Lattices] \label{exa: lattices} \em
Let $L$ be a connected semisimple Lie group  with
finite center such that its Lie algebra has no compact ideal.
Let $G \subseteq L$ be a lattice, i.e. a discrete subgroup of finite covolume.
We want to compute its $L^2$-Betti numbers.
There is a subgroup $G_0 \subseteq G$ of finite index which is torsion-free.
Since $b_p^{(2)}(G) = [G:G_0] \cdot b_p^{(2)}(G_0)$,
it suffices to treat the case $G  = G_0$, i.e. $G \subseteq L$ is a torsion-free lattice.

Let $K \subseteq L$ be a maximal compact subgroup. Put $M = G\backslash L/K$.
Then the space $L/K = \widetilde{M}$
is a symmetric space of non-compact type. We have already mentioned
in Theorem~\ref{the: l^2-inv. of sym. sp.}
that the work of Borel~\cite{Borel(1985)} implies for cocompact $G$ 
that $b_p^{(2)}(G) = b_p^{(2)}(\widetilde{M}) \not= 0$ if and only if $\frk(\widetilde{M}) =  0$
and  $2p = \dim(M)$. This is actually true without the condition ``cocompact'',
because the condition ``finite covolume'' is enough. 

Next we deal with the general case of a connected Lie group $L$.
Let $\Rad(L)$ be its radical.
One can choose a compact normal subgroup $K \subseteq L$
such that $R = \Rad(L) \times K$ is a normal subgroup of $L$ and the quotient
$L_1 = L/R$ is a semisimple Lie group such that its Lie algebra has no compact ideal.
Then $G_1 = L/L\cap R$ is a lattice in $L_1$ and $G \cap R$ is a lattice in $R$.
The group $G \cap R$ is a normal amenable subgroup of
$G$. If $G \cap R$  is infinite, we get $b_p^{(2)}(G) = 0$ for all $p \ge 0$
from Theorem~\ref{the: properties of bnull}.
If $G \cap R$  is finite, we get $b_p^{(2)}(G) = |G \cap R|^{-1} \cdot b_p^{(2)}(G_1)$ 
for all $p \ge 0$ from Theorem~\ref{the: L^2-Betti numbers and Betti numbers of groups} 
\ref{the: L^2-Betti numbers and Betti numbers of groups: normal finite sub.}.
If the center of $L_1$ is infinite, the center of $G_1$ must also be infinite
and hence $b_p^{(2)}(G_1) = 0$ for all $p \ge 0$ by
Theorem~\ref{the: properties of bnull}. Suppose that the center of $L_1$ is finite.
Then we know already how to compute the $L^2$-Betti numbers of $G_1$ from the explanation above.

Given a lattice $G$ in a connected Lie group, 
$b_1^{(2)}(G) > 0$ is true if and only if
$G$ is commensurable with a torsion-free lattice in $PSL_2(\bbR)$, or, equivalently
commensurable with a surface group for genus $\ge 2$ or a finitely generated non-abelian free group
(see  Eckmann \cite{Eckmann(2002a)} or Lott \cite[Theorem 2]{Lott(1999c)}).
\em
\end{example}


\subsection{Deficiency and  $L^2$-Betti Numbers of Groups}
\label{subsec: Deficiency and L2-Betti Numbers of Groups}

Let $G$ be a finitely presented group.
Define its \emph{deficiency}%
\index{deficiency}
$\defi(G)$%
\indexnotation{defi(G)}
to be the maximum $g(P)-r(P)$, where $P$ runs over
all presentations $P$ of $G$ and $g(P)$ is the number of generators and
$r(P)$ is the number of relations of a presentation $P$.

Next we reprove the well-known fact that the maximum
appearing in the definition of the deficiency does exist.

\begin{lemma}
\label{lem: estimate of defiency by L^2-Betti numbers}
Let $G$ be a group with finite presentation
$$P = \langle s_1, s_2, \ldots ,  s_g \mid R_1, R_2, \ldots ,  R_r\rangle$$
Let $\phi\colon  G \to K$ be any group homomorphism. Then
\begin{multline*}
g(P) - r(P) ~ \le ~ 1 - b_0^{(2)}(K \times_{\phi}EG;\caln(K))
+ b_1^{(2)}(K \times_{\phi}EG;\caln(K))
\\
- b_2^{(2)}(K \times_{\phi}EG;\caln(K)).
\end{multline*}
\end{lemma}
\proof
Given a presentation $P$ with $g$ generators and $r$ relations,
let $X$ be the associated finite $2$-dimensional $CW$-complex. It has one $0$-cell,
$g$ $1$-cells, one  for each generator, and $r$ $2$-cells, one for each
relation.  There is an obvious isomorphism
from $\pi_1(X)$ to $G$ so that we can choose a map
$f\colon  X \to BG$ which induces an isomorphism on the fundamental groups.
It induces a $2$-connected $K$-equivariant map
$\overline{f}\colon  K \times_{\phi} \widetilde{X} \to
K \times_{\phi} \widetilde{EG}$.
We conclude from
Theorem~\ref{the: properties of gen. b_p^{(2)}}
\ref{the: properties of gen. b_p^{(2)}: homotopy invariance: n-connected}
\begin{eqnarray*}
b_p^{(2)}(K \times_{\phi}\widetilde{X};\caln(K))
& = &
b_p^{(2)}(K \times_{\phi}EG;\caln(K))
\hspace{5mm} \mbox{ for } p=0,1;
\\
b_2^{(2)}(K \times_{\phi}\widetilde{X};\caln(K))
& \ge &
b_2^{(2)}(K \times_{\phi}EG;\caln(K)).
\end{eqnarray*}
We conclude from the $L^2$-Euler-Poincar\'e formula (see
Theorem \ref{the: basic properties of the l2-Euler characteristic}
\ref{the: basic properties of the l2-Euler characteristic: m and h})
\begin{eqnarray*}
g - r
& = & 1 - \chi^{(2)}(K \times_{\phi} \widetilde{X};\caln(K))
\\
& = &
1 -
b_0^{(2)}(K \times_{\phi}\widetilde{X};\caln(K))
+ b_1^{(2)}(K \times_{\phi}\widetilde{X};\caln(K))\nonumber
\\ & & \hspace{47mm}
- b_2^{(2)}(K \times_{\phi}\widetilde{X};\caln(K))
\\
& \le &
1 -
b_0^{(2)}(K \times_{\phi} EG;\caln(K))
+ b_1^{(2)}(K \times_{\phi}EG;\caln(K))\nonumber
\\ & & \hspace{43mm}
- b_2^{(2)}(K \times_{\phi}EG;\caln(K)). \qed
\end{eqnarray*}

\begin{example}[Deficiency of some groups] \label{exa: Deficiency of some groups} \em
Sometimes the deficiency is realized by the ``obvious'' presentation.
For instance the deficiency of a free group $\langle s_1, s_2, \ldots , s_g \mid \emptyset\rangle$
on $g$ letters is indeed $g$.  The cyclic group 
$\bbZ/n$ of order $n$ has the presentation $\langle t \mid t^n = 1\rangle$ and its deficiency is $0$.
The group $\bbZ/n \times \bbZ/n$ has the presentation $\langle s,t \mid s^n, t^n, [s,t] \rangle$
and its deficiency is $-1$.
\em
\end{example}

\begin{remark}[Non-additivity of the deficiency]
\label{exa: Non-additivity of the deficiency}
\em 
The deficiency is not additive under free products
by the following example which is a special case of a more general example
due to Hog, Lustig and Metzler
\cite[Theorem 3 on page 162]{Hog-Lustig-Metzler(1985)}.
The group
\mbox{$(\bbZ/2 \times \bbZ/2) \ast (\bbZ/3 \times \bbZ/3)$}
has the obvious presentation
\begin{multline*}
\langle s_0,t_0,s_1,t_1 \mid s_0^2 = t_0^2 = [s_0,t_0] =
s_1^3 = t_1^3 = [s_1,t_1] = 1 \rangle
\end{multline*}
One may think that its deficiency is $-2$. However,
it turns out that its deficiency is $-1$, realized by the following presentation
\begin{multline*}
\langle s_0,t_0,s_1,t_1 \mid s_0^2 = 1, [s_0,t_0] =
t_0^2, s_1^3 = 1, [s_1,t_1] = t_1^3, t_0^2 = t_1^3 \rangle.
\end{multline*}
This shows that it is important to get upper bounds on the deficiency of groups.
Writing down presentations gives lower bounds, but it is not clear whether 
a given presentation realizes the deficiency.
\em
\end{remark}

\begin{lemma}
\label{lem: vanishing of L^2-Betti numbers and deficiency and signature}
Let $G$ be a finitely presented group and let $\phi\colon  G \to K$ be
a homomorphism such that $b_1^{(2)}(K\times_{\phi}EG;\caln(K)) = 0$. Then
\begin{enumerate}
\item
\label{lem: vanishing of L^2-Betti numbers and deficiency and signature: defi}
$\defi(G) \le 1$;

\item
\label{lem: vanishing of L^2-Betti numbers and deficiency and signature: sign}
Let $M$ be a closed oriented $4$-manifold with $G$ as
fundamental group. Then
$$ |\sign (M)| ~ \le ~ \chi(M).$$
\end{enumerate}
\end{lemma}
\proof
\ref{lem: vanishing of L^2-Betti numbers and deficiency and signature: defi} 
This follows directly from Lemma~\ref{lem: estimate of defiency by L^2-Betti numbers}.
\\[1mm]
\ref{lem: vanishing of L^2-Betti numbers and deficiency and signature: sign}
This is a consequence of the $L^2$-Signature Theorem due to Atiyah~\cite{Atiyah(1976)}.
Details of  the proof can be found in \cite[Lemma 7.22]{Lueck(2002)}.\qed

\begin{theorem} \label{the: extensions and deficiency and signature}
Let $1 \to H \xrightarrow{i} G \xrightarrow{q}
K \to 1$ be an exact sequence of infinite groups.
Suppose that $G$ is
finitely presented  and one of the following conditions is satisfied:

\begin{enumerate}
\item
\label{the: extensions and deficiency and signature: condition 1}
$b_1^{(2)}(H) < \infty$;

\item \label{the: extensions and deficiency and signature: condition 2}
The classical first Betti number of $H$ satisfies
$b_1(H)< \infty$  and $K$ belongs to $\calb_1$.

\end{enumerate}

Then
\renewcommand{\labelenumi}{(\roman{enumi})}
\begin{enumerate}

\item $\defi(G) \le 1$;

\item Let $M$ be a closed oriented $4$-manifold with $G$ as
fundamental group. Then
$$ |\sign (M)| ~ \le ~ \chi(M).$$

\end{enumerate}
\renewcommand{\labelenumi}{(\arabic{enumi})}
\end{theorem}
\proof
If condition
\ref{the: extensions and deficiency and signature: condition 1}
is satisfied, then $b_p^{(2)}(G) =  0$  for $p = 0,1$ by Theorem
\ref{the: properties of bnull}
\ref{the: properties of bnull: Gaboriau's result},
and the claim follows
from Lemma \ref{lem: vanishing of L^2-Betti numbers and deficiency and
signature}.
\par

Suppose that condition
\ref{the: extensions and deficiency and signature: condition 2}
is satisfied.
There is a spectral sequence converging to
$H_{p+q}^{K}(K \times_q EG;\caln(K))$
with $E^2$-term
$$E^2_{p,q} ~ = ~ \Tor^{\bbC K}_p(H_q(BH;\bbC),\caln(K))$$
\cite[Theorem 5.6.4 on page 143]{Weibel(1994)}.
Since $H_q(BH;\bbC)$ is $\bbC$ with the trivial
$K$-action for $q=0$ and
finite dimensional as complex vector space
by assumption for $q = 1$, we conclude
$\dim_{\caln(K)}(E^2_{p,q}) = 0$ for $p + q = 1$
from the assumption $b_1^{(2)}(K) = 0$.
This implies $b_1^{(2)}(K \times_q EG;\caln(K)) = 0$
and the claim follows from
Lemma \ref{lem: vanishing of L^2-Betti numbers and deficiency and signature}.
\qed

Theorem \ref{the: extensions and deficiency and signature}
generalizes results in  \cite{Eckmann(1997)},
\cite{Johnson-Kotschick(1993)},
where also some additional information is given. Furthermore see 
\cite{Hausmann-Weinberger(1985)}, \cite{Kotschick(1994)}. We mention the result of Hitchin
\cite{Hitchin(1974a)} that a connected closed oriented smooth $4$-manifold
which admits an Einstein metric satisfies the stronger inequality
$|\sign(M)| \le \frac{2}{3} \cdot \chi(M)$.

Finally we mention the following result of Lott
\cite[Theorem 2]{Lott(1999c)} (see also \cite{Eckmann(2002b)}) 
which generalizes a result of Lubotzky
\cite{Lubotzky(1983)}. The statement we present here is a slight improvement
of Lott's result due to Hillman \cite{Hillman(1999)}.

\begin{theorem}[Lattices of positive deficiency]
\label{the: Lott's result on deficiency}
\index{Theorem!Lattices of positive deficiency}
Let $L$ be a connected Lie group. Let $G$ be a lattice in $L$. If
$\defi(G) > 0$, then one of the following assertions holds:
\begin{enumerate}

\item  $G$  is a lattice in $PSL_2(\bbC )$;

\item $\defi(G) = 1$. Moreover, either
$G$ is isomorphic to a torsion-free non-uniform
lattice in $\bbR \times PSL_2(\bbR)$ or $PSL_2(\bbC)$, or
$G$ is $\bbZ$ or $\bbZ^2$.
\end{enumerate}
\end{theorem}


\typeout{--------------------   Section 7 --------------------------}

\section{$G$- and $K$-Theory}
\label{sec: $G$- and $K$-Theory}

In this section we discuss the projective class group $K_0(\caln(G))$ 
of a group von Neumann algebra. We present applications of
its computation to $G_0(\bbC G)$ and the Whitehead group $\Wh(G)$ of a group $G$.


\subsection{The $K_0$-group of a Group von Neumann Algebra}
\label{subsec: The K_0-group of a Group von Neumann Algebra}

In this subsection we want to investigate the projective class group
of a group von Neumann algebra.

\begin{definition}[Definition of  $K_0(R)$ and $G_0(R)$]
\label{def: K_0(R) and G_0(R)}  Let $R$ be an (associative) ring
(with unit).
Define its \emph{projective class group}%
\index{projective class group}
$K_0(R)$%
\indexnotation{K_0(R)}
to be the abelian group whose generators are isomorphism classes $[P]$
of finitely generated projective $R$-modules $P$ and whose relations are
$[P_0] + [P_2] = [P_1]$ for any exact sequence
$0 \to P_0 \to P_1 \to P_2 \to 0$
of finitely generated projective $R$-modules. 

Define the
\emph{Grothendieck group of finitely generated modules}
\index{Grothendieck group of finitely generated modules}
$G_0(R)$%
\indexnotation{G_0(R)}
analogously but replace finitely generated projective with finitely generated.
\end{definition}

The group $K_0$ is known for any von Neumann algebra (see for instance 
\cite[Subsection 9.2.1]{Lueck(2002)}. For simplicity we only treat
 the von Neumann algebra $\caln(G)$ of a group here. 

The next result is taken from
\cite[Theorem 7.1.12 on page 462,
Proposition 7.4.5 on page 483, Theorem 8.2.8 on page 517, Proposition 8.3.10
on page 525, Theorem 8.4.3 on page 532]{Kadison-Ringrose(1986)}.

\begin{theorem}[The universal trace] \label{the: universal trace}
\index{Theorem!The universal trace}
There is a map
$$\tr^u_{\caln(G)}%
\indexnotation{tr^u_caln(G)}
 \colon  \caln(G) \to \calz(\caln(G))$$
into the center $\calz(\caln(G))$ of $\caln(G)$
called the  \emph{center valued trace}%
\index{trace!center valued}
 or \emph{universal trace}%
\index{trace!universal}
of $\caln(G)$, which is uniquely determined by the following two properties:

\begin{enumerate}

\item $\tr^u_{\caln(G)}$ is a trace with values in the center,
i.e. $\tr^u_{\caln(G)}$ is $\bbC$-linear, for $a \in \caln(G)$ with $a \ge 0$ we have
$\tr^u_{\caln(G)}(a) \ge 0$ and $\tr^u_{\caln(G)}(ab) = \tr^u_{\caln(G)}(ba)$ for all $a,b \in \caln(G)$;

\item $\tr^u_{\caln(G)}(a) = a$ for all $a \in \calz(\caln(G))$.
\\[4mm]
\hspace*{-10mm}

The map $\tr^u_{\caln(G)}$ has the following further properties:
\vspace*{2mm}

\item $\tr^u_{\caln(G)}$ is faithful, i.e. $\tr^u_{\caln(G)}(a) = 0 \Leftrightarrow a = 0$ for $a \in \caln(G), a
  \ge 0$;

\item $\tr^u_{\caln(G)}$ is normal, i.e. 
for a monotone increasing net $\{a_i \mid i \in I\}$ of positive elements $a_i$ with
supremum $a$ we have $\tr^u_{\caln(G)}(a) = \sup\{\trû(a_i) \mid i \in \}$,
or, equivalently,
$\tr^u_{\caln(G)}$ is continuous with respect to the
ultra-weak topology on $\caln(G)$;

\item $||\tr^u_{\caln(G)}(a)||\;\le\; ||a||$ for $a \in \caln(G)$;

\item $\tr^u_{\caln(G)}(ab) = a \tr^u_{\caln(G)}(b)$ for all
$a \in \calz(\caln(G))$ and $b \in \caln(G)$;

\item Let $p$ and $q$ be projections in $\caln(G)$. Then $p\sim q$, 
i.e. $p = uu^*$ and $q = u^*u$ for some element $u \in \caln(G)$, if
and only if $\tr^u_{\caln(G)}(p) = \tr^u_{\caln(G)}(q)$;

\item Any linear functional
$f \colon  \caln(G) \to \bbC$which is continuous with respect to the norm
topology on $\caln(G)$ and which is central, i.e. $f(ab)= f(ba)$
for all $a,b \in \caln(G) $ factorizes as
$$\caln(G) \xrightarrow{\tr^u_{\caln(G)}} \calz(\caln(G)) \xrightarrow{f|_{\calz(\caln(G))}} \bbC.$$

\end{enumerate}
\end{theorem}

\begin{definition}[Center valued dimension] \label{def: dim^u}
For a finitely generated projective $\caln(G)$-module $P$ define its
\emph{center valued von Neumann dimension}
\index{von Neumann dimension!center valued}
by
$$
\dim^u_{\caln(G)}(P)%
\indexnotation{dim^u_caln(G)}
:= \sum_{i=1}^n \tr^u_{\caln(G)}(a_{i,i})
\in \calz(\caln(G))^{\bbZ/2} = \{a \in \calz(\caln(G)) \mid a = a^*\}$$
for any matrix $A= (a_{i,j})_{i,j} \in M_n(\caln(G))$ with
$A^2 = A$ such that $\im(r_A\colon  \caln(G)^n \to \caln(G)^n)$
induced by right multiplication with $A$ is $\caln(G)$-isomorphic to
$P$.
\end{definition}

There is a classification of von Neumann algebras into certain types. We only need to know
what the type of a group von Neumann algebra is. 

\begin{lemma} \label{lem: type and center of N(G)}
Let $G$ be a discrete group. Let $G_f$%
\indexnotation{G_f} be the normal subgroup
of $G$ consisting of elements $g \in G$ whose centralizer
has finite index (or, equivalently,
whose conjugacy class (g) consists of
finitely many elements). Then:

\begin{enumerate}

\item \label{lem: type and center of N(G): type  I}
The group von Neumann algebra
$\caln(G)$ is of type $I_f$ if and only if $G$ is virtually abelian;

\item \label{lem: type and center of N(G): type  II}
The group von Neumann algebra $\caln(G)$ is of type $II_1$ if and only if
the index of $G_f$ in $G$ is infinite;

\item \label{lem: type and center of N(G): type for finitely generated G}
Suppose that $G$ is finitely generated. Then $\caln(G)$ is of type $I_f$
if $G$ is virtually abelian, and of type $II_1$ if $G$ is
not virtually abelian;

\item \label{lem: type and center of N(G): factor}
The group von Neumann algebra $\caln(G)$ is a \emph{factor,}%
\index{factor}
i.e. its center consists of $\{r \cdot 1_{\caln(G)} \mid r \in \bbC\}$,
if and only if $G_f$ is the trivial group.

\end{enumerate}
\end{lemma}
\proof
\ref{lem: type and center of N(G): type  I} This is proved in
\cite{Kaniuth(1969)}, \cite{Thoma(1964)}.
\\[1mm]
\ref{lem: type and center of N(G): type  II} This is proved in
\cite{Kaniuth(1969)},\cite{Mautner(1950)}.
\\[1mm]
\ref{lem: type and center of N(G): type for finitely generated G}
This follows from assertions \ref{lem: type and center of N(G): type I} and
\ref{lem: type and center of N(G): type II} since for finitely
generated $G$ the group $G_f$ has finite index in $G$ if and only if $G$ is
virtually abelian.
\\[1mm]
\ref{lem: type and center of N(G): factor}
This follows from \cite[Proposition 5 in III.7.6 on page 319]{Dixmier(1981)}.
\qed

The next result follows from
\cite[Theorem 8.4.3 on page 532, Theorem 8.4.4 on page 533]
{Kadison-Ringrose(1986)}.

\begin{theorem}[$K_0$ of finite von Neumann algebras]
\index{Theorem!K0 of finite von Neumann algebras@$K_0$ of finite von
  Neumann algebras}
\label{the: K_0 of fin. vN-alg.}
Let $G$ be a group. 

\begin{enumerate}
\item \label{the: K_0 of fin. vN-alg.: dim^u and [P] = [Q]}
The following statements are equivalent for two
finitely generated projective $\caln(G)$-modules $P$ and $Q$:

\begin{enumerate}
\item $P$ and $Q$ are $\caln(G)$-isomorphic;

\item $P$ and $Q$ are stably $\caln(G)$-isomorphic, i.e.
$P \oplus V$ and $Q \oplus V$ are $\caln(G)$-isomorphic
for some finitely generated projective
$\caln(G)$-module $V$;

\item $\dim^u_{\caln(G)}(P) = \dim^u_{\caln(G)}(Q)$;

\item $[P] = [Q]$ in $K_0(\caln(G))$;
\end{enumerate}

\item \label{the: K_0 of fin. vN-alg.: dim^u and K_0(cala)}
The center valued dimension induces an injection
$$\dim^u_{\caln(G)}\colon  K_0(\caln(G)) \to \calz(\caln(G))^{\bbZ/2} =
\{a \in \calz(\caln(G)) ~ \mid ~
a = a^*\},$$
where the group structure on $\calz(\caln(G))^{\bbZ/2}$ comes from addition.
If $\caln(G)$ is of type $II_1$, this map is an isomorphism.
\qed
\end{enumerate}
\end{theorem}
 
\begin{remark}[Group von Neumann algebras and representation theory]
\label{Group von Neumann algebras and representation theory}
\em
Theorem~\ref{the: K_0 of fin. vN-alg.} shows that the group von
Neumann algebra is the right generalization of the complex group ring
from finite groups to infinite groups if one is concerned with
representation theory of finite groups.  Namely, let $G$ be a finite group. Recall that a
finite dimensional complex $G$-representation $V$ is the same as a
finitely generated $\bbC G$-module and that $K_0(\bbC G)$ is the same
as the complex representation ring.
Moreover, two finite dimensional $G$-representations $V$ and $W$ are
linearly $G$-isomorphic if and only if they have the same character.
Recall that the character is a class function. One easily checks that
the complex vector space of class functions on a finite group $G$ is
the same as the center $\calz(\bbC G)$ and that the character of $V$ is
the same as $\dim^u_{\caln(G)}(V)$. 
\em
\end{remark}

\begin{remark}[Factors]
\label{rem: Factors}
\em
Suppose that $\caln(G)$ is a factor, i.e. 
its center consists of $\{r \cdot 1_{\caln(G)} \mid r \in \bbC\}$.
By Lemma~\ref{lem: type and center of N(G)}
\ref{lem: type and center of N(G): factor}
this is the case if and only if 
$G_f$ is the trivial group. Then $\dim_{\caln(G)} = \dim^u_{\caln(G)}$  and two 
finitely generated projective $\caln(G)$-modules $P$ and $Q$ are $\caln(G)$-isomorphic if and only
if $\dim_{\caln(G)}(P) = \dim_{\caln(G)}(Q)$ holds. This has the consequence
that for a free $G$-$CW$-complex $X$ of finite type the $p$-th $L^2$-Betti number
determines the isomorphism type of $\bfP H_p^G(X;\caln(G))$. In particular
we must have $\bfP H_p^G(X;\caln(G)) \cong_{\caln(G)} \caln(G)^n$
provided that 
$n = b_p^{(2)}(X;\caln(G))$ is an integer. If one prefers to work with 
reduced $L^2$-homology, this is equivalent to the statement that
$H_p^{(2)}(X;l^2(G))$ is isometrically $G$-linearly isomorphic to
$l^2(G)^n$ provided that 
$n = b_p^{(2)}(X;\caln(G))$ is an integer.
\em
\end{remark}

\begin{remark}[The reduced $L^2$-cohomology of torsion-free groups]
\label{rem: Rigidity of reduced $L^2$-cohomology of groups}
\em
Let $G$ be a torsion-free group. Suppose that it satisfies the 
Atiyah Conjecture~\ref{con: Atiyah Conjecture} for $(G,1,\bbQ)$.
Suppose that there is a model for $BG$ of finite type. 
Then we get for all $p$ that 
$\bfP H_p^G(EG;\caln(G)) ~\cong_{\caln(G)} \caln(G)^n$, or, equivalently, that
$H_p^{(2)}(X;l^2(G))$ is isometrically $G$-linearly isomorphic to $l^2(G)^n$
if the integer $n$ is given by $n = b_p^{(2)}(X;\caln(G))$.
This claim is proved in  \cite[solution to Exercise 10.11 on page 546]{Lueck(2002)}.
\em
\end{remark}

We mention that the inclusion $i \colon \caln(G) \to \calu(G)$ induces an isomorphism
$$K_0(\caln(G)) \xrightarrow{\cong} K_0(\calu(G)).$$ 

The Farrell-Jones Conjecture for $K_0(\bbC G)$, the Bass Conjecture and the passage in $K_0$ from
$\bbZ G$ to $\bbC G$ and to $\caln(G)$ is discussed in \cite[Section 9.5.2]{Lueck(2002)}
and \cite{Lueck-Reich(2003)}. 



\subsection{The $K_1$-group and the $L$-groups of a Group von Neumann Algebra}
\label{subsec: The $K_1$-group and the $L$-groups of a Group von Neumann Algebra}

A complete calculation of the $K_1$-group and of the $L$-groups of any von Neumann algebra
and of the associated algebra of affiliated operators 
can be found in \cite[Section 9.3 and Section 9.4]{Lueck(2002)},
\cite{Lueck-Roerdam(1993)} and
\cite{Reich(2001)}.


\subsection{Applications to $G$-theory of Group Rings}
\label{subsec: Applications to $G$-theory of Group Rings}

\begin{theorem}[Detecting $G_0(\bbC G)$ by $K_0(\caln(G))$ for amenable groups]
\index{Theorem!Detecting $G_0(\bbC G)$ by $K_0(\caln(G))$ for amenable groups}
\label{the: G_0(C) to K_0(N(G))}
If $G$ is amenable, the map
$$l\colon  G_0(\bbC G) \to K_0(\caln(G)),
\hspace*{10mm} [M] \mapsto [\bfP \caln(G) \otimes_{\bbC G} M]$$
is a well-defined homomorphism. If $f\colon  K_0(\bbC G) \to G_0(\bbC G)$
is the forgetful map sending $[P]$ to $[P]$ and $i_*\colon  K_0(\bbC G) \to
K_0(\caln(G))$ is induced by the inclusion $i\colon  \bbC G \to \caln(G)$, then
the composition $l \circ f$ agrees with $i_*$.
\end{theorem}
\proof This is essentially a consequence of 
the dimension-flatness of $\caln(G)$ over $\bbC G$
(see Theorem \ref{the: dimension of higher Tor-s vanish in the amenable case}).
Details of the proof can be found in 
\cite[Theorem 9.64]{Lueck(2002)}. \qed

Now one can combine Theorem~\ref{the: K_0 of fin. vN-alg.} and
Theorem~\ref{the: G_0(C) to K_0(N(G))}  to detect elements in $G_0(\bbC G)$ for amenable
$G$. In particular one can show  
\begin{eqnarray}
\dim_{\bbQ}\left(\bbQ \otimes_{\bbZ} G_0(\bbC G)\right) & \ge & |\con(G)_{f,cf}|,
\label{estimate for rk_{bbZ}(G_0(bbC G))}
\end{eqnarray}
where $\con(G)_{f,cf}$ is the the set of conjugacy classes $(g)$ of elements $g \in G$
such that $g$ has finite order and $(g)$ contains only finitely many elements.
Notice that $\con(G)_{f,cf}$ contains at least one element, namely the unit element
$e$. 

\begin{remark}[The non-vanishing of {$[RG]$} in $G_0(RG)$ for amenable groups] 
\label{rem: The non-vanishing of [RG] in G_0(RG) for amenable groups}
\em
A direct consequence of Theorem~\ref{the: G_0(C) to K_0(N(G))}
is that for an amenable group $G$ the class $[\bbC G]$ in $G_0(\bbC G)$
generates an infinite cyclic subgroup. Namely, the dimension induces a well-defined
homomorphism
$$\dim_{\caln(G)} \colon G_0(\bbC G) \to \bbR, \quad 
[M] ~ \mapsto ~\dim_{\caln(G)}\left(\caln(G) \otimes_{\bbC G} M\right),$$
which sends $[\bbC G]$ to $1$.
This result has been extended by Elek \cite{Elek(2002)} 
to finitely generated amenable groups 
and arbitrary fields $F$, i.e. there is a well-defined homomorphism
$G_0(FG) \to \bbR$, which sends $[FG]$ to $1$ and is given by a certain rank function
on finitely generated $FG$-modules.
\em
\end{remark}

The class $[RG]$ in $K_0(RG)$ is never zero for a commutative integral domain $R$
with quotient field $R_{(0)}$.
The augmentation $RG \to R$ and the map 
$K_0(R) \to \bbZ, \quad [P] \mapsto \dim_{R_{(0)}}(R_{(0)} \otimes_R P)$ 
together induce a homomorphism
$K_0(RG) \to \bbZ$ which sends $[RG]$ to $1$.
A decisive difference between $K_0(RG)$ and $G_0(RG)$ is that 
$[RG] = 0$ is possible in $G_0(RG)$ as the following example shows. 

\begin{example}[The vanishing of {$[RG]$} in $G_0(RG)$ for groups G containing
  $\bbZ \ast \bbZ$] 
\label{exa: The Vanishing of [RG] in G_0(RG) for groups G containing bbZ ast bbZ}
\em
We abbreviate $F_2 = \bbZ \ast \bbZ$. Suppose that $G$ contains $F_2$ as a subgroup.
Let $R$ be a ring. Then
$$[RG] ~ = ~ 0 \quad \in G_0(RG)$$
holds by the following argument.
Induction with the inclusion $F_2 \to G$
induces a homomorphism
$G_0(RF_2) \to G_0(RG)$
which sends $[RF_2]$ to $[RG]$. Hence it suffices
to show $[RF_2] = 0$ in $G_0(RF_2)$.
The cellular chain complex of the universal covering
of $S^1\vee S^1$ yields an exact sequence of $RF_2$-modules
$0 \to (RF_2)^2  \to RF_2 \to R
\to 0$, where $R$ is equipped with the trivial $F_2$-action.
This implies $[RF_2] = - [R]$ in $G_0(RF_2)$.
Hence it suffices to show $[R] = 0$ in $G_0(RF_2)$.
Choose an epimorphism $f\colon  F_2 \to \bbZ$.
Restriction with $f$ defines a homomorphism
$G_0(R\bbZ) \to G_0(RF_2)$. It
sends the class of $R$ viewed as trivial $R\bbZ$-module to the class of $R$ viewed as
trivial $RF_2$-module. Hence it remains to show
$[R] = 0$ in $G_0(R\bbZ)$.
This follows from the exact sequence
$0 \to R\bbZ \xrightarrow{s-1} R\bbZ
\to R \to 0$ for $s$ a generator of $\bbZ$
which comes from the cellular $R\bbZ$-chain complex of $\widetilde{S^1}$.
\em
\end{example}

Remark~\ref{rem: The non-vanishing of [RG] in G_0(RG) for amenable groups} and
Example~\ref{exa: The Vanishing of [RG] in G_0(RG) for groups G containing bbZ ast bbZ} 
give some evidence for

\begin{conjecture}{\bf (Amenability and the regular
representation in $G$-theory).}
\index{Conjecture!Amenability and the regular representation in $G$-theory}
\label{con: [RG] in G_0(RG)}
Let $R$ be a commutative integral domain. Then a group $G$ is amenable if and only if
$[RG] \not= 0$ in $G_0(RG)$.
\end{conjecture}

\begin{remark}[The Atiyah Conjecture for amenable groups and $G_0(\bbC G)$]
\label{rem: The Atiyah Conjecture for amenable groups and G_0(bbC G)} \em
Assume that $G$ is amenable and that there
is an upper bound on the orders of finite subgroups
of $G$. Then the 
Atiyah Conjecture~\ref{con: Atiyah Conjecture} for $(G,d,\bbC)$ is true 
if and only if the image of the map
$$\dim_{\caln(G)}\colon G_0(\bbC G) \to \bbR, \hspace{5mm} [M] \mapsto
\dim_{\caln(G)}(\caln(G) \otimes_{\bbC G} M)$$
is contained in  $\{r \in \bbR \mid d \cdot r \in \bbZ\}$.
\em
\end{remark}

\begin{example}[$K_0(\bbC G) \to G_0(\bbC G)$ is not necessarily surjective]
\label{exa: oplus_n in Z ZZ/2}\em
Let $A= \bigoplus_{n \in \bbZ} \bbZ/2$. This abelian group is locally finite. Hence
the map 
$$\bigoplus_{H \subseteq A, |H| < \infty}~ K_0(\bbC H) \to K_0(\bbC A)$$ 
is surjective and the image of
$$\dim_{\caln(G)} \colon K_0(\bbC A) \to \bbR, \quad 
[P] \mapsto \dim_{\caln(A)}\left(\caln(A) \otimes_{\bbC A} P\right)$$
is $\bbZ[1/2]$. On the other hand the argument in
\cite[Example 10.13]{Lueck(2002)} shows that the map
$$\dim_{\caln(G)} \colon G_0(\bbC A) \to \bbR, \quad 
[M] \mapsto \dim_{\caln(A)}\left(\caln(A) \otimes_{\bbC A} M\right)$$
is surjective. In particular the obvious map
$K_0(\bbC A) \to G_0(\bbC A)$ is not surjective.
\em
\end{example}


\subsection{Applications to the Whitehead Group}
\label{subsec: Applications to the Whitehead Group}

The \emph{Whitehead group}%
\index{Whitehead group} 
$\Wh(G)$%
\indexnotation{Wh(G)} of a group $G$ is the quotient of $K_1(\bbZ G)$ by the subgroup
which consists of elements given by units of the shape $\pm g \in \bbZ G$ 
for $g \in G$. Let $i\colon  H \to G$ be the inclusion of a normal subgroup $H \subseteq G$.
It induces a homomorphism $i_0\colon  \Wh(H) \to \Wh(G)$.
The conjugation action of $G$ on $H$ and on $G$ induces
a $G$-action on $\Wh(H)$ and on $\Wh(G)$ which turns out to be trivial
on $\Wh(G)$. Hence $i_0$ induces  homomorphisms
\begin{eqnarray}
i_1\colon  \bbZ \otimes_{\bbZ G} \Wh(H) & \to & \Wh(G);
\label{the: Wh(H) to Wh(G): i_1}
\\
i_2\colon  \Wh(H)^G  & \to & \Wh(G).
\label{the: Wh(H) to Wh(G): i_2}
\end{eqnarray}

\begin{theorem}[Detecting elements in $\Wh(G)$]
\index{Theorem!Detecting elements in $\Wh(G)$}
\label{the: Wh(H) to Wh(G)}
Let $i\colon  H \to G$ be the inclusion of a normal finite subgroup $H$
into an arbitrary group $G$. Then the maps
$i_1$ and $i_2$ defined in
\eqref{the: Wh(H) to Wh(G): i_1}
and
\eqref{the: Wh(H) to Wh(G): i_2}
have finite kernel.
\end{theorem}
\proof See \cite[Theorem 9.38]{Lueck(2002)}. \qed

We emphasize that Theorem~\ref{the: Wh(H) to Wh(G)} above  holds for all groups $G$.
It seems to be related to the Farrell-Jones Isomorphism Conjecture.


\typeout{--------------------   Section 8 --------------------------}

\section{$L^2$-Betti Numbers and Measurable Group Theory}
\label{sec: L^2-Betti Numbers and Measurable Group Theory}

In this section we want to discuss an interesting relation between $L^2$-Betti numbers and
measurable group theory. We begin with formulating the main result.

\begin{definition}[Measure equivalence] \label{def: Measure Equivalence}
Two countable groups $G$ and
$H$ are called
\emph{measure equivalent}%
\index{group!measure equivalent}
\index{measure equivalent groups}
if there exist commuting measure-preserving free
actions of $G$ and $H$ on some standard Borel space $(\Omega,\mu)$
with non-zero Borel measure $\mu$
such that the actions  of both $G$ and $H$ admit  measure fundamental
domains $X$ and $Y$ of finite measure.

The triple $(\Omega,X,Y)$ is called a \emph{measure coupling}%
\index{measure coupling} of $G$ and $H$. The index%
\index{index!of a coupling triple} 
of $(\Omega,X,Y)$ is the quotient
$\frac{\mu(X)}{\mu(Y)}$.
\end{definition}

Here are some explanations.
A \emph{Polish space}%
\index{Polish Space}
\index{space!Polish}
is a separable topological space which is metrizable by a complete metric.
A \emph{measurable space}%
\index{measurable space}
\index{space!measurable}
$\Omega = (\Omega,\cala)$ is a set $\Omega$ together with a $\sigma$-algebra $\cala$. 
It is called a \emph{standard Borel space}%
\index{standard Borel space}
\index{space!standard Borel}
if it is isomorphic to a Polish space with its Borel $\sigma$-algebra. (The Polish space
is not part of the structure, only its existence is required.)
More information about this notion of measure equivalence can be found for instance in
\cite{Furman(1999a)}, \cite{Furman(1999b)} and \cite[0.5E]{Gromov(1993)}.

The following result is due to Gaboriau
\cite[Theorem 6.3]{Gaboriau(2001)}. We will discuss its applications and
sketch the proof based on homological algebra and the dimension function
 due to R. Sauer \cite{Sauer(2002)}.

\begin{theorem}[Measure equivalence and $L^2$-Betti numbers]%
\index{Theorem!Measure equivalence and L2-Betti numbers@Measure equivalence and $L^2$-Betti numbers}
\label{the: Gaboriau's result on measure equivalent groups}
Let $G$ and $H$ be two countable groups which are measure equivalent.
If $C>0$ is the index of a measure coupling, then we get for all $p \ge 0$
$$b_p^{(2)}(G) ~ = ~ C \cdot b_p^{(2)}(H).$$
\end{theorem}

The general strategy of the proof of
Theorem~\ref{the: Gaboriau's result on measure equivalent groups} is as follows.
In the first step one introduces the notion of a standard action $G \action X$ and 
of a weak orbit equivalence of
standard actions of index $C$ and shows
that two groups $G$ and $H$ are measure equivalent of index $C$ if and only if
there exist standard actions $G \action X$ and $H \action Y$ which are weakly 
orbit equivalent with index $C$. In the second step one assigns to
a standard action $G\action X$ $L^2$-Betti numbers $b_p^{(2)}(G \action X)$,
which involve only data that is invariant under orbit equivalence.
Hence $b_p^{(2)}(G \action X)$ itself depends only on the orbit equivalence class
of $G \action X$. In order to deal with weak orbit equivalences, one has to investigate
the behaviour of the $L^2$-Betti numbers of $b_p^{(2)}(G\action X)$ under restriction.
Finally one proves that the $L^2$-Betti numbers of a standard action $G \action X$
agree with the $L^2$-Betti numbers of $G$ itself.

A version of Theorem~\ref{the: Gaboriau's result on measure equivalent groups}
for the $L^2$-torsion is presented in Conjecture~\ref{con: Measure equivalence and L^2-torsion}.



\subsection{Measure Equivalence and Quasi-Isometry}
\label{subsec: Measure Equivalence and Quasi-Isometry}

\begin{remark}[Measure equivalence is the measure theoretic version of quasi-isometry]
\em 
The notion of measure equivalence can be viewed as
the measure theoretic analogue of the metric notion of quasi-isometric
groups. Namely, two finitely generated groups $G_0$ and $G_1$ are
\emph{quasi-isometric}%
\index{group!quasi-isometric}
\index{quasi-isometric groups}
 if and only if there exist commuting proper (continuous)
actions of $G_0$ and $G_1$ on some locally compact space such that
each action has a cocompact fundamental domain
\cite[0.2 $C_2'$ on page 6]{Gromov(1993)}. 
\em
\end{remark}

\begin{example}[Infinite amenable groups] \label{exa: Infinite amenable groups} \em
Every countable infinite amenable group is measure equivalent to $\bbZ$
(see \cite{Ornstein-Weiss(1980)}).
Since obviously all the $L^2$-Betti numbers of $\bbZ$ vanish, 
Theorem~\ref{the: Gaboriau's result on measure equivalent groups} implies the result
of  Cheeger and Gromov that all the $L^2$-Betti numbers of an infinite amenable group vanish.
\em
\end{example}

\begin{remark}[$L^2$-Betti numbers and quasi-isometry] 
\label{rem: ^L2-Betti numbers and quasi-isometry}
\em
If the finitely generated groups $G_0$ and $G_1$
are quasi-isometric and there exist finite models for $BG_0$ and $BG_1$, 
then $b_p^{(2)}(G_0) = 0 \Leftrightarrow
b_p^{(2)}(G_1)= 0$ holds (see \cite[page 224]{Gromov(1993)}, \cite{Pansu(1995)}).
But in general it is not true
that there is a constant $C > 0$ such that
$b_p^{(2)}(G_0) ~ = ~ C \cdot b_p^{(2)}(G_1)$ holds for all $p \ge 0$
(cf.  \cite[page 7]{Gaboriau(2001b)}, \cite[page 233]{Gromov(1993)},
\cite{Whyte(1999)}).
\em
\end{remark}

\begin{remark}[Measure equivalence versus quasi-isometry] 
\label{rem: Measure equivalence versus quasi-isometry}\em
If $F_g$ denotes the free group on $g$ generators, then
define $G_n := (F_3 \times F_3) \ast F_n$ for $n \ge 2$.
The groups $G_m$ and $G_n$ are
quasi-isometric for $m,n \ge 2$ (see \cite[page 105 in IV-B.46]{delaHarpe(2000)},
\cite[Theorem 1.5]{Whyte(1999)})
and have finite models for their classifying spaces.
One easily checks using 
Theorem~\ref{the: L^2-Betti numbers and Betti numbers of groups}
that $b_1^{(2)}(G_n) = n$ and $b_2^{(2)}(G_n) = 4$.

Theorem \ref{the: Gaboriau's result on measure equivalent groups} due to Gaboriau
implies that $G_m$ and $G_n$ are measure equivalent if and only if $m = n$ holds.
Hence there are finitely presented groups which are quasi-isometric
but not measure equivalent.

The converse is also true. The groups $\bbZ^n$ and $\bbZ^m$ are infinite
amenable and hence measure equivalent. But they are not quasi-isometric for different $m$ and $n$
since $n$ is the growth rate of
$\bbZ^n$ and the growth rate is a quasi-isometry invariant.

Notice that
Theorem \ref{the: Gaboriau's result on measure equivalent groups} implies that the
sign of the Euler characteristic of a group $G$ is an invariant under measure equivalence,
which is not true for quasi-isometry by the example of the groups $G_n$ above.

Let the two groups $G$ and $H$ act on the same metric space $X$ properly and cocompactly by isometries.
If $X$ is second countable and proper, then $G$ and $H$ are measure equivalent.
\cite[Theorem~2.36]{Sauer(2002)}.
If $X$ is a geodesic and proper, then  $G$ and $H$ are quasi-isometric.
\em
\end{remark}

\begin{remark}[Kazhdan's property (T)] \label{rem: Kazhdan's property (T)} \em
Kazhdan's property $(T)$ is an invariant under
measure equivalence \cite[Theorem 8.2]{Furman(1999a)}.
There exist quasi-isometric finitely generated groups
$G_0$ and $G_1$ such that $G_0$ has Kazhdan's property $(T)$ and
$G_1$ not (see \cite[page 7]{Gaboriau(2001b)}). Hence $G_0$ and $G_1$
are quasi-isometric but not measure equivalent.
\em
\end{remark}

The rest of this section is devoted to an outline of the proof
of Theorem~\ref{the: Gaboriau's result on measure equivalent groups}
 due to R. Sauer \cite{Sauer(2002)} which is simpler and more algebraic 
than the orginal one of Gaboriau \cite{Gaboriau(2001)}
and may have the potential to apply also to $L^2$-torsion.


\subsection{Discrete Measured Groupoids}
\label{subsec: Discrete Measured Groupoids}

A \emph{groupoid}%
\index{groupoid}
is a small category in which all morphisms are isomorphisms. We will identify
a groupoid $\underline{G}$ with its set of morphisms. Then the set of objects
$\underline{G}^0$ can be considered as a subset of $\underline{G}$ via the identity morphisms.
There are four canonical maps, 
$$\begin{array}{lllll}
\text{source map} & s \colon \underline{G} & \to & \underline{G}^0, &
(f \colon x  \to y) ~ \mapsto ~ x;
\\
\text{target map} & t \colon \underline{G} & \to & \underline{G}^0, &
(f \colon x  \to y) ~ \mapsto ~ y;
\\
\text{inverse map} & i \colon \underline{G} & \to & \underline{G}, &
f ~ \mapsto ~ f^{-1};
\\
\text{composition} & \circ \colon \underline{G}^2 & \to & \underline{G}, & (f,g) \mapsto
f \circ g,
\end{array}
$$
where $ \underline{G}^2$ is $\{(f,g) \in \underline{G} \times \underline{G} \mid s(f) =
t(g)\}$. We will often abbreviate $f \circ g$ by $fg$.

A \emph{discrete measurable groupoid}%
\index{discrete measurable groupoid}
\index{groupoid!discrete measurable} 
is a groupoid $\underline{G}$ equipped with the structure of a standard Borel space such
that the inverse map and the composition are measurable maps and $s^{-1}(x)$ is countable for
all objects $x \in \underline{G}^0$. Then $\underline{G}^0 \subseteq \underline{G}$ is a
Borel subset, the source and the target maps are
measurable and $t^{-1}(x)$ is countable for all objects $x \in \underline{G}^0$.

Let $\mu$ be a probability measure on $\underline{G}^0$. Then for each measurable subset
$A \subseteq \underline{G}$ the function
$$\underline{G}^0 \to \bbC, \quad x \mapsto |s^{-1}(x) \cap A|$$
is measurable and we obtain a $\sigma$-finite measure $\mu_s$ on $\underline{G}$ by
$$\mu_s(A) ~ := ~ \int_{\underline{G}^0} |s^{-1}(x) \cap A| ~ d\mu(x).$$
It is called the left counting measure of $\mu$. The right counting measure
$\mu_t$ is defined analogously replacing the source map $s$ by the target map $t$.
We call $\mu$ invariant if $\mu_s = \mu_t$, or, equivalently, if
$i_*\mu_s = \mu_s$. A discrete measurable groupoid $\underline{G}$ together with
an invariant measure $\mu$ on $\underline{G}^0$ is called a 
\emph{discrete measured groupoid}.%
\index{discrete measured groupoid}
\index{groupoid!discrete measured}
Given  a Borel subset $A \subseteq \underline{G}^0$ with $\mu(A) > 0$,
there is the \emph{restricted discrete measured groupoid}%
\index{discrete measured groupoid!restricted}
$\underline{G}|_A = s^{-1}(A) \cap t^{-1}(A)$, which is equipped with the 
\emph{normalized measure}%
\index{normalized measure}
$\frac{1}{\mu(A)}\cdot \mu|_A$. 

An \emph{isomorphism of discrete measured groupoids}
$f \colon \underline{G} \to \underline{H}$ is an isomorphisms
of groupoids which preserves the measures. Given measurable subsets 
$A \subseteq \underline{G}^0$ and
$B \subseteq \underline{H}^0$ such that $t(s^{-1}(A))$ and $t(s^{-1}(B))$ have full
measure in $\underline{G}^0$ and $\underline{H}^0$ respectively,  
we call an isomorphism of discrete measured groupoids
$f \colon \underline{G}_A \to \underline{H}_B$ a 
\emph{weak isomorphism of discrete measured groupoids}.

\begin{example}[Orbit equivalence relation] \label{exa: orbit equivalence relation} 
\em
Consider the countable group $G$ with an action $G\action X$ on a standard Borel space $X$ with probability
measure $\mu$ by $\mu$-preserving isomorphisms. The \emph{orbit equivalence relation}
\index{orbit equivalence relation}
$$\calr(G \action X)%
\indexnotation{calr(G action X)}
 ~ := ~ \{(x,gx) \mid x \in X, g \in G\} ~ \subseteq X \times X$$
becomes a discrete measured groupoid by the obvious groupoid structure and measure.
\em
\end{example}

An action $G \action X$ of a countable group $G$ is called \emph{standard}%
\index{standard action}
\index{action!standard}
if $X$ is a standard Borel measure space with a probability measure $\mu$, the action is by $\mu$-preserving
Borel isomorphisms and the action is \emph{essentially free},%
\index{essentially free action}
\index{action!essentially free}
i.e. the stabilizer of almost every $x \in X$ is trivial. Every countable group $G$
admits a standard action, which is given by the shift action on 
$\prod_{g \in G} [0,1]$. Notice that this $G$-action is not free but essentially free.

Two standard actions $G\action X$ and $H \action Y$ are \emph{weakly orbit equivalent}%
\index{weakly orbit equivalent standard actions}
\index{standard action!weakly orbit equivalent}
if there are Borel subsets $A \subseteq X$ and $B \subseteq Y$, which meet almost every
orbit and have positive measure in $X$ and $Y$ respectively, 
and a Borel isomorphism $f \colon A \to B$, which preserves the normalized
measures on $A$ and $B$ and satisfies
$$f(G \cdot x \cap A) ~ = ~ H \cdot f(x) \cap B$$
for almost all $x \in A$. If $A$ has full measure in $X$ and $B$ has full measure in $Y$,
then the two standard actions are called \emph{orbit equivalent}.%
\index{orbit equivalent standard actions}
\index{standard action!orbit equivalent}
The map $f$ is called a \emph{weak orbit equivalence}%
\index{standard action!weak orbit equivalence of}
\index{weak orbit equivalence of standard actions}
 or \emph{orbit equivalence}%
\index{standard action!orbit equivalence of}
\index{orbit equivalence of standard actions}
respectively. The 
\emph{index of a weak orbit equivalence}%
\index{index!of a weak orbit equivalence}
of $f$ is the quotient $\frac{\mu(A)}{\mu(B)}$. The next result is due to
Furman \cite[Theorem 3.3]{Furman(1999b)}.

\begin{theorem}[Measure equivalence and weak orbit equivalence]
\label{the: Measure equivalence and weak orbit equivalence}
Two countable groups are measure equivalent with respect to a measure coupling  of index $C>0$
if and only if there exist standard actions of $G$ and $H$ which are weakly orbit
equivalent with index $C$.
\end{theorem}


\subsection{Groupoid Rings}
\label{subsec: Groupoid Rings}

Let $\underline{G}$ be a discrete measured groupoid with invariant
measure $\mu$ on $\underline{G}^0$. For a function $\phi \colon \underline{G} \to \bbC$ and 
$x \in \underline{G}^0$ put
\begin{eqnarray*}
S(\phi)(x) & := & |\{g \in \underline{G} \mid \phi(g) \not= 0, s(g) = x\} \quad 
\in \{0,1,2 \ldots\} \amalg \{\infty\};
\\
T(\phi)(x) & := & |\{g \in \underline{G} \mid \phi(g) \not= 0, t(g) = x\} \quad 
\in \{0,1,2 \ldots\} \amalg \{\infty\}.
\end{eqnarray*}
Let $\mu_{\underline{G}} = \mu_s = \mu_t$ be the measure on $\underline{G}$ induced by $\mu$.
Let $L^{\infty}(\underline{G}) = L^{\infty}(\underline{G};\mu_{\underline{G}})$ 
be the $\bbC$-algebra of equivalence classes of essentially bounded
measurable functions $\underline{G} \to \bbC$. Define
$L^{\infty}(\underline{G}^0) = L^{\infty}(\underline{G}^0;\mu)$  analogously.
  Define the \emph{groupoid ring}%
\index{groupoid ring}
\index{ring!groupoid ring} of $\underline{G}$ as the subset 
\begin{eqnarray}
\hspace{-8mm} \bbC \underline{G}%
\indexnotation{C underline G}
 & := & \{\phi \in L^{\infty}(\underline{G}) \mid
S(\phi) \text{ and } T(\phi) \text{ are essentially bounded on } \underline{G}\}.
\label{def of Cunderline G}
\end{eqnarray}
The addition comes from the  pointwise addition in $ L^{\infty}(\underline{G}) $.
Multiplication comes from the convolution product
$$(\phi \cdot \psi)(g) ~ = ~
\sum_{\substack{g_1,g_2 \in \underline{G}\\g_2 \circ g_1 = g}} \phi(g_1) \cdot
  \psi(g_2).$$
An involution of rings on $\bbC \underline{G}$ is defined by $(\phi^*)(g) := \overline{\phi(i(g))}$. 
Define the augmentation homomorphism $\epsilon \colon \bbC \underline{G} \to 
L^{\infty}(\underline{G}^0)$ by sending $\phi$ to 
$\epsilon(\phi) \colon \underline{G}^0 \to \bbC, 
\quad x \mapsto \sum_{g \in s^{-1}(x)} \phi(g)$.
Notice that $\epsilon$ is in general not a ring homomorphism, it is only compatible
with the additive structure. It becomes a homomorphism of $\bbC \underline{G}$-modules
if we equip $L^{\infty}(\underline{G}^0)$ with the following $\bbC \underline{G}$-module
structure
$$\phi \cdot f := \epsilon(\phi \cdot j(f)) \quad \text{ for } \phi \in \bbC \underline{G},
f \in L^{\infty}(\underline{G}^0),$$
where $j \colon L^{\infty}(\underline{G}^0) \to \bbC \underline{G}$ is the inclusion of
rings, which is given by extending a function on $\underline{G}^0$ to $\underline{G}$ 
by putting it to be zero outside $\underline{G}^0$. 

Given a group $G$ and a ring $R$ together with a homomorphism $c \colon G \to \aut(R)$,
define the \emph{crossed product ring}%
\index{ring!crossed product}
\index{crossed product ring}
$R \ast_c G$%
\indexnotation{R ast_c G}
as the free $R$-module with $G$ as $R$-basis and the multiplication given by
$$\left(\sum_{g \in G} r_g \cdot g\right) \cdot \left(\sum_{g \in G} s_g \cdot g\right)
~ = ~ \sum_{g \in G} \left(\sum_{\substack{g_1,g_2 \in G,\\g = g_1g_2}} r_{g_1}\cdot c(g_1)(s_{g_2})\right)
\cdot g.$$
Given a standard action $G \action X$, let $L^{\infty}(X) \ast G$ be the crossed product
ring $L^{\infty}(X) \ast_c G$ with respect to the group homomorphism
$c \colon G \to \aut\left(L^{\infty}(X)\right)$ sending $g$ to the automorphism
given by composition with $l_{g^{-1}} \colon X \to X, \quad x \mapsto g^{-1} x$.
We obtain an injective ring homomorphism
$$k \colon L^{\infty}(X) \ast G \to \bbC\calr(G \action X)$$
which sends $\sum_{g \in G}f_g \cdot g$ to the function
$(gx,x) ~ \mapsto ~ f_g(gx)$.
In the sequel we will regard $L^{\infty}(X) \ast G$
as a subring of $\bbC\calr(G \action X)$
using $k$.

Next we briefly explain how one can associate to the groupoid ring
$\bbC \underline{G}$ of a discrete measured groupoid $\underline{G}$ a
von Neumann algebra $\caln(\underline{G})$, which is finite, or, equivalently, which possesses 
a faithful finite normal trace. One can define on $\bbC \underline{G}$ an inner product
\begin{eqnarray*}
\scal{\phi}{\psi}  ~ = ~ 
\int_{\underline{G}} \phi(g) \cdot \overline{\psi(g)} ~ d\mu_{\underline{G}}.
\end{eqnarray*}
Then $\bbC \underline{G}$ as a $\bbC$-algebra with involution and the scalar product above
satisfies the axioms of a \emph{Hilbert algebra}%
\index{Hilbert algebra} $A$, i.e. we have $\scal{y}{x} = \scal{x^*}{y^*}$
for $x,y \in A$, $\scal{xy}{z} = \scal{y}{x^*z}$ for $x,y,z \in A$ and 
the map $A \to A, ~ y \mapsto yx$ is continuous for all $x \in A$. Let $H_A$ be the Hilbert
space completion of $A$ with respect to the given inner product. Define the von 
Neumann algebra $\caln(A)$ associated to $A$ by the $\bbC$-algebra with involution $\calb(H_A)^A$
which consists of all bounded left $A$-invariant operators $H_A \to H_A$. The 
standard trace is given by
$$\tr_{\caln(A)} \colon \caln(A) \to \bbC, \quad f ~ \mapsto ~ \scal{f(1_A)}{1_A}.$$
 We do get a dimension function 
as in Theorem~\ref{the: prop. ext. dim.} for $\caln(A)$.

Our main example will be $\caln(G \action X) := \caln(\bbC \calr(G\action X))$ for a standard
action  $G \action X$ of $G$.

If $G$ is a countable group and $\underline{G} = G$ is the associated discrete measured groupoid with one
object, then $\bbC \underline{G} = \bbC G$, $l^2(G) = H_{\bbC \underline{G}}$ and
the definition of  $\caln(\underline{G})$ and $\tr_{\caln(\underline{G})}$ above
agrees with the previous Definition~\ref{def: group von Neumann algebra}
of $\caln(G)$ and $\tr_{\caln(G)}$.

\begin{remark}[Summary and Relevance of the algebraic structures associated to a standard action] 
\label{rem: Summary and Relevance of the algebraic structures associated to a standard action}
\em
Let $G\action X$ be a standard action. 
We have the following commutative diagram of inclusions of rings
$$
\begin{CD}
\bbC @>>> \bbC G @> = >> \bbC G @>>> \caln(G)
\\
@VVV @VVV @VVV @VVV 
\\
L^{\infty}(X) @>>> L^{\infty}(X) \ast G @>>> \bbC\calr(G \action X) @>>> \caln(G \action X)
\end{CD}
$$
There is a $\bbC \underline{G}$-module structure on $L^{\infty}(\calr(G \action X)^0) = L^{\infty}(X)$.
Its restriction to $L^{\infty}(X) \ast G ~ \subseteq ~ \bbC\calr(G \action X)$ 
is the obvious $L^{\infty}(X) \ast G$-module structure on $L^{\infty}(X)$.

The following observation will be crucial. Given two standard actions
$G \action X$ and $G \action Y$, an orbit equivalence $f$ from $G \action X$ to
$H \action Y$ induces isomorphisms of rings, all denoted by $f_*$, such that the following 
diagram with inclusions as horizontal maps commutes
$$
\begin{CD}
L^{\infty}(X) @>>> \bbC\calr(G \action X) @>>> \caln(G \action X)
\\
@V f_* V \cong V @V f_* V \cong V @V f_* V \cong V 
\\
L^{\infty}(Y) @>>> \bbC\calr(H \action Y) @>>> \caln(G \action Y)
\end{CD}
$$
It is \emph{not} true that $f$ induces a ring map 
$L^{\infty}(X) \ast G \to L^{\infty}(Y) \ast H$, since we only require that
$f$ maps orbits to orbits but nothing is demanded about equivariance of $f$ with 
respect to some homomorphism of groups from $G \to H$. The crossed
product ring $L^{\infty}(X) \ast G$  contains too much information about the group
$G$ itself. Hence we shall only involve
$L^{\infty}(X)$, $\bbC\calr(G \action X)$, and $\caln(G \action X)$ in any algebraic
construction which is designed to be invariant under orbit equivalence.
\em
\end{remark}


\subsection{$L^2$-Betti Numbers of Standard Actions}
\label{subsec: L^2-Betti Numbers of Standard Actions}

\begin{definition} \label{def: L^2-Betti numbers of measaured groupoids}
Let $\underline{G}$ be a discrete measured groupoid. Define its
\emph{$p$-th $L^2$-Betti number}%
\index{L2-Betti number@$L^2$-Betti number!of a discrete measured groupoid}
by
$$b_p^{(2)}(\underline{G})%
\indexnotation{b_p^(2)(underline G)}
 ~ = ~ \dim_{\caln(\underline{G})}\left(
\Tor_p^{\bbC \underline{G}}\left(\caln(\underline{G}),L^{\infty}(\underline{G}^0)\right)\right).$$

Given a standard  action $G \action X$, define its 
\emph{$p$-th $L^2$-Betti number}%
\index{L2-Betti number@$L^2$-Betti number!of a standard action}
as the $p$-th $L^2$-Betti number of the associated orbit equivalence relation
$\calr(G \action X)$, i.e. 
$$
b_p^{(2)}(G \action X)%
\indexnotation{b_p^(2)(G action X)}
~ = ~ \dim_{\caln(G\action X)}\left(
\Tor_p^{\bbC \calr(G\action X)}\left(\caln(G \action X),L^{\infty}(X)\right)\right).
$$
\end{definition}

Notice that Theorem~\ref{the: Gaboriau's result on measure equivalent groups}
is true if we can prove the following three lemmas.

\begin{lemma} \label{lem: Invar. of L^2-Betti numb. of standard actions under orbit equiv.}
If two standard actions $G \action X$ and $H \action Y$ are orbit equivalent, then
they have the same $L^2$-Betti numbers.
\end{lemma}

\begin{lemma} \label{lem. restriction of groupoids} 
Let $\underline{G}$ be a discrete measured groupoid. Let $A \subseteq \underline{G}$ be
a Borel subset such that $t(s^{-1}(A))$ has full measure in $\underline{G}^0$. 
Then we get for all $p \ge 0$
$$b_p^{(2)}(\underline{G}) ~ = ~ \mu(A) \cdot b_p^{(2)}(\underline{G}|_A).$$
\end{lemma}

\begin{lemma} \label{lem: b-P^(2)(G) = b_p^(2)(G action X)}
Let $G \action X$ be a standard action. Then we get for all $p \ge 0$
$$b_p^{(2)}(G \action X) ~ = ~ b_p^{(2)}(G).$$
\end{lemma}

Lemma ~\ref{lem: Invar. of L^2-Betti numb. of standard actions under orbit equiv.} 
follows directly from 
Remark~\ref{rem: Summary and Relevance of the algebraic structures associated to a standard action}.
The hard part of the proof of 
Theorem~\ref{the: Gaboriau's result on measure equivalent groups} is indeed
the proof of the remaining two Lemmas~\ref{lem. restriction of groupoids}  and
\ref{lem: b-P^(2)(G) = b_p^(2)(G action X)}. This is essentially done by developing
some homological algebra over finite von Neumann algebras taking the dimension for
arbitrary modules into account.


\subsection{Invariance of $L^2$-Betti Numbers under Orbit Equivalence}
\label{subsec: Invariance of $L^2$-Betti Numbers under Orbit Equivalence}

As an illustration we sketch the proof of 
Lemma~\ref{lem: b-P^(2)(G) = b_p^(2)(G action X)}. It follows from the following chain
of equalities which we explain briefly below.
\begin{eqnarray}
\hspace{-9mm} b_p^{(2)}(G) & = & \dim_{\caln(G)}\left(\Tor_p^{\bbC G}\left(\caln(G),\bbC\right)\right)
\label{eq 1}
\\
& = & 
\dim_{\caln(G\action X)}\left(
\caln(G\action X) \otimes_{\caln(G)} \Tor_p^{\bbC G}\left(\caln(G),\bbC\right)\right)
\label{eq 2}
\\
& = & 
\dim_{\caln(G\action X)}\left(\Tor_p^{\bbC G}\left(\caln(G\action X),\bbC\right)\right)
\label{eq 3}
\\
& = & 
\dim_{\caln(G\action X)}\left(\Tor_p^{L^{\infty}(X) \ast G}
\left(\caln(G\action X),L^{\infty}(X) \ast G \otimes_{\bbC G} \bbC\right)\right)
\label{eq 4}
\\
& = & 
\dim_{\caln(G\action X)}\left(\Tor_p^{L^{\infty}(X) \ast G}
\left(\caln(G\action X),L^{\infty}(X)\right)\right)
\label{eq 5}
\\
& = & 
\dim_{\caln(G\action X)}\left(\Tor_p^{\bbC\calr(G\action X)}
\left(\caln(G\action X),L^{\infty}(X)\right)\right)
\label{eq 6}
\\
& = & b_p^{(2)}(G \action X).
\label{eq 7}
\end{eqnarray}
Equations \eqref{eq 1} and \eqref{eq 7} are true by definition. The inclusion of von
Neumann algebras $\caln(G) \to \caln(G\action X)$ preserves the traces. This implies
that the functor $\caln(G \action X) \otimes_{\caln(G)} -$ from $\caln(G)$-modules to
$\caln(G \action X)$-modules is faithfully flat and preserves dimensions.
The proof of this fact is completely analogous to the proof of
Theorem~\ref{the: Induction and Dimension}. This shows \eqref{eq 2} and \eqref{eq 3}.
For every $\bbC G$-module $M$ there is a natural $L^{\infty}(X) \ast G$-isomorphism
$$L^{\infty}(X) \ast G \otimes_{\bbC G} M \xrightarrow{\cong} L^{\infty}(X)
\otimes_{\bbC} M.$$
This shows that $L^{\infty}(X) \ast G$ is flat as $\bbC G$-module and that
\eqref{eq 4} and \eqref{eq 5} are true. The hard part is now to prove 
\eqref{eq 6}, which is the decisive step, since here one eliminates
$L^{\infty}(X) \ast G$ from the picture and stays with terms which depend
only on the orbit equivalence class of $G \action X$. Its proof
involves homological algebra and dimension theory. It is not true
that the relevant Tor-terms are isomorphic, they only have the same dimension.

This finishes the outline of
the proof of Lemma~\ref{lem: b-P^(2)(G) = b_p^(2)(G action X)} and
of  Theorem~\ref{the: Gaboriau's result on measure equivalent groups}. The complete proof
can be found in \cite{Sauer(2002)}.


\typeout{--------------------   Section 9 --------------------------}

\section{The Singer Conjecture}
\label{sec: The Singer Conjecture}

In this section we briefly discuss the following conjecture.

\begin{conjecture}[Singer Conjecture]%
\index{Conjecture!Singer Conjecture}
\label{con: Singer Conjecture}
If $M$ is an aspherical closed manifold, then
$$b_p^{(2)}(\widetilde{M}) ~ = ~ 0 \hspace{10mm}
\mbox{if } 2p \not= \dim(M).$$
If $M$ is a closed connected Riemannian manifold with negative sectional
curvature, then
$$b_p^{(2)}(\widetilde{M}) ~  \left\{
\begin{array}{lll}
= 0 & & \mbox{if } 2p \not= \dim(M);\\
> 0 & & \mbox{if } 2p = \dim(M).
\end{array}\right.$$
\end{conjecture}

We mention that all the explicit computations presented in
Section~\ref{sec: Computations of L^2-Betti Numbers} 
are compatible with the Singer Conjecture~\ref{con: Singer Conjecture}.
A version of the Singer Conjecture for $L^2$-torsion will be presented
in Conjecture~\ref{con: $L^2$-torsion for aspherical manifolds}.


\subsection{The Singer Conjecture and the Hopf Conjecture}
\label{subsec: The Singer Conjecture and the Hopf Conjecture}

Because of the Euler-Poincar\'e formula $\chi(M) = \sum_{p \ge 0}
(-1)^p \cdot b_p^{(2)}(\widetilde{M})$ (see 
Theorem~\ref{the: properties of gen. b_p^{(2)}}
\ref{the: properties of gen. b_p^{(2)}: Euler-Poincar'e formula})
the Singer Conjecture~\ref{con: Singer Conjecture} implies the following conjecture in
case $M$ is aspherical or has negative sectional curvature.

\begin{conjecture}[Hopf Conjecture]%
\index{Conjecture!Hopf Conjecture}
\label{con: Hopf Conjecture}
If  $M$ is an aspherical closed manifold of even dimension, then
$$(-1)^{\dim(M)/2} \cdot \chi(M) \ge 0.$$
If  $M$ is a  closed Riemannian manifold of even dimension with
sectional curvature $\sec(M)$, then
$$\begin{array}{rlllllll}
(-1)^{\dim(M)/2} \cdot \chi(M) & > & 0 & &
\mbox{ if } & \sec(M) & < & 0;
\\
(-1)^{\dim(M)/2} \cdot \chi(M) & \ge   & 0 & &
\mbox{ if } & \sec(M) & \le & 0;
\\
\chi(M) & = & 0 & &
\mbox{ if } & \sec(M) & = & 0;
\\
\chi(M) & \ge & 0 & &
\mbox{ if } & \sec(M) & \ge & 0;
\\
\chi(M) & > & 0 & &
\mbox{ if } & \sec(M) & > & 0.
\end{array}$$
\end{conjecture}

In original  versions of the Singer Conjecture \ref{con: Singer Conjecture}
and the Hopf Conjecture \ref{con: Hopf Conjecture} the statements
for aspherical manifolds did not appear.
Every Riemannian manifold with non-positive sectional
curvature is aspherical by Hadamard's Theorem.


\subsection{Pinching Conditions}
\label{subsec: Pinching Conditions}

The following two results are taken from the paper by Jost and Xin
\cite[Theorem 2.1 and Theorem 2.3]{Jost-Xin(2000)}.

\begin{theorem} \label{the: Jost-Xin-Ricci}
Let $M$ be a closed connected Riemannian manifold of dimension $\dim(M) \ge
3$. Suppose that there are real numbers $a > 0$ and $b > 0$ such that the
sectional curvature satisfies $-a^2 \le \sec(M) \le 0$ and the Ricci
curvature is bounded from above by $-b^2$.   If the non-negative
integer $p$ satisfies $2p \not= \dim(M)$ and $2pa \le b$, then
$$b_p^{(2)}(\widetilde{M}) ~ = ~ 0.$$
\end{theorem}

\begin{theorem} \label{the: Jost-Xin-sec}
Let $M$ be a closed connected Riemannian manifold of dimension $\dim(M) \ge 4$.
Suppose that there are real numbers $a > 0$ and $b > 0$ such that the
sectional curvature satisfies $-a^2 \le \sec(M) \le -b^2$. If the non-negative
integer $p$ satisfies $2p \not= \dim(M)$ and
$(2p-1)\cdot a \le  (\dim(M)-2)\cdot b$, then
$$b_p^{(2)}(\widetilde{M}) ~ = ~ 0.$$
\end{theorem}

The next result is a consequence of a result of Ballmann and Br\"uning
\cite[Theorem B on page 594]{Ballmann-Bruening(2001)}.

\begin{theorem} \label{the: Ballmann-Bruening}
Let $M$ be a closed connected Riemannian manifold.
Suppose that there are real numbers $a > 0$ and $b > 0$ such that the
sectional curvature satisfies $-a^2 \le \sec(M) \le -b^2$. If the non-negative
integer $p$ satisfies $2p < \dim(M)-1$ and
$p \cdot a < (\dim(M) -1 -p)\cdot b$, then
\begin{eqnarray*}
b_p^{(2)}(\widetilde{M}) & = & 0.
\end{eqnarray*}
\end{theorem}

Theorem \ref{the: Jost-Xin-sec} and Theorem \ref{the: Ballmann-Bruening}
are improvements of the
older results by Donnelly and Xavier \cite{Donnelly-Xavier(1984)}.

\begin{remark}[Right angled Coxeter groups and Coxeter complexes]  \label{rem: Davis-Okun} \em
Next we mention the work of Davis and Okun \cite{Davis-Okun(2001)}.
A simplicial complex $L$ is called a \emph{flag complex}%
\index{flag complex}
if each finite non-empty set of vertices
which pairwise are connected by edges spans a simplex of $L$. To such a flag complex
they  associate a right-angled Coxeter group $W_L$ defined by the following presentation
\cite[Definition 5.1]{Davis-Okun(2001)}.
Generators are the vertices $v$ of $L$. Each generator $v$ satisfies $v^2 = 1$. If
two vertices $v$ and $w$ span an edge, there is
the relation $(vw)^2 = 1$. Given a finite flag complex
$L$, Davis and Okun associate to it a finite proper $W_L$-$CW$-complex $\Sigma_L$,
which turns out to be a model for the
classifying space of the family of finite subgroups $\EGF{W_L}{\calfin}$
\cite[6.1, 6.1.1  and 6.1.2]{Davis-Okun(2001)}.
Equipped with a specific metric,
$\Sigma_L$ turns out to be non-positive curved in a combinatorial sense,
namely, it is a CAT(0)-space \cite[6.5.3]{Davis-Okun(2001)}.
If $L$ is a generalized rational homology $(n-1)$-sphere, i.e. a homology $(n-1)$-manifold
with the same rational homology as
$S^{n-1}$, then $\Sigma_L$ is a polyhedral
homology $n$-manifold with rational coefficients
\cite[7.4]{Davis-Okun(2001)}.
So $\Sigma_L$ is a reminiscence of the universal covering of a
closed $n$-dimensional  manifold with non-positive sectional curvature
and fundamental group $W_L$. In view of the Singer Conjecture
\ref{con: Singer Conjecture} the conjecture makes sense that
$b_p^{(2)}(\Sigma_L;\caln(W_L)) = 0$ for $2p \not= n$ provided that the underlying topological space
of $L$  is $S^{n-1}$ (or, more generally,  that it is a homology $(n-1)$-sphere)
\cite[Conjecture 0.4 and 8.1]{Davis-Okun(2001)}.
Davis and Okun show that the conjecture is true in dimension $n \le 4$
and that it is true in dimension $(n+1)$ if it holds
in dimension $n$ and $n$ is odd \cite[Theorem 9.3.1 and Theorem
10.4.1]{Davis-Okun(2001)}.
\em
\end{remark}


\subsection{The Singer Conjecture and K\"ahler Manifolds}

\label{subsec: The Singer Conjecture an Kaehler Manifolds}

\begin{definition} \label{def: d(bounded) form}
Let $(M,g)$ be a connected Riemannian manifold. A $(p-1)$-form
$\eta \in \Omega^{p-1}(M)$ is \emph{bounded}
if $||\eta||_{\infty} := \sup\{||\eta||_x \mid x \in M\} < \infty$ holds, where
$||\eta||_x$ is the norm on $\Alt^{p-1}(T_xM)$ induced by $g_x$.
A $p$-form $\omega \in \Omega^p(M)$ is called \emph{$d$(bounded)}
if $\omega = d(\eta)$ holds for some
bounded $(p-1)$-form $\eta \in \Omega^{p-1}(M)$.
A $p$-form $\omega \in \Omega^p(M)$ is called \emph{$\widetilde{d}$(bounded)}
if its lift
$\widetilde{\omega} \in \Omega^p(\widetilde{M})$ to the universal covering
$\widetilde{M}$ is $d$(bounded).
\end{definition}

The next definition is taken from \cite[0.3 on page 265]{Gromov(1991)}.
\begin{definition}[K\"ahler hyperbolic manifold]
\label{def: Kaehler hyperbolic manifold}
A \emph{K\"ahler hyperbolic manifold}%
\index{Kaehler hyperbolic manifold@K\"ahler hyperbolic manifold}
\index{manifold!Kaehler hyperbolic manifold@K\"ahler hyperbolic manifold}
is a closed connected K\"ahler manifold $(M,h)$ whose fundamental form $\omega$
is $\widetilde{d}$(bounded).
\end{definition}

\begin{example}[Examples of K\"ahler hyperbolic manifolds]
\label{exa: Kaehler hyperbolic manifolds} \em
The following list of examples of K\"ahler hyperbolic manifolds is
taken from \cite[Example 0.3]{Gromov(1991)}:

\begin{enumerate}

\item $M$ is a closed K\"ahler manifold which is homotopy equivalent to
a Riemannian manifold with negative sectional curvature;

\item $M$ is a closed K\"ahler manifold such that $\pi_1(M)$ is
word-hyperbolic in the sense of \cite{Gromov(1987)} and $\pi_2(M) = 0$;

\item $\widetilde{M}$ is a symmetric Hermitian space
of non-compact type;

\item $M$ is a complex submanifold of a K\"ahler hyperbolic manifold;

\item $M$ is a product of two K\"ahler hyperbolic manifolds.

\end{enumerate}
\em
\end{example}
The following result is due to Gromov
\cite[Theorem 1.2.B and Theorem 1.4.A on page 274]{Gromov(1987)}.
A detailed discussion of the proof and the consequences of this
theorem can also be found in \cite[Chapter 11]{Lueck(2002)}.

\begin{theorem}{\bf ($L^2$-Betti numbers and Novikov-Shubin invariants of
K\"ahler hyperbolic manifolds).}%
\index{Theorem!L2-Betti numbers and Novikov-Shubin invariants of
Kaehler hyperbolic manifolds@$L^2$-Betti numbers and Novikov-Shubin invariants
of K\"ahler hyperbolic manifolds}
 \label{the: b_p and alpha for Kaehler hyperbolic}
Let $M$ be a K\"ahler hyperbolic manifold of complex dimension $m$ and
real dimension $n=2m$. Then
\begin{eqnarray*}
b_p^{(2)}(\widetilde{M}) & = & 0 \hspace{10mm}  \mbox{ if } p \not= m;
\\
b_m^{(2)}(\widetilde{M}) & > & 0;
\\
(-1)^m \cdot \chi(M) & > & 0.
\end{eqnarray*}
\end{theorem}


\typeout{--------------------   Section 10 --------------------------}

\section{The Approximation Conjecture}
\label{sec: Approximation Conjecture}

This section  is devoted to the following conjecture.

\begin{conjecture}[Approximation Conjecture]
\label{con: Approximation Conjecture}%
\index{Conjecture!Approximation Conjecture}
A group $G$ satisfies the \emph{Approximation Conjecture}
if the following holds:

Let $\{G_i \mid i \in I\}$ be an inverse
system of normal subgroups of $G$ directed by inclusion over the directed set
$I$. Suppose that $\bigcap_{i \in I} G_i = \{1\}$.
Let $X$ be a $G$-$CW$-complex of finite type. Then $G_i\backslash X$
is a $G/G_i$-$CW$-complex of finite type and
$$b_p^{(2)}(X;\caln(G)) ~ = ~
\lim_{i \in I} b_p^{(2)}(G_i\backslash X;\caln(G/G_i)).$$
\end{conjecture}

\begin{remark}[The Approximation Conjecture for subgroups of finite index]
\label{rem: The Approximation Conjecture for subgroups of finite index}
\em 
Let us consider the special case where the inverse system $\{G_i \mid i \in I\}$ is
given by a nested sequence of normal subgroups of finite index
$$G = G_0 \supset G_1 \supset G_2 \supset G_3 \supset \ldots.$$
Notice that then  $b_p^{(2)}(G_i\backslash X;\caln(G/G_i)) ~ = ~
\frac{b_p(G_i\backslash X)}{[G : G_i]}$, where $b_p(G_i\backslash X)$
is the classical $p$-th Betti number of the finite $CW$-complex
$G_i\backslash X$. In this special case 
Conjecture \ref{con: Approximation Conjecture} was formulated by Gromov
\cite[pages 20, 231]{Gromov(1993)}
and proved in \cite[Theorem 0.1]{Lueck(1994c)}.
Thus we get an asymptotic relation between
the $L^2$-Betti numbers and Betti numbers, namely
$$b_p^{(2)}(X;\caln(G)) ~ = ~ 
\lim_{i \to \infty} \frac{b_p(G_i\backslash X)}{[G : G_i]},$$
although the Betti numbers of a connected finite $CW$-complex $Y$ and
the $L^2$-Betti numbers of its universal covering $\widetilde{Y}$ have nothing in common
except the fact that their alternating sum equals $\chi(Y)$ (see
Example \ref{exa: Independence of $L^2$-Betti numbers and Betti numbers}).

Interesting variations of this result for not necessarily normal subgroups of finite index and 
Betti-numbers with coefficients in representations can be found in the paper by Farber
\cite{Farber(1997)}.
\em
\end{remark}

\begin{definition} \label{def: class calg}
Let $\calg$%
\indexnotation{class calg}
be the smallest class of groups which contains the trivial group
and is closed under the following operations:
\begin{enumerate}

\item \label{def: class calg: amenable quotient}
Amenable quotient\\
Let $H \subseteq G$ be a (not necessarily normal) subgroup. Suppose that $H \in \calg$ and
the quotient $G/H$ is
an amenable discrete homogeneous space. (For the precise definition of
amenable discrete homogeneous space see for instance \cite[Definition 13.8]{Lueck(2002)}.
If $H\subseteq G$ is normal and $G/H$ is amenable, then $G/H$ is
an amenable discrete homogeneous space.)

Then $G \in \calg$;

\item \label{def: class calg: direct limit}
Colimits\\
If $G = \colim_{i \in I} G_i$ is the colimit of the
directed system $\{G_i \mid i \in I\}$ of groups indexed by the
directed set $I$ and each $G_i$
belongs to $\calg$, then $G$ belongs to $\calg$;

\item \label{def: class calg: inverse limit}
Inverse limits\\
If $G = \lim_{i \in I} G_i$ is the limit of the
inverse system $\{G_i \mid i \in I\}$ of groups indexed by the
directed set $I$ and each $G_i$
belongs to $\calg$, then $G$ belongs to $\calg$;

\item \label{def: class calg: subgroups}
Subgroups\\
If $H$ is isomorphic to a subgroup of the group $G$ with $G \in \calg$,
then $H \in \calg$;

\item \label{def: class calg: quotient with finite kernel}
Quotients with finite kernel\\
Let $1 \to K \to G \to Q \to 1$ be an exact sequence of groups. If $K$
is finite and $G$ belongs to $\calg$, then $Q$ belongs to $\calg$.
\end{enumerate}
\end{definition}

Next we provide some information about the class $\calg$.
Notice that in the original definition of $\calg$ due to 
Schick \cite[Definition 1.12]{Schick(2001b)} the resulting class is slightly smaller:
there it is required  that the class
contains the trivial subgroup and is closed under operations
\ref{def: class calg: amenable quotient}, \ref{def: class calg: direct limit},
\ref{def: class calg: inverse limit}  and \ref{def: class calg: subgroups}, but not
necessarily under operation \ref{def: class calg: quotient with finite kernel}.
The proof of the next lemma can be found in  \cite[Lemma 13.11]{Lueck(2002)}.

\begin{lemma} \label{lem: prop. of calg}
\begin{enumerate}

\item \label{lem: prop. of calg: fin. gen. sub.}
A group $G$ belongs to $\calg$ if and only if every finitely generated
subgroup of $G$ belongs to $\calg$;

\item \label{lem: prop. of calg: res. closed}
The class $\calg$ is residually closed, i.e. if there is a nested
sequence of normal subgroups $G = G_0 \supset G_1 \supset G_2 \supset \ldots$
such that $\bigcap_{i \ge 0} G_i = \{1\}$ and each quotient $G/G_i$ belongs to
$\calg$, then $G$ belongs to $\calg$;

\item \label{lem: prop. of calg: res. f. and am.}
Any residually amenable and in particular any residually finite
group belongs to $\calg$;

\item \label{lem: prop. of calg: mapping torus}
Suppose that $G$ belongs to $\calg$
and $f\colon G \to G$ is an endomorphism. Define the ``mapping torus
group'' $G_f$ to be the quotient of $G \ast \bbZ$ obtained by
introducing the relations $t^{-1}gt = f(g)$ for $g \in G$ and
$t \in \bbZ$ a fixed generator. Then $G_f$ belongs to $\calg$;

\item \label{lem: prop. of calg: direct sum and products}
Let $\{G_j \mid j \in J\}$ be a set of groups with $G_j \in
\calg$. Then the direct sum $\bigoplus_{j \in J} G_j$ and
the direct product $\prod_{j \in J} G_j$ belong
to $\calg$.

\end{enumerate}
\end{lemma}

The proof of the next result can  be found in \cite[Theorem 13.3]{Lueck(2002)}.
It is a mild generalization of the results of Schick~\cite{Schick(2000c)} and
\cite{Schick(2001b)}, where the original proof of the Approximation Conjecture for subgroups
of finite index was generalized to the much more general setting above 
and then applied to the Atiyah Conjecture. The connection between the Approximation
Conjecture and the Atiyah Conjecture for torsion-free groups comes from the obvious fact
that a convergent series of integers has an integer as limit.

\begin{theorem}[Status of the Approximation Conjecture]
\label{the: Status of the Approximation Conjecture}
\index{Theorem!Status of the Approximation Conjecture}
Every group $G$ which belongs to the class $\calg$
(see Definition~\ref{def: class calg})
 satisfies the Approximation
Conjecture~\ref{con: Approximation Conjecture}.
\end{theorem}


\typeout{--------------------   Section 11 --------------------------}


\section{$L^2$-Torsion}
\label{sec: L^2-Torsion}

Recall that $L^2$-Betti numbers are modelled on Betti numbers.
Analogously one can generalize the classical notion of 
Reidemeister torsion to an $L^2$-setting, which will lead to the notion
of $L^2$-torsion. The $L^2$-torsion may be viewed as a secondary $L^2$-Betti number
just as the Reidemeister torsion can be viewed as a secondary Betti number.
Namely, the Reidemeister torsion is only defined if all the Betti numbers (with coefficients
in a suitable representation) vanish, and similarly the $L^2$-torsion is defined only if
the $L^2$-Betti numbers vanish. Both invariants give valuable information about
the spaces in question.


\subsection{The Fuglede-Kadison Determinant}
\label{subsec: The Fuglede-Kadison Determinant}

In this subsection we briefly explain the notion of the Fuglede-Kadison 
determinant. We have extended the notion of the (classical) dimension of a 
finite  dimensional complex vector space to the von Neumann dimension of a 
finitely generated projective $\caln(G)$-module (and later even to arbitrary
$\caln(G)$-modules). Similarly we want to generalize the 
classical determinant  of an endomorphism
of a finite dimensional complex vector space to 
the Fuglede-Kadison determinant of an
$\caln(G)$-endomorphism $f \colon P \to P$ of
a finitely generated projective $\caln(G)$-module $P$ and of
an $\caln(G)$-map $f \colon \caln(G)^m \to \caln(G)^n$ of based finitely
generated $\caln(G)$-modules.
This is necessary since for the definition of Reidemeister torsion
one needs determinants and hence 
for the definition of $L^2$-torsion one has to develop an appropriate $L^2$-analogue.

\begin{definition}[Spectral density function]
\label{def: spectral density function}
Let $f \colon \caln(G)^m \to \caln(G)^n$ be an $\caln(G)$-homomorphism.
Let $\nu(f) \colon l^2(G)^m \to l^2(G)^n$  be the associated
bounded $G$-equivariant operator (see 
Remark~\ref{rem: elementary proof of semihereditary}).
Denote by $\{E_{\lambda}^{f^*f} \mid \lambda \in \bbR\}$ 
the (right-continuous)
family of spectral projections of the positive operator $\nu(f^*f)$.
Define the {\em spectral density function of}%
\index{spectral density function}
$f$ by
$$F_f%
\indexnotation{F_f}
\colon \bbR \to [0,\infty) \quad \lambda \mapsto
\dim_{\caln(G)}\left(\im(E_{\lambda^2}^{f^*f})\right).$$
\end{definition}

The spectral density function is monotone and right-continuous.
It takes values in $[0,m]$. Here and in the sequel
$|x|$ denotes the norm of an element
$x$ of a Hilbert space and 
$\|T\|$ the operator norm of a bounded operator $T$.
Since $\nu(f)$ and $\nu(f^{\ast}f)$ have the same kernel,
$\dim_{\caln(G)}(\ker(f))= F_f(0)$.

\begin{example}[Spectral density function for finite $G$]
 \label{spectral density function for finite G}\em
Suppose that $G$ is finite. Then 
$\bbC G = \caln(G) = l^2(G)$ and $\nu(f) = f$.
Let $0 \le \lambda_0 < \ldots < \lambda_r$ be the eigenvalues
of $f^{\ast}f$ and $\mu_i$ be the multiplicity of $\lambda_i$, i.e.
the dimension of the eigenspace of $\lambda_i$. Then the spectral
density function is a right continuous step function which is zero
for $\lambda < 0$ and has a step of height $\frac{\mu_i}{|G|}$
at each $\sqrt{\lambda_i}$. 
\em 
\end{example}

\begin{example}[Spectral density function for $G = \bbZ^n$]
\label{exa: spectral density function of a morphism for G=Z^n} \em
Let $G = \bbZ^n$. We use the  identification
$\caln(\bbZ^n) = L^{\infty}(T^n)$ of Example~\ref{exa: group von neumann algebra of Z^n}.
For $f \in L^{\infty}(T^n)$
the spectral density function $F_{M_f}$ of $M_f\colon L^2(T^n) \to
L^2(T^n), ~ g \mapsto g\cdot f$
sends $\lambda$ to the volume of the set
$\{z \in T^n \mid |f(z)| \le \lambda\}$. \em
\end{example}

\begin{definition}{\bf (Fuglede-Kadison determinant of $\caln(G)$-maps
$\caln(G)^m \to \caln(G)^n$)}.
 \label{def: Fuglede-Kadison determinant}
Let $f\colon \caln(G)^m \to \caln(G)^n$  be an $\caln(G)$-map.
Let $F_f(\lambda)$ be the spectral density
function of Definition \ref{def: spectral density function}
which is a monotone non-decreasing right-continuous function. Let $dF$
be the unique measure on the Borel $\sigma$-algebra on $\bbR$
which satisfies
$dF(]a,b]) = F(b)-F(a)$ for $a < b$. Then define the 
\emph{Fuglede-Kadison determinant}
\index{Fuglede-Kadison determinant!of $\caln(G)$-maps
$\caln(G)^m \to \caln(G)^n$}
$${\det}_{\caln(G )}(f) \in [0,\infty)$$
by the positive real number
$${\det}_{\caln(G )}(f)%
\indexnotation{det_{caln(G)}(f),based}
 ~ = ~
\exp\left(\int_{0+}^{\infty} \ln(\lambda) ~ dF\right)$$
if the Lebesgue integral
$\int_{0+}^{\infty} \ln(\lambda) ~ dF$ converges to a real number
and by $0$ otherwise. 
\end{definition}

Notice that in the definition above we do not require
$m = n$ or that $f$ is injective or $f$ is surjective.

\begin{example}[Fuglede-Kadison determinant for finite $G$] 
\label{exa: det for finite G} \em
To illustrate this definition, we look at
the example where $G$ is finite. We essentially get
the classical determinant $\det_{\bbC}$. Namely,
we have computed  the spectral density function for finite $G$
in Example \ref{spectral density function for finite G}.
Let $\lambda_1$, $\lambda_2$, $\ldots$, $\lambda_r$ be the non-zero
eigenvalues of $f^{\ast}f$ with multiplicity $\mu_i$. Then
one obtains, if $\overline{f^{\ast}f}$ is the automorphism
of the orthogonal complement of the kernel of $f^{\ast}f$ induced
by $f^{\ast}f$,
\begin{multline*}
{\det}_{\caln(G)}(f) 
~ = ~
\exp\left(
\sum_{i = 1}^r \frac{\mu_i}{|G|} \cdot \ln(\sqrt{\lambda_i})\right)
~ = ~
 \prod_{i=1}^r \lambda_i^{\frac{\mu_i}{2\cdot |G|}}
~ = 
\left({\det}_{\bbC}\left(\overline{f^{\ast}f}\right)\right)^{\frac{1}{2\cdot |G|}}.
\end{multline*}
If $f \colon \bbC G^m \to \bbC G^m$ is an automorphism, we get
$${\det}_{\caln(G)}(f) ~ = ~
\left|{\det}_{\bbC}(f)\right|^{\frac{1}{|G|}}. $$
\em
\end{example}

\begin{example}[Fuglede-Kadison determinant over $\caln(\bbZ^n)$]
\label{exa: Fuglede- Kadision determinant for G=Z^n} \em 
Let $G = \bbZ^n$. We use the  identification
$\caln(\bbZ^n) = L^{\infty}(T^n)$ of Example~\ref{exa: group von neumann algebra of Z^n}.
For $f \in L^{\infty}(T^n)$ we conclude from
Example~\ref{exa: spectral density function of a morphism for G=Z^n}
$${\det}_{\caln(\bbZ^n)}\left(M_f\colon  L^2(T^n) \to L^2(T^n)\right)
= 
\exp\left(\int_{T^n} \ln(|f(z)|) \cdot
\chi_{\{u \in S^1\mid f(u) \not= 0\}}~ dvol_z\right)$$
using the convention  $\exp(-\infty) = 0$.
\em
\end{example}

Here are some basic properties of this notion.
A morphism $f \colon \caln(G)^m \to \caln(G)^n$ has dense image if the closure
$\overline{\im(f)}$ of its image in $\caln(G)^n$ in the sense
of Definition~\ref{def: closure, K(M), P(M)} is $\caln(G)^n$.
The adjoint $A^*$ of a matrix $A = (a_{i,j})  \in M(m,n;\caln(G))$ is the matrix
in $M(n,m;\caln(G))$ given by $(a_{j,i}^*)$, where 
$\ast \colon \caln(G) \to \caln(G)$ sends an operator $a$ to its adjoint $a^*$.
The adjoint $f^* \colon \caln(G)^n \to \caln(G)^m$
of $f \colon \caln(G)^m \to \caln(G)^n$ is given by the matrix
$A^*$ if $f$ is given by the matrix $A$. The proof of the next result can be found
in \cite[Theorem 3.14]{Lueck(2002)}.

\begin{theorem}[Fuglede-Kadison determinant]
\index{Theorem!Fuglede-Kadison determinant}
\label{the: main properties of det}\ \\
\begin{enumerate}
\item Composition\\[1mm]
\label{the: main properties of det: composition}
Let $f\colon  \caln(G)^l \to \caln(G)^m$ and $g\colon  \caln(G)^m \to \caln(G)^n$ be
$\caln(G)$-homo\-mor\-phisms such that $f$ has dense image and
$g$ is injective. Then
$${\det}_{\caln(G)}(g \circ f) ~ = ~ {\det}_{\caln(G)}(g) \cdot {\det}_{\caln(G)}(f);$$

\item Additivity\\[1mm]
\label{the: main properties of det: additivity}
Let $f_1\colon  \caln(G)^{m_1} \to \caln(G)^{n_1}$, $f_2\colon  \caln(G)^{m_2} \to \caln(G)^{n_2}$ and
$f_3\colon  \caln(G)^{m_3} \to \caln(G)^{n_3}$ be $\caln(G)$-homomorphisms
such that  $f_1$ has dense image and $f_2$ is injective. Then
$${\det}_{\caln(G)}\squarematrix{f_1}{f_3}{0}{f_2}
~ = ~
{\det}_{\caln(G)}(f_1) \cdot {\det}_{\caln(G)}(f_2);$$

\item Invariance under adjoint map\\[1mm]

\label{the: main properties of det: det(f) = det(f^*)}
Let $f\colon  \caln(G)^m \to \caln(G)^n$ be an $\caln(G)$-homomorphism. Then
$${\det}_{\caln(G)}(f) = {\det}_{\caln(G)}(f^*);$$

\item Induction\\[1mm]
\label{the: main properties of det: induction}
Let $i\colon  H \to G$ be an injective group homomorphism and let
$f\colon  \caln(H)^m\to \caln(H)^n$ be an$\caln(H)$-homomorphism.
Then
$${\det}_{\caln(G)}(i_*f) ~ = ~ {\det}_{\caln(H)}(f).$$

\end{enumerate}
\end{theorem}

\begin{definition}{\bf (Fuglede-Kadison determinant of $\caln(G)$-endo\-mor\-phisms
of finitely generated projective modules).}
 \label{def: Fuglede-Kadison determinant of endos}
Let $f \colon P \to P$ be an endomorphism of a finitely generated projective
$\caln(G)$-module $P$.  Choose a finitely generated projective $\caln(G)$-module $Q$ 
and an $\caln(G)$-isomorphism $u \colon \caln(G)^n \xrightarrow{\cong} P \oplus Q$.
Define the \emph{Fuglede-Kadison determinant}%
\index{Fuglede-Kadison determinant!of endomorphisms
of finitely generated projective $\caln(G)$-modules}
$${\det}_{\caln(G)}(f)
 ~ \in ~ [0,\infty)$$
by the Fuglede-Kadison determinant in the sense
of Definition~\ref{def: Fuglede-Kadison determinant}
$${\det}_{\caln(G)}\left(u^{-1} \circ (f \oplus \id_Q) \circ u\right).$$
\end{definition}

This definition is independent of the choices of $Q$ and $u$ by
Theorem~\ref{the: main properties of det}. 
Notice that in Definition~\ref{def: Fuglede-Kadison determinant of endos}
no $\caln(G)$-basis appear but that it works only for endomorphisms, whereas in 
Definition~\ref{def: Fuglede-Kadison determinant} we work with finitely generated free
based modules but do not require that the source and target of $f$ are isomorphic.
There is an obvious analogue of Theorem~\ref{the: main properties of det}
for the Fuglede-Kadison determinant of endomorphisms of finitely generated projective
$\caln(G)$-modules.


\subsection{The Determinant Conjecture}
\label{subsec: The Determinant Conjecture}

It will be important for applications to geometry to study the Fuglede-Kadison
determinant of $\caln(G)$-maps $f \colon \caln(G)^m \to \caln(G)^n$ which
come by induction from $\bbZ G$-maps or $\bbC G$-maps. The following example is taken from
\cite[Example 3.22]{Lueck(2002)}.

\begin{example}[Fuglede-Kadison determinant of maps coming from elements in {$\bbC[\bbZ]$}] \em
\label{exa: Fuglede-Kadison determinant of maps coming from elements CZ} 
Consider a non-trivial element $p \in \bbC[\bbZ] = \bbC[z,z^{-1}]$.
We can write
$$p(z) ~ = ~ C \cdot z^n \cdot \prod_{k = 1}^l (z - a_k)$$
for non-zero complex numbers $C, a_1, \ldots , a_l$ and non-negative
integers $n,l$. Let $r_p \colon \caln(\bbZ) \to \caln(\bbZ)$
be the $\caln(\bbZ)$-map given by right multiplication with $p$. Then
\begin{eqnarray*}
{\det}_{\caln(\bbZ)}(r_p)
& = &
|C| \cdot
\prod_{\substack{1 \le k \le l, \\|a_k| > 1}} |a_k|.
\end{eqnarray*}
\em
\end{example}

\begin{definition}[Determinant class] \label{def: determinant class}
A group $G$ is of \emph{$\det \ge 1$-class}
\index{det ge 1 class@$\det \ge 1$-class}
\index{group!det ge 1 class@$\det \ge 1$-class} 
if for each $A \in M(m,n;\bbZ G)$ the Fuglede-Kadison determinant
(see Definition~\ref{def: Fuglede-Kadison determinant}) of the
morphism $r_A\colon \caln(G)^m \to \caln(G)^n$ given by right multiplication
with $A$ satisfies
$${\det}_{\caln(G)}(r_A) \ge 1.$$
\end{definition}

\begin{conjecture}[Determinant Conjecture]
\label{con: Determinant Conjecture}
\index{Conjecture!Determinant Conjecture}
Every group $G$ is of $\det \ge 1$-class.
\end{conjecture}

The proof of the next result can  be found in \cite[Theorem 13.3]{Lueck(2002)}.
It is a mild generalization of the results of Schick~\cite{Schick(2000c)} and
\cite{Schick(2001b)}.

\begin{theorem}[Status of the Determinant Conjecture]
\label{the: Status of the Determinant Conjecture}
\index{Theorem!Status of the Determinant Conjecture}
Every group $G$ which belongs to the class $\calg$ 
(see Definition~\ref{def: class calg}) satisfies the Determinant 
Conjecture~\ref{con: Determinant Conjecture}.
\end{theorem}

One easily checks that the Fuglede-Kadison determinant defines a homomorphism
of abelian groups 
\begin{eqnarray}
\Phi^G: \Wh(G) &\to & (0,\infty) = \{r \in \bbR \mid r > 0\}
\label{def: Phi^G}
\end{eqnarray}
with respect to the group structure given by multiplication
of positive real numbers on the target. We mention the following conjecture.

\begin{conjecture}[Triviality of the map induced by the Fuglede-Kadison determinant on $\Wh(G)$]
\index{Conjecture!Triviality of the map induced by the Fuglede-Kadison determinant on $\Wh(G)$} 
\label{con: Triviality of the map induced by the Fuglede-Kadison determinant on Wh(G)}
The map $\Phi^G\colon  \Wh(G) ~ \to ~ (0,\infty)$ is trivial.
\end{conjecture}

\begin{lemma} \label{Det implies approxi and triviality of phi^G}
\begin{enumerate} 

\item \label{Det implies approxi and triviality of phi^G: phi}
If $G$ satisfies the Determinant Conjecture~\ref{con: Determinant Conjecture},
then $G$ satisfies 
Conjecture~\ref{con: Triviality of the map induced by the Fuglede-Kadison determinant on Wh(G)};

\item \label{Det implies approxi and triviality of phi^G: approxi}
The Approximation Conjecture \ref{con: Approximation Conjecture} for $G$ and the inverse system
$\{G_i \mid i \in I\}$ is true if each group
$G_i$ is of $\det \ge 1$-class.

\end{enumerate}
\end{lemma}
\proof See \cite[Theorem 13.3 (1) and Lemma 13.6]{Lueck(2002)}. \qed


\subsection{Definition and Basic Properties of $L^2$-Torsion}
\label{subsec: Definition and Basic Properties of L^2-Torsion}

We will consider  $L^2$-torsion only for universal coverings and in the $L^2$-acyclic  case.
A more general setting is  treated in \cite[Section 3.4]{Lueck(2002)}.

\begin{definition}[$\det$-$L^2$-acyclic]  \label{def: det-l^2-acyclic}
Let $X$ be a finite connected $CW$-complex with fundamental group $\pi = \pi_1(X)$.
 Let $C_*^{\caln}(\widetilde{X})$ be the 
$\caln(\pi)$-chain complex $\caln(G) \otimes_{\bbZ G} C_*(\widetilde{X})$ 
with $p$-th differential $c^{\caln}_p = \id_{\caln(G)} \otimes_{\bbZ } c_p$. We say that
$X$ is \emph{$\det$-$L^2$-acyclic}%
\index{det_L2-acyclic@$\det$-$L^2$-acyclic}
if for each $p$ we get for the Fuglede-Kadison determinant of $c_p^{\caln}$
and for the $p$-th $L^2$-Betti number of $\widetilde{X}$
\begin{eqnarray*}
{\det}_{\caln(\pi)}\left(c_p^{\caln}\right)  & > & 0;
\\
b_p^{(2)}(\widetilde{X}) & = & 0.
\end{eqnarray*}

If $X$ is $\det$-$L^2$-acyclic, we define the \emph{$L^2$-torsion}%
\index{L2-torsion@$L^2$-torsion}
of $\widetilde{X}$ by
$$\rho^{(2)}(\widetilde{X})%
\indexnotation{rho^(2)(widetilde X)}
~ = ~  
- \sum_{p \ge 0} (-1)^p \cdot \ln\left({\det}_{\caln(G)}\left(c_p^{\caln}\right)\right)
\quad \in \bbR.$$

If $X$ is a finite $CW$-complex, we call it $\det$-$L^2$-acyclic if each component $C$ is
$\det$-$L^2$-acyclic. In this case we define 
$$\rho^{(2)}(\widetilde{X}) ~ := ~ \sum_{C \in \pi_0(X)} ~ \rho^{(2)}(\widetilde{C}).$$ 
\end{definition}

The condition that $\widetilde{X}$ is $L^2$-acyclic is not needed  for the definition of
$L^2$-torsion, but is necessary to ensure the basic and useful properties 
which we will discuss below.

\begin{remark}[$L^2$-torsion in terms of the Laplacian]
\label{rem: L^2-torsionin terms of the Laplacian}
\em
One can express the $L^2$-torsion also in terms of the Laplacian
which is closer to the  notions of analytic torsion and analytic $L^2$-torsion.
After a choice of cellular $\bbZ \pi$-basis, every $\caln(\pi)$-chain module
$C_p^{\caln}(\widetilde{X})$ looks like $\caln(\pi)^{n_p}$ for appropriate non-negative integers
$n_p$. Hence we can assign to 
$c_p^{\caln} \colon C_p^{\caln}(\widetilde{X}) \to C_{p-1}^{\caln}(\widetilde{X})$ its adjoint
$\left(c_p^{\caln}\right)^* \colon C_{p-1}^{\caln}(\widetilde{X}) \to  C_p^{\caln}(\widetilde{X})$
which is given  by the matrix $(a_{j,i}^*)$ if $c_p^{(2)}$ is given by the matrix $(a_{i,j})$.
Define the $p$-th Laplace homomorphism 
$\Delta_p \colon C_p^{\caln}(\widetilde{X}) \to C_p^{\caln}(\widetilde{X})$ 
to be the $\caln(\pi)$-homomorphism
$\left(c_p^{\caln}\right)^* \circ c_p^{\caln} + c_{p+1}^{\caln} \circ \left(c_{p+1}^{\caln}\right)^*$.
Then $\widetilde{X}$ is $\det$-$L^2$-acyclic if and only if
$\Delta_p$ is injective and has dense image, i.e. the closure of its image in 
$C_p^{\caln}(\widetilde{X})$ is $C_p^{\caln}(\widetilde{X})$, and $\det_{\caln(\pi)}(\Delta_p) > 0$.
In this case we get
$$\rho^{(2)}(\widetilde{X}) ~ = ~ - \frac{1}{2} \cdot \sum_{p \ge 0}(-1)^p \cdot p \cdot 
\ln\left({\det}_{\caln(\pi)}(\Delta_p)\right).$$
This follows from \cite[Lemma 3.30]{Lueck(2002)}.
\em
\end{remark}

The next theorem presents the basic properties of
$\rho^{(2)}(\widetilde{X})$ and is proved in 
\cite[Theorem 3.96]{Lueck(2002)}.
Notice the formal analogy between the behaviour
of $\rho^{(2)}(\widetilde{X})$ and the classical
Euler characteristic $\chi(X)$.

\begin{theorem}{\bf (Cellular $L^2$-torsion for
  universal coverings).}
\index{Theorem!Cellular $L^2$-torsion for universal coverings}
\label{the: main properties of rho2(widetildeX)}

\begin{enumerate}

\item \label{the: main properties of rho2(widetildeX): homotopy invariance}
Homotopy invariance
\\[1mm]
Let $f\colon  X \to Y$ be a homotopy equivalence of finite $CW$-complexes.
Let $\tau(f) \in \Wh(\pi_1(Y))$ be its Whitehead torsion (see
\cite{Cohen(1973)}). Suppose that
$\widetilde{X}$ or $\widetilde{Y}$ is $\det$-$L^2$-acyclic.
Then both $\widetilde{X}$ and  $\widetilde{Y}$ are $\det$-$L^2$-acyclic and
$$\rho^{(2)}(\widetilde{Y}) - \rho^{(2)}(\widetilde{X}) ~ = ~
\Phi^{\pi_1(Y)}(\tau(f)),$$
where
$\Phi^{\pi_1(Y)}\colon  \Wh(\pi_1(Y)) = \bigoplus_{C \in \pi_0(Y)} \Wh(\pi_1(C))
\to \bbR$
is the sum of the maps $\Phi^{\pi_1(C)}$ of
\eqref{def: Phi^G};

\item \label{the: main properties of rho2(widetildeX): sum formula}
Sum formula\\[1mm]
Consider the pushout of finite
$CW$-complexes such that $j_1$ is an inclusion of $CW$-complexes,
$j_2$ is cellular and $X$ inherits its $CW$-complex structure from $X_0$, $X_1$ and $X_2$
$$\comsquare{X_0}{j_1}{X_1}{j_2}{i_1}{X_2}{i_2}{X}$$
Assume
$\widetilde{X_0}$, $\widetilde{X_1}$ and $\widetilde{X_2}$ are
$\det$-$L^2$-acyclic and that
for $k=0,1,2$ the map
$\pi_1(i_k)\colon  \pi_1(X_k) \to \pi_1(X)$ induced by the obvious
map $i_k\colon  X_k \to X$ is injective for all base points in $X_k$.

Then $\widetilde{X}$ is $\det$-$L^2$-acyclic and we get
$$\rho^{(2)}(\widetilde{X}) ~ = ~
\rho^{(2)}(\widetilde{X_1}) + \rho^{(2)}(\widetilde{X_2}) -
\rho^{(2)}(\widetilde{X_0});$$

\item \label{the: main properties of rho2(widetildeX): Poincar'e duality}
Poincar\'e duality\\[1mm]
Let $M$ be a closed manifold of even dimension.
Equip it with some $CW$-complex structure. Suppose that
$\widetilde{M}$ is $\det$-$L^2$-acyclic. Then 
$$\rho^{(2)}(\widetilde{M}) ~ = ~ 0;$$

\item \label{the: main properties of rho2(widetildeX): product formula}
Product formula\\[1mm]
Let $X$ and $Y$ be finite $CW$-complexes. Suppose that
$\widetilde{X}$ is $\det$-$L^2$-acyclic. Then
$\widetilde{X \times Y}$ is $\det$-$L^2$-acyclic and
$$\rho^{(2)}(\widetilde{X \times Y}) ~ = ~
\chi(Y) \cdot \rho^{(2)}(\widetilde{X});$$

\item \label{the: main properties of rho2(widetildeX): multiplicativity}
Multiplicativity\\[1mm]
Let $X \to Y$ be a finite covering of finite $CW$-complexes
with $d$ sheets. Then $\widetilde{X}$ is $\det$-$L^2$-acyclic
if and only if $\widetilde{Y}$ is $\det$-$L^2$-acyclic and in this case
$$\rho^{(2)}(\widetilde{X}) ~ = ~ d \cdot \rho^{(2)}(\widetilde{Y}).$$

\end{enumerate}
\end{theorem}

The next three results are taken from \cite[Corollary 3.103, Theorem 3.105
and Theorem 3.111]{Lueck(2002)}. There is also a more general version of
Theorem~\ref{the: rho^2 and fiber bundles} for fibrations (see
\cite[Theorem 3.100]{Lueck(2002)}).

\begin{theorem}[$L^2$-torsion and fiber bundles]
 \label{the: rho^2 and fiber bundles}
Suppose that $F \to E \xrightarrow{p} B$ is a (locally trivial)
fiber bundle of finite $CW$-complexes with $B$ connected.
Suppose that for one (and hence all) $b \in B$ the inclusion
of the fiber $F_b$ into $E$ induces an injection on the fundamental groups
for all base points in $F_b$ and $\widetilde{F_b}$ is $\det$-$L^2$-acyclic.
Then $\widetilde{E}$ is $\det$-$L^2$-acyclic and
$$\rho^{(2)}(\widetilde{E})
~ = ~
\chi(B) \cdot \rho^{(2)}(\widetilde{F}).$$
\end{theorem}

\begin{theorem}[$L^2$-torsion and $S^1$-actions]
\index{Theorem!L2-torsion and S^1-actions@$L^2$-torsion and $S^1$-actions}
\label{the: $S^1$-actions and $L^2$-torsion}
Let $X$ be a connected $S^1$-$CW$-complex of finite type.
Suppose that for one orbit $S^1/H$ (and hence for all orbits) the inclusion
into $X$ induces a map on $\pi_1$ with infinite image.
(In particular the $S^1$-action has no fixed points.)
 Then $\widetilde{X}$ is $\det$-$L^2$-acyclic and 
\begin{eqnarray*}
\rho^{(2)}(\widetilde{X}) & = & 0.
\end{eqnarray*}
\end{theorem}

\begin{theorem}[$L^2$-torsion on aspherical closed $S^1$-manifolds]
\index{Theorem!$L^2$-torsion on aspherical closed $S^1$-manifolds}
\label{the: $L^2$-torsion on aspherical closed $S^1$-manifolds}
Let $M$ be an aspherical closed manifold with
non-trivial $S^1$-action. Then the action has no fixed
points and the inclusion of any orbit into $M$ induces an
injection on the fundamental groups. Moreover,
$\widetilde{M}$ is $\det$-$L^2$-acyclic  and
\begin{eqnarray*}
\rho^{(2)}(\widetilde{M}) & = & 0.
\end{eqnarray*}
\end{theorem}

The assertion for the $L^2$-torsion in the theorem below is the main result of
\cite{Wegner(2000)} (see also \cite{Wegner(2001)}).
Its proof  is based on localization techniques.

\begin{theorem}[$L^2$-torsion and aspherical $CW$-complexes]
\index{Theorem!L2-torsion and aspherical CW-complexes@$L^2$-torsion and aspherical $CW$-complexes}
\label{the: l2-tors. for asp. CW-com. with e.m-am sub pi}
Let $X$ be an aspherical  finite  $CW$-complex. Suppose that
its fundamental group $\pi_1(X)$ contains an elementary amenable infinite normal
subgroup $H$ and $\pi_1(X)$ is of $\det \ge 1$-class.
Then $\widetilde{X}$ is $\det$-$L^2$-acyclic and
\begin{eqnarray*}
\rho^{(2)}(\widetilde{X}) & = & 0.
\end{eqnarray*}
\end{theorem}

\begin{remark}[Homotopy invariance of $L^2$-torsion] 
\label{rem: Homotopy invariance of L^2-torsion} \em
Notice that Conjecture~
\ref{con: Triviality of the map induced by the Fuglede-Kadison determinant on Wh(G)}
implies because of Theorem~\ref{the: main properties of rho2(widetildeX)}
\ref{the: main properties of rho2(widetildeX): homotopy invariance} the homotopy 
invariance of the $L^2$-torsion. i.e. for two homotopy equivalent $\det$-$L^2$-acyclic
finite $CW$-complexes $X$ and $Y$ we have 
$\rho^{(2)}(\widetilde{X}) = \rho^{(2)}(\widetilde{Y})$.  
\em
\end{remark}


\subsection{Computations of $L^2$-Torsion}
\label{subsec: Computations of L^2-Torsion}

\begin{remark}[Analytic $L^2$-torsion] \em
\label{rem: anayltic $L^2$-torsion} 
It is important to know for the following specific calculations
that there is an analytic version of $L^2$-torsion in terms of the
heat kernel due to  Lott \cite {Lott(1992a)} and Mathai
\cite{Mathai(1992)} and that a deep result of Burghelea, Friedlander, Kappeler and McDonald
\cite{Burghelea-Friedlander-Kappeler-McDonald(1996a)} says that the
analytic one agrees with the one presented here.
\em
\end{remark}

The following result is due to  Hess and Schick \cite{Hess-Schick(1998)}.

\begin{theorem}[Analytic $L^2$-torsion of hyperbolic manifolds]
\index{Theorem!Analytic $L^2$-torsion of hyperbolic manifolds}
 \label{the: analytic L2-torsion of hyperbolic manifolds}
\ \\
Let $d = 2n+1$ be an odd integer. To $d$ one can associate an explicit  real number
$C_d >0$ with the following property:

For every closed hyperbolic $d$-dimensional manifold $M$ we have
$$\rho^{(2)}(\widetilde{M}) ~ = ~ (-1)^n \cdot C_d \cdot \vol(M),$$
where $\vol(M)$%
\indexnotation{vol(M)}
 is the volume of $M$.
\end{theorem}

The existence of a real number $C_d$ with 
$\rho^{(2)}(\widetilde{M}) ~ = ~ (-1)^n \cdot C_d \cdot \vol(M)$ follows from
the version of the Proportionality Principle for $L^2$-Betti numbers
(see Theorem~\ref{the: proportionality principle for $L^2$-Betti numbers})
for $L^2$-torsion (see \cite[Theorem 3.183]{Lueck(2002)}).
The point is that this number $C_d$ is given explicitly. 
For instance $C_3 = \frac{1}{6\pi}$ and $C_5 = \frac{31}{45\pi^2}$.
For each odd $d$ there exists a rational number $r_d$ such that
$C_d = \pi^{-n} \cdot r_d$
holds. The proof of this result is based on calculations involving 
the heat kernel on hyperbolic space.

\begin{remark}[$L^2$-torsion of symmetric spaces of non-compact type]
\label{rem: L2-torsion od symmetric spaces of non-compact type} \em
More generally, the $L^2$-torsion $\rho^{(2)}(\widetilde{M})$
for an aspherical closed manifold $M$ whose universal covering $\widetilde{M}$ is a symmetric space
is computed by Olbricht \cite{Olbrich(2000)}. \em
\end{remark}

The following result is proved in \cite[Theorem 0.6]{Lueck-Schick(1999)}.

\begin{theorem}[$L^2$-torsion of $3$-manifolds]
\index{Theorem!L2-torsion of 3-manifolds@$L^2$-torsion of $3$-manifolds}
\label{the: $L^2$-torsion of irreducible $3$-manifold}
Let $M$ be a compact connected orientable prime $3$-manifold
with infinite fundamental group such that the boundary
of $M$ is empty or a disjoint union of incompressible tori. Suppose that
$M$ satisfies Thurston's Geometrization Conjecture, i.e.
there is a geometric toral splitting 
along disjoint incompressible $2$-sided tori
in $M$ whose pieces are Seifert manifolds or hyperbolic manifolds.
Let $M_1$, $M_2$, $\ldots$,  $M_r$ be the hyperbolic pieces.
They all have finite volume \cite[Theorem B on page 52]{Morgan(1984)}.
Then $\widetilde{M}$ is $\det$-$L^2$-acyclic and
$$\rho^{(2)}(\widetilde{M}) ~ = ~
-\frac{1}{6\pi}\cdot \sum_{i=1}^r \vol(M_i).$$
In particular,
$\rho^{(2)}(\widetilde{M})$ is $0$ if and and only if
there are no hyperbolic pieces. 
\end{theorem}


\subsection{Some Open Conjectures about $L^2$-Torsion}
\label{subsec: Some Open Conjectures about $L^2$-Torsion}

All the computations and results above 
give evidence and are compatible with the following conjectures about
$L^2$-torsion taken from \cite[Theorem 11.3]{Lueck(2002)}.

\begin{conjecture}[$L^2$-torsion for aspherical manifolds]%
\index{Conjecture!L2-torsion for aspherical manifolds@
$L^2$-torsion for aspherical manifolds}
\label{con: $L^2$-torsion for aspherical manifolds}
If $M$ is an aspherical closed manifold of odd dimension, then
$\widetilde{M}$ is $\det$-$L^2$-acyclic  and
$$(-1)^{\frac{\dim(M)-1}{2}} \cdot \rho^{(2)}(\widetilde{M}) ~ \ge ~
0.$$
If $M$ is a closed connected Riemannian manifold of odd dimension with negative
sectional curvature, then
$\widetilde{M}$ is $\det$-$L^2$-acyclic  and
$$(-1)^{\frac{\dim(M)-1}{2}} \cdot \rho^{(2)}(\widetilde{M}) ~ > ~
0.$$
If $M$ is an aspherical closed manifold
whose fundamental group contains an amenable infinite normal subgroup,
then $\widetilde{M}$ is $\det$-$L^2$-acyclic and
$$\rho^{(2)}(\widetilde{M}) ~ = ~ 0.$$
\end{conjecture}

Consider a closed orientable manifold $M$ of dimension $n$.
Let $[M;\bbR]$ be the image of the fundamental class $[M] \in H_n^{\sing}(M;\bbZ)$ 
under the change of coefficient map
$H_n^{\sing}(M;\bbZ) \to H_n^{\sing}(M;\bbR)$. Define the $L^1$-norm on $C_n^{\sing}(M;\bbR)$ by
sending $\sum_{i=1}^s r_i \cdot [\sigma_i \colon \Delta_n \to M]$ to
$\sum_{i = 1}^s |r_i|$. It induces a seminorm on $H_n(M;\bbR)$.
Define the \emph{simplicial volume}%
\index{simplicial volume}
$||M||%
\indexnotation{||M||}
 \in \bbR$
to be the seminorm of $[M;\bbR]$. More information about the simplicial volume
can be found for instance in
\cite{Gromov(1982)}, \cite{Gromov(1999a)} and  \cite{Ivanov(1987)}, and 
in \cite[Chapter 14]{Lueck(2002)}, where also the following conjecture is discussed.

\begin{conjecture}[Simplicial volume and $L^2$-invariants]
\label{con: simplicial volume and L^2-invariants}
\index{Conjecture!Simplicial volume and $L^2$-invariants}
Let $M$ be an aspherical closed orientable manifold of dimension $\ge 1$.
Suppose that its simplicial volume $||M||$ vanishes. Then
$\widetilde{M}$ is $\det$-$L^2$-acyclic and 
\begin{eqnarray*}
\rho^{(2)}(\widetilde{M}) & = & 0.
\end{eqnarray*}
\end{conjecture}

The simplicial volume is a special invariant concerning bounded cohomology. 
The point of this conjecture is that it suggests a connection between bounded
cohomology and $L^2$-invariants such as $L^2$-cohomology and $L^2$-torsion.

We have already seen that $L^2$-Betti numbers are up to scaling invariant
under measure equivalence. The next conjecture is interesting
because it would give a sharper invariant in case all the $L^2$-Betti numbers vanish,
namely the vanishing of the $L^2$-torsion.

\begin{conjecture}[Measure equivalence and $L^2$-torsion]
\label{con: Measure equivalence and L^2-torsion}
\index{Conjecture!Measure equivalence and $L^2$-torsion}
Let $G_i$ for $i = 0,1$ be a group such that there is a finite $CW$-model for $BG_i$ and
$EG_i$ is $\det$-$L^2$-acyclic. Suppose that $G_0$ and $G_1$ are
measure equivalent. Then 
$$\rho^{(2)}(EG_0;\caln(G_0)) = 0 ~ \Leftrightarrow ~ \rho^{(2)}(EG_1;\caln(G_1)) = 0.$$
\end{conjecture}


\subsection{$L^2$-Torsion of Group Automorphisms}
\label{subsec: $L^2$-Torsion of Group Automorphisms}

In this section we explain that for a group automorphism $f\colon  G \to G$
the $L^2$-torsion applied to the $(G\rtimes_f \bbZ)$-$CW$-complex
$E(G\rtimes_f \bbZ)$  gives an  interesting new invariant, provided that
$G$ is of $\det \ge 1$-class and satisfies certain finiteness assumptions.
It seems to be worthwhile to investigate it further.
The following definition and theorem are taken from 
\cite[Definition 7.26 and Theorem 7.27]{Lueck(2002)}.

\begin{definition}[$L^2$-torsion of group automorphisms] \label{def: rho(2)(G_f)}
Let $f\colon  G \to G$ be a group automorphism. Suppose that there is a
finite $CW$-model for $BG$ and $G$ is of $\det \ge 1$-class. Define the
\emph{$L^2$-torsion}%
\index{L2-torsion@$L^2$-torsion!of a group automorphism}
of $f$ by
$$\rho^{(2)}(f\colon  G \to G)%
\indexnotation{rho^{(2)}(f: G to G)}
  ~ := \rho^{(2)}(\widetilde{B(G\rtimes_f\bbZ)}) \hspace{5mm}  \in ~ \bbR.$$
\end{definition}

Next we present the basic properties of this invariant. Notice that its behaviour
is similar to the Euler characteristic $\chi(G) := \chi(BG)$.

\begin{theorem} \label{the: basic properties of L^2-torsion}
Suppose that all groups appearing below have finite $CW$-models for their classifying
spaces and are of $\det \ge 1$-class.
\begin{enumerate}

\item \label{the: basic properties of L^2-torsion: amalgamated products}
Suppose that $G$ is the amalgamated product $G_1\ast_{G_0} G_2$
for subgroups $G_i \subseteq G$ and the automorphism
$f\colon  G \to G$ is the amalgamated product $f_1 \ast_{f_0} f_2$
for automorphisms $f_i\colon  G_i \to G_i$. Then
$$\rho^{(2)}(f) ~ = ~ \rho^{(2)}(f_1) + \rho^{(2)}(f_2) - \rho^{(2)}(f_0);$$

\item \label{the: basic properties of L^2-torsion: trace property}
Let $f\colon  G \to H$ and $g\colon  H \to G$ be isomorphisms of groups.  Then
$$\rho^{(2)}(f \circ g) ~ =  ~ \rho^{(2)}(g \circ f).$$
In particular $\rho^{(2)}(f)$ is invariant under conjugation
with automorphisms;

\item \label{the: basic properties of L^2-torsion: additivity}
Suppose that the following diagram of groups
$$
\begin{CD}
1 @>>> G_1 @>>> G_2 @>>> G_3 @>>> 1
\\
 & & @V f_1 VV  @V f_2 VV  @V \id VV
\\
1 @>>> G_1 @>>> G_2 @>>> G_3 @>>> 1
\end{CD}
$$
commutes, has exact rows and its vertical arrows are automorphisms. Then
$$\rho^{(2)}(f_2) ~ = ~ \chi(BG_3) \cdot \rho^{(2)}(f_1);$$

\item \label{the: basic properties of L^2-torsion: multiplicativity}
Let $f\colon  G \to G$ be an automorphism of a group. Then
for all integers $n \ge 1$
$$\rho^{(2)}(f^n) ~ = ~ n \cdot \rho^{(2)}(f);$$

\item\label{the: basic properties of L^2-torsion: subgroups of finite index}
Suppose that $G$ contains a subgroup $G_0$ of finite index $[G:G_0]$.
Let $f\colon  G \to G$ be an automorphism with $f(G_0) = G_0$.
Then
$$\rho^{(2)}(f) ~ = ~ \frac{1}{[G:G_0]}\cdot \rho^{(2)}(f|_{G_0});$$

\item \label{the: basic properties of L^2-torsion: G acyclic}
We have $\rho^{(2)}(f) = 0$ if $G$ satisfies one of the following conditions:

\begin{enumerate}

\item All the $L^2$-Betti numbers of $G$ vanish;

\item $G$ contains an amenable infinite normal subgroup.

\end{enumerate}

\end{enumerate}
\end{theorem}

\begin{example}[Automorphisms of surfaces] \label{exa: Automorphisms of surfaces}
\em
Using Theorem~\ref{the: $L^2$-torsion of irreducible $3$-manifold} one can compute the
$L^2$-torsion of the automorphism $\pi_1(f)$ for an automorphism
$f \colon S \to S$ of a compact connected orientable surface, possibly with boundary.
Suppose that $f$ is irreducible. Then the following statements are
equivalent:
i.)  $f$ is pseudo-Anosov, ii.) The mapping torus $T_f$ has a hyperbolic structure and iii.)
$\rho^{(2)}(\pi_1(f)) <  0$. Moreover, $f$ is periodic if and only if $\rho^{(2)}(\pi_1(f)) =  0$.
(see \cite[Subsection 7.4.2]{Lueck(2002)}.

The $L^2$-torsion of a Dehn twist is always zero since the associated mapping torus 
contains no hyperbolic pieces in his Jaco-Shalen-Johannson-Thurston splitting.
\em
\end{example}

\begin{remark}[Weaker finiteness conditions] \label{rem: Weaker finiteness conditions} \em
The definition of the $L^2$-torsion of a group automorphism above still makes sends
and has still most of the  properties above, if one weakens the condition that there is a finite model for $BG$
to the assumption that there is a finite model for the classifying space of proper $G$-actions
$\underline{E}G = \EGF{G}{\calfin}$. This is explained in \cite[Subsection 7.4.4]{Lueck(2002)}.
\em
\end{remark}


\typeout{--------------------   Section 12 --------------------------}

\section{Novikov-Shubin Invariants}
\label{sec: Novikov-Shubin Invariants}

In this section we briefly discuss Novikov-Shubin invariants. They
were originally defined in terms of heat kernels. We will focus on
their algebraic definition and aspects.


\subsection{Definition of Novikov-Shubin Invariants}
\label{subsec: Definition of Novikov-Shubin Invariants}

Let $M$ be a finitely presented $\caln(G)$-module. 
Choose some exact sequence
$\caln(G)^m \xrightarrow{f} \caln(G)^n \to M \to 0$. Let
$F_f$ be the  spectral density function of $f$
(see Definition~\ref{def: spectral density function}). Recall that $F_f$ is a monotone increasing
right continuous function $[0,\infty) \to [0,\infty)$. Define the \emph{Novikov-Shubin invariant}%
\index{Novikov-Shubin invariant!of a finitely generated $\caln(G)$-module}
of $M$ by 
$$\alpha(M)%
\indexnotation{alpha(M)}
~  = ~ \liminf_{\lambda \to 0^+}
\frac{\ln(F_f(\lambda) - F_f(0))}{\ln(\lambda)} \; \in [0, \infty ],$$
provided that $F_f(\lambda) > F_f(0)$ holds for all $\lambda > 0$. Otherwise,
one puts formally 
$$\alpha(M) ~ = ~ \infty^+.$$
It measures how fast $F_f(\lambda)$ approaches $F_f(0)$ for $\lambda\to 0^+$. 
For instance, if $F_f(\lambda) = \lambda^{\alpha}$ for $\lambda > 0$, then $\alpha(M) = \alpha$.
The proof that $\alpha(M)$ is independent of the choice of $f$ is analogous to the proof
of \cite[Theorem 2.55 (1)]{Lueck(2002)}. 

\begin{definition}[Novikov-Shubin invariants] 
\label{def: Novikov-Shubin invariant}
Let $X$ be a $G$-$CW$-complex of finite type. Define its \emph{$p$-th Novikov-Shubin invariant}%
\index{Novikov-Shubin invariant!$p$-th Novikov-Shubin invariant of a 
$G$-$CW$-complex of finite type} 
by
$$\alpha_p(X;\caln(G)) ~ = ~ \alpha\left(H_{p-1}^{(2)}(X;\caln(G))\right) 
\quad \in [0,\infty] \amalg \{\infty^+\}.$$
If the group $G$ is clear from the context, we abbreviate 
$\alpha_p(X) ~ = ~ \alpha_p(X;\caln(G))$.
\end{definition}

Notice that $H_{p-1}^{(2)}(X;\caln(G))$ is finitely presented since $\caln(G)$ is
semihereditary (see Theorem~\ref{the: von Neumann algebras are semihereditary}) 
and $C_k^{(2)}(X)$ is a finitely  generated free $\caln(G)$-module
for all $k \in \bbZ$ because $X$ is by assumption of finite type.

\begin{remark}[Analytic definition of Novikov-Shubin invariants]
 \label{rem: Analytic definition of Novikov-Shubin invariants} \em
Novi\-kov-Shubin invariants were originally analytically defined by Novikov and Shubin 
(see \cite{Novikov-Shubin(1986b)}, \cite{Novikov-Shubin(1986a)}).
For a cocompact smooth $G$-manifold $M$ without boundary and with $G$-invariant
Riemannian metric one can assign to its $p$-th Laplace operator $\Delta_p$ a density
function $F_{\Delta_p}(\lambda) = \tr_{\caln(G)}(E_{\lambda})$ for $\{E_{\lambda} \mid
\lambda \in [0,\infty)\}$ the spectral family associated to the essentially selfadjoint
operator $\Delta_p$. Define 
$\alpha_p^{\Delta}(M;\caln(G))%
\indexnotation{alpha_p^Delta(M)}
 \in [0,\infty] \amalg\{\infty^+\}$ 
by the same expression as appearing in the definition of $\alpha(M)$ above, only replace
 $F_f$ by $F_{\Delta_p}$. Then $\alpha_p^{\Delta}(M)$ agrees with $\frac{1}{2} \cdot \min\{\alpha_p(K),
\alpha_{p+1}(K)\}$ for any equivariant triangulation $K$ of $M$. For a proof of this
equality see \cite{Efremov(1991a)} or \cite[Section 2.4]{Lueck(2002)}. One can define
the analytic Novikov-Shubin invariant $\alpha_p^{\Delta}(M;\caln(G))$ also in terms
of heat kernels. It measures how fast the function 
$\int_{\calf} \tr_{\bbC}(e^{-t\Delta_p}(x,x)) ~ dvol_x$ 
approaches for $t \to \infty$ its limit
\begin{eqnarray*}
b_p^{(2)}(M)
& = &
\lim_{t \to \infty} \int_{\calf} \tr_{\bbC}(e^{-t\Delta_p}(x,x)) ~ dvol_x .
\end{eqnarray*}
The ``thinner'' the spectrum of $\Delta_p$ is at zero, the larger is 
$\alpha_p^{\Delta}(M;\caln(G))$. 

In view of this original analytic definition the result due to Gromov and Shubin
\cite{Gromov-Shubin(1991)} that
the Novikov-Shubin invariants are homotopy invariants, is rather surprising. 
\em
\end{remark}

\begin{remark}[Analogy to finitely generated $\bbZ$-modules] 
\label{rem: Analogy to finitely generated Z-modules} \em
Recall Slogan~\ref{slo: caln(G) like Z} that the
 group von Neumann algebra $\caln(G)$ behaves like the ring of integers $\bbZ$, provided
one ignores the properties integral domain and Noetherian. Given a finitely
generated abelian group $M$, the $\bbZ$-module $M/\tors(M)$ is finitely generated free,
there is a $\bbZ$-isomorphism $ M \cong M/\tors(M) \oplus \tors(M)$  
and the rank as abelian group of $M$ is $\dim_{\bbQ}(\bbQ \otimes_{\bbZ} M)$
and of $\tors(M)$ is $0$.
In analogy, given a finitely
generated $\caln(G)$-module $M$, then $\bfP M := M/\bfT M$ is a finitely generated projective$\caln(G)$-module,
there is an $\caln(G)$-isomorphism $ M \cong \bfP M \oplus \tors(M)$  and we get
$\dim_{\caln(G)}(M) = \dim_{\calu(G)}(\calu(G) \otimes_{\caln(G)} M)$ and
$\dim_{\caln(G)}(\bfT M) = 0$. Define the so called \emph{capacity}%
\index{capacity}
$c(M) \in [0,\infty] \cup \{0^-\}$ of a finitely presented $\caln(G)$-module $M$ by 
$$c(M) ~ = ~  \left\{\begin{array}{lll} 
\frac{1}{\alpha(M)} & & \text{ if } \alpha(M) \in (0,\infty);
\\
\infty & & \text{ if } \alpha(M) = 0;
\\
0 & & \text{ if } \alpha(M) = \infty;
\\
0^- & & \text{ if } \alpha(M) = \infty^+.
\end{array}\right.$$
Then the  capacity $c(M)$ contains the same information as $\alpha(M)$ and
corresponds under the dictionary between $\bbZ$ and $\caln(G)$ to the order of the finite group $\tors(M)$. 
Notice for a finitely presented
$\caln(G)$-module $M$ that $M = 0$ is true if and only if
both $\dim_{\caln(G)}(M) = 0$ and $c(M) = 0^-$ hold.
The capacity is at least subadditive, i.e.
for an exact sequence $1 \to M_0 \to M_1 \to M_2 \to 0$ of finitely presented
$\caln(G)$-modules we have $c(M_1) \le c(M_0) + c(M_2)$ (with the obvious interpretation
of $+$ and $\le$). In particular we get $c(M) \le c(N)$ 
for an inclusion of finitely presented $\caln(G)$-modules $M \subseteq N$. 
\em
\end{remark}

\begin{remark}[Extension to arbitrary $\caln(G)$-modules and $G$-spaces]
\label{rem: Extension to arbitrary caln(G)-modules and G-spaces}
\em
The algebraic approach presented above has been independently developed 
in \cite{Farber(1996)} and \cite{Lueck(1997a)}. The notion of capacity has been extended 
by L\"uck-Reich-Schick \cite{Lueck-Reich-Schick(1999)} to so called
cofinal-measurable $\caln(G)$-modules, i.e. $\caln(G)$-modules such that each finitely
generated $\caln(G)$-submodule is a quotient of a finitely presented $\caln(G)$-module with
trivial von Neumann dimension. This allows to define
Novikov-Shubin invariants for arbitrary $G$-spaces and also for arbitrary groups $G$.
\em
\end{remark}


\subsection{Basic Properties of Novikov-Shubin Invariants}
\label{subsec: Basic Properties of Novikov-Shubin Invariants}

We briefly list some properties of Novikov-Shubin invariants.
The proof of the following theorem can be found in \cite[Theorem 2.55]{Lueck(2002)} and 
\cite[Lemma 13.45]{Lueck(2002)}.

\begin{theorem}[Novikov-Shubin invariants]\
\index{Theorem!Novikov-Shubin invariants}
\label{the: properties of cellular Novikov-Shubin invariants}
\begin{enumerate}
\item Homotopy invariance
\label{the: properties of cellular Novikov-Shubin invariants:
homotopy invariance}
\\[1mm]
Let $f\colon  X \rightarrow Y$ be a $G$-map of free
$G$-$CW$-complexes of finite type. Suppose that the map 
$H_p(f;\bbC)\colon H_p(X;\bbC) \to H_p(Y;\bbC)$
induced on homology with complex coefficients
is an isomorphism for $p \le d-1$.
Then we get
$$\alpha_p(X;\caln(G)) ~ = ~ \alpha_p(Y;\caln(G)) \quad \text{ for }p \le d.$$

In particular
we get $\alpha_p(X;\caln(G)) ~ = ~ \alpha_p(Y;\caln(G))$ for all $p \ge 0$
if $f$ is a weak homotopy equivalence;

\item \label{the: properties of cellular Novikov-Shubin invariants:
Poincar'e duality}
Poincar\'e duality
\\[1mm]
Let $M$ be a  cocompact free proper $G$-manifold of dimension $n$
which is orientable. Then $\alpha_p(M;\caln(G)) ~ = ~ \alpha_{n+1-p}(M,\partial
M;\caln(G))$ for $p \ge 1$;

\item First Novikov-Shubin invariant
\label{the: properties of cellular Novikov-Shubin invariants:
First Novikov-Shubin invariant}\\[1mm]
Let $X$ be a connected free $G$-$CW$-complex of finite type. Then
$G$ is finitely generated and
\begin{enumerate}

\item
\label{the: properties of cellular Novikov-Shubin invariants:
First Novikov-Shubin invariant: 1}
$\alpha_1(X)$ is
finite if and only if  $G$ is infinite and virtually nilpotent. In this case
$\alpha_1(X)$ is the growth rate of $G$;

\item \label{the: properties of cellular Novikov-Shubin invariants:
First Novikov-Shubin invariant: 2}
$\alpha_1(X)$
is $\infty^+$ if and only if $G$ is finite or non-amenable;

\item \label{the: properties of cellular Novikov-Shubin invariants:
First Novikov-Shubin invariant: 3}
$\alpha_1(X)$
is $\infty$ if and only if $G$ is amenable and not virtually nilpotent;

\end{enumerate}

\item Restriction to subgroups of finite index
\label{the: properties of cellular Novikov-Shubin invariants: restriction}
\\[1mm]
Let $X$ be a free $G$-$CW$-complex of finite type and $H \subseteq G$
a subgroup of finite index.
Then $\alpha_p(X;\caln(G)) ~ = ~ \alpha_p(\res_G^H X;\caln(H))$ holds  for $p \ge 0$;

\item \label{the: properties of cellular Novikov-Shubin invariants: normal finite sub.} 
Extensions with finite kernel\\[1mm]
Let $1 \to H \to G \to Q \to 1$ be an extension of groups such that $H$ is finite.
Let $X$ be a free $Q$-$CW$-complex of finite type. Then we get
$\alpha_p(p^*X;\caln(G))  ~ = ~ \alpha_p(X;\caln(Q))$ for all $p \ge 1$;

\item \label{the: properties of cellular Novikov-Shubin invariants:
induction}
Induction
\\[1mm]
Let $H$ be a subgroup of $G$ and let
$X$ be a free $H$-$CW$-complex of finite type.
Then $\alpha_p(G \times_H X;\caln(G)) ~ = ~ \alpha_p(X;\caln(H))$
holds for all $p \ge 1$.

\end{enumerate}
\end{theorem}

A product formula and a formula for connected sums can also be found in
\cite[Theorem 2.55]{Lueck(2002)}. If $X$ is a finite $G$-$CW$-complex
such that $b_p^{(2)}(X;\caln(G)) = 0$ for $p \ge 0$ and
$\alpha_p(X;\caln(G)) > 0$ for $p \ge 1$, then $X$ is
$\det$-$L^2$-acyclic \cite[Theorem 3.93 (7)]{Lueck(2002)}.


\subsection{Computations of Novikov-Shubin Invariants}
\label{subsec: Computations of Novikov-Shubin Invariants}

\begin{example}[Novikov-Shubin invariants of $\widetilde{T^n}$]
\label{exa: Novikov-Shubin invariants of widetilde{T^n}} \em
The product formula can be used to show
$\alpha_p(\widetilde{T^n}) = n$ if $1 \le p \le n$,
 and $\alpha_p(\widetilde{T^n}) = \infty^+$ otherwise (see \cite[Example 2.59]{Lueck(2002)}.)
\em
\end{example}

\begin{example}[Novikov-Shubin invariants for finite groups]
 \label{exa: Novikov-Shubin invarioants for finite groups} \em
If $G$ is finite, then $\alpha_p(X;\caln(G)) = \infty^+$ for each $p
\ge 1$ and $G$-$CW$-complex $X$ of finite type. This follows from 
Example ~\ref{spectral density function for finite G}. This shows that
the Novikov-Shubin invariants are interesting only for infinite groups
$G$ and have no classical analogue in contrast to $L^2$-Betti numbers
and $L^2$-torsion.
\em
\end{example}

\begin{example}[Novikov-Shubin invariants for $G = \bbZ$]
\label{exa: Novikov-Shubin invariants for G = Z} \em
Let $X$ be a free $\bbZ$-$CW$-complex of finite type.
Since $\bbC[\bbZ]$ is a principal ideal domain, we get $\bbC[\bbZ]$-isomorphisms
\begin{eqnarray*}
H_p(X;\bbC)
& \cong &
\bbC[\bbZ]^{n_p} \oplus \left(\bigoplus_{i_p=1}^{s_p}
\bbC[\bbZ]/((z-a_{p,i_p})^{r_{p,i_p}})\right)
\end{eqnarray*}
for $a_{p,i_p} \in \bbC$
and $n_p, s_p, r_{p,i_p} \in \bbZ$ with $n_p,s_p \ge 0$
and $ r_{p,i_p} \ge 1$,
where $z$ is a fixed generator of $\bbZ$. Then we get from \cite[Lemma 2.58]{Lueck(2002)}
\begin{eqnarray*}
b^{(2)}_p(X;\caln(\bbZ))
& = &
n_p.
\end{eqnarray*}
If $s_p \ge 1$ and $\{i_p = 1,2 \ldots, s_p, |a_{p,i_p}| = 1\} \not= \emptyset$, then
\begin{eqnarray*}
\alpha_{p+1}(X;\caln(\bbZ)) & = &
\min\left\{\frac{1}{r_{p,i_p}}\mid i_p = 1,2 \ldots, s_p, |a_{p,i_p}| = 1\right\},
\end{eqnarray*}
and otherwise
\begin{eqnarray*}
\alpha_{p+1}(X;\caln(\bbZ)) & = & \infty^+.
\end{eqnarray*}
\em
\end{example}

\begin{remark}[Novikov-Shubin invariants and $S^1$-actions]
\label{Novikov-Shubin invariants and S^1-actions} \em
Under the conditions of Theorem~\ref{the: $S^1$-actions and general $L^2$-Betti numbers} 
and of Theorem~\ref{the: fixed point free S^1-actions on Eilenberg MacLane spaces on
aspherical closed manifolds and L^2-Betti numbers}
one can show $\alpha_p(\widetilde{X}) \ge 1$ for all $p \ge 1$
(see \cite[Theorem 2.61 and Theorem 2.63]{Lueck(2002)}). 
\em
\end{remark}

\begin{remark}[Novikov-Shubin invariants of symmetric spaces of non-compact type]
\label{the: Nov.Shu. inv. of sym. sp.} \em
The Novikov-Shubin invariants of symmetric spaces of non-compact type 
with cocompact free $G$-action are
computed by Olbricht \cite[Theorem 1.1]{Olbrich(2000)}, the result is also stated
in \cite[Section 5.3]{Lueck(2002)}. \em
\end{remark}

\begin{remark}[Novikov-Shubin invariants of universal coverings of $3$-manifolds]
\label{the: Nov.Shu. inv. of uni. cov. of .3-mani} \em
Partial results about the computation of the Novikov-Shubin invariants
of universal coverings of compact orientable $3$-manifolds can be found in
\cite{Lott-Lueck(1995)} and \cite[Theorem 4.2]{Lueck(2002)}.
\em
\end{remark}


\subsection{Open Conjectures about Novikov-Shubin invariants}

\label{subsec: Open Conjectures about Novikov-Shubin invariants}

The following conjecture is taken from
\cite[Conjecture 7.1 on page 56]{Lott-Lueck(1995)}.

\begin{conjecture}{\bf (Positivity and rationality of Novikov-Shubin invariants).}%
\index{Conjecture!Positivity of Novikov-Shubin invariants}
\label{con: rationality and positivity of Novikov-Shubin-invariants}
Let $G$ be a group. Then for any
free $G$-$CW$-complex $X$ of finite type its Novikov-Shubin invariants
$\alpha_p(X)$ are positive rational numbers unless they are
$\infty$ or $\infty^+$.
\end{conjecture}

This conjecture is equivalent to the statement that for every finitely
presented $\bbZ G$-module $M$ the Novikov-Shubin invariant of
$\caln(G) \otimes_{\bbZ G} M$  is a positive rational number, $\infty$
or $\infty^+$.

Here is some evidence for Conjecture
\ref{con: rationality and positivity of Novikov-Shubin-invariants}.
Unfortunately, all the evidence comes from computations, no convincing
conceptual reason is known. Conjecture
\ref{con: rationality and positivity of  Novikov-Shubin-invariants}
is true for $G = \bbZ$ by the explicit computation appearing in 
Example~\ref{exa: Novikov-Shubin invariants for G = Z}.
Conjecture \ref{con: rationality and positivity of Novikov-Shubin-invariants}
is true for virtually abelian $G$ by
\cite[Proposition 39 on page 494]{Lott(1992a)}.
Conjecture \ref{con: rationality and positivity of Novikov-Shubin-invariants}
is also true for a free group $G$.
Details of the proof appear in the Ph.D. thesis of Roman Sauer
\cite{Sauer(2002)} following ideas of Voiculescu. The essential
ingredients are non-commutative power series and the question
whether they are algebraic or rational. 
All the computations mentioned above are compatible and give evidence for Conjecture
\ref{con: rationality and positivity of Novikov-Shubin-invariants}.

\begin{conjecture}[Zero-in-the-spectrum Conjecture]
\label{con: Zero-in-the-spectrum Conjecture}
Let $G$ be a group such that $BG$ has a closed manifold as model. Then there is $p \ge
0$ with $H_p^G(EG;\caln(G)) \not= 0.$
\end{conjecture}

\begin{remark}[Original zero-in-the-spectrum Conjecture]
\label{rem: original zero-in-the-spectrum Conjecture}
\em The original zero-in-the-spectrum Conjecture,
which appears for the first time
in Gromov's article \cite[page 120]{Gromov(1986a)}, says the following:
Let $\widetilde{M}$ be a complete Riemannian manifold.
Suppose that $\widetilde{M}$ is the universal covering of an
aspherical closed Riemannian manifold $M$ (with the Riemannian metric
coming from $M$). Then for some $p \ge 0$ zero is  in the
spectrum of the minimal closure
$$(\Delta_p)_{\min}\colon  \dom\left((\Delta_p)_{\min}\right) \subseteq
L^2\Omega^p(\widetilde{M}) \to L^2\Omega^p(\widetilde{M})$$ of the
Laplacian acting on smooth $p$-forms on $\widetilde{M}$.

It follows from \cite[Lemma 12.3]{Lueck(2002)} that this formulation
is equivalent to the homological algebraic formulation appearing in 
Conjecture~\ref{con: Zero-in-the-spectrum Conjecture}. \em
\end{remark}

\begin{remark}[Status of the zero-in-the-spectrum Conjecture] 
\label{rem: Status of the zero-in-the-spectrum Conjecture}\em
The zero-in-the-spectrum Conjecture is true for $G$
if there is a closed manifold model for $BG$ which is K\"ahler hyperbolic
\cite{Gromov(1991)}, or whose universal covering is
hyperEuclidean \cite{Gromov(1986a)} or is uniformly contractible with finite asymptotic
dimension \cite{Yu(1998a)}. The zero-in-the-spectrum Conjecture is true for $G$
if the strong Novikov Conjecture holds for $G$ \cite{Lott(1996b)}. 
More information about zero-in-the-spectrum Conjecture 
can be found for instance in \cite{Lott(1996b)} and
\cite[Section 12]{Lueck(2002)}.
\em
\end{remark}

\begin{remark}[Variations of the zero-in-the-spectrum Conjecture]
\label{rem: Variations of the zero-in-the-spectrum Conjecture}
\em
One may ask whether one can weaken the condition in 
Conjecture~\ref{con: Zero-in-the-spectrum Conjecture} that $BG$ has a
closed manifold model to the condition that there is a finite
$CW$-complex model for $BG$. This would rule out Poincar\'e duality
from the picture. Or one could only require that $BG$ is of finite
type. Without any finiteness conditions on $G$ 
Conjecture~\ref{con: Zero-in-the-spectrum Conjecture}  is not true in general. For instance
$H_p^G(EG;\caln(G)) = 0$ holds for all $p \ge 0$ if $G$ is $\prod_{i=1}^{\infty} \bbZ \ast \bbZ$.

The condition aspherical cannot be dropped. Farber and Weinberger 
\cite{Farber-Weinberger(2001)} proved the existence of a closed non-aspherical
manifold $M$ with fundamental group $\pi$ a product of three copies of $\bbZ \ast \bbZ$
such that $H_p^{\pi}(\widetilde{M};\caln(\pi))$ vanishes for all $p\ge 0$.
Later Higson-Roe-Schick \cite{Higson-Roe-Schick(2001)} proved that one can find for every  finitely
presented group $\pi$, for which $H_p^{\pi}(E\pi;\caln(\pi)) = 0$ holds for $p = 0,1,2$,
a closed manifold $M$ with $\pi$ as fundamental group such that 
$H_p^{\pi}(\widetilde{M};\caln(\pi))$ vanishes for all $p\ge 0$.
\em
\end{remark}

\begin{remark}[Novikov-Shubin invariants and quasi-isometry] \em
\label{rem: Novikov-Shubin invariants and quasi-isometry} 
Since $\alpha_1(\bbZ^n) = n$ for $n \ge 1$, the Novikov-Shubin invariants are not
invariant under measure equivalence. It is not known whether they are invariant under
quasi-isometry. At least it is known that two quasi-isometric amenable groups $G_1$ and $G_2$
which possess finite models for $BG_1$ and $BG_2$ have the same
Novikov-Shubin invariants \cite{Sauer(2002)}. Compare also
Theorem~\ref{the: Gaboriau's result on measure equivalent groups},
Remark~\ref{rem: ^L2-Betti numbers and quasi-isometry} and
Conjecture~\ref{con: Measure equivalence and L^2-torsion}.
\em
\end{remark}


\typeout{--------------------   Section 13 --------------------------}

\section{A Combinatorial Approach to $L^2$-Invariants}
\label{sec: A Combinatorial Approach to $L^2$-Invariants}

In this section we want to give a more combinatorial approach
to the $L^2$-invariants such as $L^2$-Betti numbers, Novikov-Shubin
invariants and $L^2$-torsion. The point is that it is in general very
hard to compute the spectral density function of an $\caln(G)$-map 
$f \colon \caln(G)^m \to \caln(G)^n$. However in the geometric
situation these morphisms are induced by
matrices over the integral group ring $\bbZ G$. We  want to exploit
this information to get an algorithm which produces
a sequence of rational numbers converging to the $L^2$-Betti number
or the  $L^2$-torsion in question.

Let $A \in M(m,n;\bbC G)$ be an $(m,n)$-matrix over $\bbC G$.
It induces by
right multiplication an $\caln(G)$-homomorphism 
$r_A\colon  \caln(G)^m \to \caln(G)^n$.
We define an involution  of rings on $\bbC G$ by sending
$\sum_{g \in G} \lambda_g \cdot g$ to 
$\sum_{g \in G} \overline{\lambda_g} \cdot g^{-1}$,
where $\overline{\lambda_g}$ is the complex conjugate of $\lambda_g$.
Denote by $A^*$ the $(m,n)$-matrix obtained from $A$ by
transposing and applying the involution above to each entry. Define the
\emph{$\bbC G$-trace}%
\index{trace!over the complex group ring}
of an element
$u = \sum_{g \in G} \lambda_g \cdot g \in \bbC G$ by
the complex number $\tr_{\bbC G}(u) := \lambda_e$
for $e$ the unit element in $G$. This extends to a trace of 
square $(n,n)$-matrices $A$ over $\bbC G$ by
\begin{eqnarray}
\tr_{\bbC G}(A)%
\indexnotation{tr_{cc G}(A)} & := & \sum_{i=1}^n \tr_{\bbC G}(a_{i,i}) \quad \in \bbC.
\label{definition of tr_{cc G}(A)}
\end{eqnarray}
We get directly from the definitions that the $\bbC G$-trace
$\tr_{\bbC G}(u)$ for $u \in \bbC G$ agrees with the von Neumann trace
$\tr_{\caln(G)}(u)$ introduced
in Definition \ref{def: trace of the group von Neumann algebra}.

Let $A \in M(m,n;\bbC G)$ be an $(m,n)$-matrix over $\bbC G$. In the
sequel let $K$ be any positive real number satisfying
$K~ \ge ~ ||r_A^{(2)}||$,
where $||r_A^{(2)}||$ is the operator norm of the
bounded  $G$-equivariant operator $r_A^{(2)} : l^2(G)^m \to l^2(G)^n$
induced by right multiplication with $A$.
For $u = \sum_{g \in G} \lambda_g \cdot g \in \bbC G$ define
$||u||_1$ by
$\sum_{g \in G} |\lambda_g|$.
Then a possible choice for $K$ is 
$$K = \sqrt{(2n-1)m} \cdot \max\left\{ ||a_{i,j}||_1
\mid 1 \le i \le n, 1 \le j \le m\right\}.$$

\begin{definition} \label{def: characteristic sequence}
The \emph{characteristic sequence}%
\index{characteristic sequence}
of a matrix $A \in M(m,n;\bbC G)$ and a non-negative real number $K$
satisfying $K \ge ||r_A^{(2)}||$ is the sequence of real numbers given by
$$c(A,K)_p := \tr_{\bbC G}\left(\left(1  - K^{-2}\cdot
AA^*\right)^p\right).$$
\end{definition}

We have defined $\dim_{\caln(G)}(\ker(r_A))$ in Definition
\ref{def: Von Neumann dimension for arbitrary caln(G)-modules}
and $\det_{\caln(G)}(r_A)$ in Definition
\ref{def: Fuglede-Kadison determinant}. The proof of the following
result can be found in \cite{Lueck(1994a)} or \cite[Theorem 3.172]{Lueck(2002)}.

\begin{theorem}{\bf (Combinatorial computation of $L^2$-invariants).}\\
\index{Theorem!Combinatorial computation of $L^2$-invariants}
\label{L^2-invariants and characteristic sequences}
Let $A \in M(m,n;\bbC G)$ be an $(m,n)$-matrix over $\bbC G$. Let $K$ be a
positive real number satisfying
$K \ge ||r_A^{(2)}||$. Then:

\begin{enumerate}

\item \label{L^2-invariants and characteristic sequences:  monotonicity}
Monotony\\[1mm]
The characteristic sequence $\left(c(A,K)_p\right)_{p \ge 1}$ is a monotone decreasing
sequence of non-negative real numbers;

\item \label{L^2-invariants and characteristic sequences: b^{(2)}}
Dimension of the kernel\\[1mm]
We have
$$\dim_{\caln(G)}(\ker(r_A)) ~ = ~ \lim_{p \to \infty} c(A,K)_p;$$

\item \label{L^2-invariants and characteristic sequences: beta}
Novikov-Shubin invariants of the cokernel\\[1mm]
Define $\beta (A) \in [0,\infty]$ by
$$\beta (A):=   \sup\left\{ \beta \in [0,\infty) ~ \left|~ 
\lim_{p \to \infty} p^{\beta} \cdot \left(c(A,K)_p -
\dim_{\caln(G)}(\ker(r_A)) \right) \right.= 0\right\}.$$
If $\alpha(\coker(r_A)) < \infty$, then $\alpha(\coker(r_A)) \le \beta(A)$
and if $\alpha(\coker(r_A)) \in \{\infty, \infty^+\}$,
then $\beta(A) = \infty$;

\item Fuglede-Kadison determinant\\[1mm]
\label{L^2-invariants and characteristic sequences: det}
The sum of positive real numbers
$$\sum_{p=1}^{\infty}
\frac{1}{p} \cdot \left(c(A,K)_p - \dim_{\caln(G)}(\ker(r_A))\right) $$
converges if and only if $r_A$ is of determinant class
and in this case
\begin{multline*}
\ln(\det(r_A)) 
~ = ~ (n - \dim_{\caln(G)}(\ker(r_A))) \cdot \ln(K)
\\
~ - ~
\frac{1}{2} \cdot \sum_{p=1}^{\infty} \frac{1}{p} \cdot
\left(c(A,K)_p - \dim_{\caln(G)}(\ker(r_A)) \right);
\end{multline*}
\item   \label{L^2-invariants and characteristic sequences: speed}
Speed of convergence\\[1mm]
Suppose $\alpha(\coker(r_A)) > 0$. Then $r_A$ is of determinant class.
Given a real number $\alpha$ satisfying $0 < \alpha < \alpha(\coker(r_A))$,
there is a real number $C$ such that we have for all $L \ge 1$
$$0 \le c(A,K)_L - \dim_{\caln(G)}(\ker(r_A))) \le \frac{C}{L^\alpha}$$
and
\begin{eqnarray*}
 0 & \le &
- \ln(\det(r_A)) + (n - \dim_{\caln(G)}(\ker(r_A)))\cdot \ln(K)
\\
& & \hspace{10mm}
- \frac{1}{2} \cdot \sum_{p=1}^{L} \frac{1}{p}\cdot \left(c(A,K)_p - \dim_{\caln(G)}(\ker(r_A))\right)
~ \le ~
\frac{C}{L^{\alpha}}.
\end{eqnarray*}
\end{enumerate}
\end{theorem}

\begin{remark}[Vanishing of $L^2$-Betti numbers and the Atiyah Conjecture]
\label{rem: Vanishing of L^2-Betti numbers and the Atiyah Conjecture} \em
Suppose that the Atiyah Conjecture~\ref{con: Atiyah Conjecture} is satisfied for
$(G,d,\bbC)$. If we want to show the vanishing of $\dim_{\caln(G)}(\ker(r_A))$,
it suffices to  show that for some $p \ge 0$ we have 
$c(A,K)_p < \frac{1}{d}$. It is possible that a computer program spits out such a value
after a reasonable amount of calculation time.
\em
\end{remark}


\typeout{--------------------   Section 14 --------------------------}

\section{Miscellaneous}
\label{sec: Miscelleaneous}

The analytic aspects of $L^2$-invariants are also very interesting. We have
already mentioned that $L^2$-Betti numbers were originally defined by
Atiyah \cite{Atiyah(1976)} in context with his $L^2$-index theorem. 
Other $L^2$-invariants are the $L^2$-Eta-invariant and the
$L^2$-Rho-invariant (see Cheeger-Gromov 
\cite{Cheeger-Gromov(1985)}, \cite{Cheeger-Gromov(1985a)}).
The $L^2$-Eta-invariant  appears in the $L^2$-index theorem for manifolds with boundary
due to Ramachandran \cite{Ramachandran(1993)}.
These index theorems have generalizations to a $C^*$-setting due to
Mi{\v{s}}{\v{c}}enko-Fomenko \cite{Mishchenko-Fomenko(1979)}.
There is also an $L^2$-version of the signature. It plays an important
role in the work of Cochran, Orr and Teichner
\cite{Cochran-Orr-Teichner(1999)} who show that
there are non-slice knots in $3$-space whose Casson-Gordon
invariants are all trivial. Chang and Weinberger 
\cite{Chang-Weinberger(2003)} show using
$L^2$-invariants that for a closed oriented smooth
manifold $M$ of dimension $4k+3$ for $k \ge 1$ whose fundamental
group has torsion there are infinitely many smooth manifolds 
which are homotopy equivalent to M (and even simply and tangentially homotopy
equivalent to $M$) but not homeomorphic to M.
The $L^2$-cohomology has also been investigated for complete non-necessarily
compact Riemannian manifolds without a group action. 
For instance algebraic and arithmetic varieties have been studied.
In particular, the Cheeger-Goresky-MacPherson Conjecture
\cite{Cheeger-Goresky-McPherson(1982)}
and the Zucker Conjecture \cite{Zucker(1983)}
have created a lot of activity. They link the $L^2$-cohomology of the
regular part with the intersection homology of an algebraic variety.

Finally we mention other survey articles which deal with
$L^2$-invariants: \cite{Eckmann(2000)}, \cite{Gaboriau(2001b)},
\cite[Section 8]{Gromov(1993)},  \cite{Lott(1996b)},
\cite{Lueck(1997c)}, \cite{Lueck(2001c)}, \cite{Lueck(2002a)}, \cite{Mathai(1998)} and
\cite{Pansu(1996a)}.


\typeout{-------------------- References -------------------------------}


\typeout{-------------------------- references  --------------------------}


\typeout{-------------------- Notation -------------------------------}
\twocolumn
\section*{Notation}
\addcontentsline{toc}{section}{Notation}

\entry{$BG$}{BG}\\
\entry{$b^{(2)}_p(G)$}{b_p^(2)(G) for arbitrary G}\\
\entry{$b_p^{(2)}(\underline{G})$}{b_p^(2)(underline G)}\\
\entry{$b_p^{(2)}(G \action X)$}{b_p^(2)(G action X)}\\
\entry{$b_p^{(2)}(X;\caln(G))$}{b_p^{(2)}(X;caln(G)) for arbitrary X}\\
\entry{$\bbC \underline{G}$}{C underline G}\\
\entry{$C_*^{(2)}(X)$}{C_*^(2)(X)}\\
\entry{${\det}_{\caln(G )}(f)$}{det_{caln(G)}(f),based}\\
\entry{$\dim_{\caln(G)}(M)$}{dim_{caln(G)}(M)}\\
\entry{$\dim^u_{\caln(G)}(P)$}{dim^u_caln(G)}\\
\entry{$\dim_{\calu(G)}(M)$}{dim_{calu(G)}(M)}\\
\entry{$EG$}{EG}\\
\entry{$F^d_g$}{F_g^d}\\
\entry{$F_f$}{F_f}\\
\entry{$\frk(M)$}{frk(M)}\\
\entry{$G_f$}{G_f}\\
\entry{$G_0(R)$}{G_0(R)}\\
\entry{$H_p^G(X;\caln(G))$}{H_p^G(X;caln(G))}\\
\entry{$H_p^{(2)}(X;l^2(G))$}{H_p^{(2)}(X;l^2(G))}\\
\entry{$h^{(2)}(G)$}{h(2)(X;caln(G))}\\
\entry{$h^{(2)}(X;\caln(G))$}{h2G}\\
\entry{$i_*M$}{i_*M}\\
\entry{$I(X)$}{I(X)}\\
\entry{$K_0(R)$}{K_0(R)}\\
\entry{$l^2(G)$}{l^2(G)}\\
\entry{$L^2(T^n)$}{L^2(T^n)}\\
\entry{$L^{\infty}(T^n)$}{L^{infty}(T^n)}\\
\entry{$m(X;G)$}{m(X;caln(G))}\\
\entry{$\overline{M}$}{overline{M}}\\
\entry{$||M||$}{||M||}\\
\entry{$R \ast_c G$}{R ast_c G}\\
\entry{$\bfP M$}{bfP M}\\
\entry{$\bfT M$}{bfT M}\\
\entry{$T_f$}{T_f}\\
\entry{$\tr_{\bbC G}(A)$}{tr_{cc G}(A)}\\
\entry{$\tr_{\caln(G)}$}{tr_{caln(G)}}\\
\entry{$\tr^u_{\caln(G)}$}{tr^u_caln(G)}\\
\entry{$\vol(M)$}{vol(M)}\\
\entry{$\Wh(G)$}{Wh(G)}\\
\entry{$\widetilde{X}$}{widetilde{X}}\\
\entry{$\alpha(M)$}{alpha(M)}\\
\entry{$\alpha_p^{\Delta}(M;\caln(G))$}{alpha_p^Delta(M)}\\
\entry{$\chi^{(2)}(G)$}{chi^{(2)}(G)}\\
\entry{$\chi^{(2)}(X;\caln(G))$}{chi^{(2)}(X;caln(G))}\\
\entry{$\chi_{\virt}(X)$}{chivirtuell} \\
\entry{$\cala$}{cala}\\
\entry{$\calb(H)$}{calb(H)}\\
\entry{$\calb_{d}$}{calb_{d}}\\
\entry{$\calc$}{calc}\\
\entry{$\cald$}{cald}\\
\entry{$\caleam$}{caleam}\\
\entry{$\calfin(G)$}{calfin(G)}\\
\entry{$\caln(G)$}{caln(G)}\\
\entry{$\caln(i)$}{induced ring homomorphism caln(i)}\\
\entry{$\calr(G \action X)$}{calr(G action X)}\\
\entry{$\calu(G)$}{calu(G)}\\
\entry{$[0,\infty]$}{[0,infty]}\\
\entry{$\frac{1}{|\calfin(G)|}\bbZ$}{frac{1}{midcalfin(G)mid}zz}
\onecolumn


\typeout{-------------------- Index ---------------------------------}

\flushbottom
\addcontentsline{toc}{section}{Index}
\printindex                                  


\addcontentsline{toc}{section}{References}

\begin{thebibliography}{99}

\bibitem{Ara-Goldstein(1993)}
P.~Ara and D.~Goldstein.
\newblock A solution of the matrix problem for {R}ickart ${C}\sp *$-algebras.
\newblock {\em Math. Nachr.}, 164:259--270, 1993.

\bibitem{Atiyah(1976)}
M.~F. Atiyah.
\newblock Elliptic operators, discrete groups and von {N}eumann algebras.
\newblock {\em Ast\'erisque}, 32-33:43--72, 1976.

\bibitem{Ballmann-Bruening(2001)}
W.~Ballmann and J.~Br{\"u}ning.
\newblock On the spectral theory of manifolds with cusps.
\newblock {\em J. Math. Pures Appl. (9)}, 80(6):593--625, 2001.

\bibitem{Bekka-Valette(1997)}
M.~E.~B. Bekka and A.~Valette.
\newblock Group cohomology, harmonic functions and the first ${L}\sp 2$-{B}etti
  number.
\newblock {\em Potential Anal.}, 6(4):313--326, 1997.

\bibitem{Berberian(1972)}
S.~K. Berberian.
\newblock {\em Baer *-rings}.
\newblock Springer-Verlag, New York, 1972.
\newblock Die Grundlehren der mathematischen Wissenschaften, Band 195.

\bibitem{Borel(1985)}
A.~Borel.
\newblock The ${L}\sp 2$-cohomology of negatively curved {R}iemannian symmetric
  spaces.
\newblock {\em Ann. Acad. Sci. Fenn. Ser. A I Math.}, 10:95--105, 1985.

\bibitem{Brin-Squier(1985)}
M.~G. Brin and C.~C. Squier.
\newblock Groups of piecewise linear homeomorphisms of the real line.
\newblock {\em Invent. Math.}, 79(3):485--498, 1985.

\bibitem{Brown-Geoghegan(1984)}
K.~S. Brown and R.~Geoghegan.
\newblock An infinite-dimensional torsion-free ${\rm {f}{p}}\sb{\infty }$
  group.
\newblock {\em Invent. Math.}, 77(2):367--381, 1984.

\bibitem{Burghelea-Friedlander-Kappeler-McDonald(1996a)}
D.~Burghelea, L.~Friedlander, T.~Kappeler, and P.~McDonald.
\newblock Analytic and {R}eidemeister torsion for representations in finite
  type {H}ilbert modules.
\newblock {\em Geom. Funct. Anal.}, 6(5):751--859, 1996.

\bibitem{Cannon-Floyd-Parry(1996)}
J.~W. Cannon, W.~J. Floyd, and W.~R. Parry.
\newblock Introductory notes on {R}ichard {T}hompson's groups.
\newblock {\em Enseign. Math. (2)}, 42(3-4):215--256, 1996.

\bibitem{Chang-Weinberger(2003)}
S.~Chang and S.~Weinberger.
\newblock On invariants of {H}irzebruch and {C}heeger.
\newblock Preprint, 2003.

\bibitem{Cheeger-Goresky-McPherson(1982)}
J.~Cheeger, M.~Goresky, and R.~MacPherson.
\newblock ${L}\sp{2}$-cohomology and intersection homology of singular
  algebraic varieties.
\newblock In {\em Seminar on Differential Geometry}, pages 303--340. Princeton
  Univ. Press, Princeton, N.J., 1982.

\bibitem{Cheeger-Gromov(1985)}
J.~Cheeger and M.~Gromov.
\newblock Bounds on the von {N}eumann dimension of ${L}\sp 2$-cohomology and
  the {G}auss-{B}onnet theorem for open manifolds.
\newblock {\em J. Differential Geom.}, 21(1):1--34, 1985.

\bibitem{Cheeger-Gromov(1985a)}
J.~Cheeger and M.~Gromov.
\newblock On the characteristic numbers of complete manifolds of bounded
  curvature and finite volume.
\newblock In {\em Differential geometry and complex analysis}, pages 115--154.
  Springer-Verlag, Berlin, 1985.

\bibitem{Cheeger-Gromov(1986)}
J.~Cheeger and M.~Gromov.
\newblock ${L}\sb 2$-cohomology and group cohomology.
\newblock {\em Topology}, 25(2):189--215, 1986.

\bibitem{Cochran-Orr-Teichner(1999)}
T.~Cochran, K.~Orr, and P.~Teichner.
\newblock Knot concordance, {W}hitney towers and ${L^2}$-signatures.
\newblock Preprint, to appear in \emph{Ann. of Math.}, 1999.

\bibitem{Cohen(1973)}
M.~M. Cohen.
\newblock {\em A course in simple-homotopy theory}.
\newblock Springer-Verlag, New York, 1973.
\newblock Graduate Texts in Mathematics, Vol. 10.

\bibitem{Connes-Shlyakhtenko(2003)}
A.~Connes and D.~Shlyakhtenko.
\newblock {$L^2$}-homology for von {N}eumann algebras.
\newblock Preprint, 2003.

\bibitem{Davis-Leary(2001z)}
M.~Davis and I.~Leary.
\newblock The ${L}^2$-cohomology of {A}rtin groups.
\newblock Preprint, to appear in LMS, 2001.

\bibitem{Davis-Okun(2001)}
M.~W. Davis and B.~Okun.
\newblock Vanishing theorems and conjectures for the $\ell\sp 2$-homology of
  right-angled {C}oxeter groups.
\newblock {\em Geom. Topol.}, 5:7--74 (electronic), 2001.

\bibitem{delaHarpe(2000)}
P.~de~la Harpe.
\newblock {\em Topics in geometric group theory}.
\newblock University of Chicago Press, Chicago, IL, 2000.

\bibitem{Dicks-Schick(2002)}
W.~Dicks and T.~Schick.
\newblock The spectral measure of certain elements of the complex group ring of
  a wreath product.
\newblock {\em Geom. Dedicata}, 93:121--137, 2002.

\bibitem{Dixmier(1981)}
J.~Dixmier.
\newblock {\em Von {N}eumann algebras}.
\newblock North-Holland Publishing Co., Amsterdam, 1981.
\newblock With a preface by E. C. Lance, translated from the second French
  edition by F. Jellett.

\bibitem{Dodziuk(1977)}
J.~Dodziuk.
\newblock De {R}ham-{H}odge theory for ${L}\sp{2}$-cohomology of infinite
  coverings.
\newblock {\em Topology}, 16(2):157--165, 1977.

\bibitem{Dodziuk(1979)}
J.~Dodziuk.
\newblock ${L}\sp{2}$\-harmonic forms on rotationally symmetric {R}iemannian
  manifolds.
\newblock {\em Proc. Amer. Math. Soc.}, 77(3):395--400, 1979.

\bibitem{Donnelly-Xavier(1984)}
H.~Donnelly and F.~Xavier.
\newblock On the differential form spectrum of negatively curved {R}iemannian
  manifolds.
\newblock {\em Amer. J. Math.}, 106(1):169--185, 1984.

\bibitem{Eckmann(2002a)}
B.~Eckmann.
\newblock Lattices and {$l_2$}-{B}etti numbers.
\newblock Preprint, Z{\"u}rich, 2002.

\bibitem{Eckmann(2002b)}
B.~Eckmann.
\newblock Lattices with vanishing first {$L^2$}-Betti number and deficiency
  one.
\newblock Preprint, Z{\"u}rich, 2002.

\bibitem{Eckmann(1997)}
B.~Eckmann.
\newblock $4$-manifolds, group invariants, and $l\sb 2$-{B}etti numbers.
\newblock {\em Enseign. Math. (2)}, 43(3-4):271--279, 1997.

\bibitem{Eckmann(2000)}
B.~Eckmann.
\newblock Introduction to $l\sb 2$-methods in topology: reduced $l\sb
  2$-homology, harmonic chains, $l\sb 2$-{B}etti numbers.
\newblock {\em Israel J. Math.}, 117:183--219, 2000.
\newblock Notes prepared by G.~Mislin.

\bibitem{Efremov(1991a)}
A.~Efremov.
\newblock Combinatorial and analytic {N}ovikov {S}hubin invariants.
\newblock Preprint, 1991.

\bibitem{Elek(2002)}
G.~Elek.
\newblock The rank of finitely generated modules over group rings.
\newblock Preprint, 2003.

\bibitem{Farber(1996)}
M.~Farber.
\newblock Homological algebra of {N}ovikov-{S}hubin invariants and {M}orse
  inequalities.
\newblock {\em Geom. Funct. Anal.}, 6(4):628--665, 1996.

\bibitem{Farber(1997)}
M.~Farber.
\newblock Geometry of growth: approximation theorems for $L^2$-invariants.
\newblock {\em Math. Ann.}, 311(2):335--375, 1998.

\bibitem{Farber-Weinberger(2001)}
M.~Farber and S.~Weinberger.
\newblock On the zero-in-the-spectrum conjecture.
\newblock {\em Ann. of Math. (2)}, 154(1):139--154, 2001.

\bibitem{Furman(1999a)}
A.~Furman.
\newblock Gromov's measure equivalence and rigidity of higher rank lattices.
\newblock {\em Ann. of Math. (2)}, 150(3):1059--1081, 1999.

\bibitem{Furman(1999b)}
A.~Furman.
\newblock Orbit equivalence rigidity.
\newblock {\em Ann. of Math. (2)}, 150(3):1083--1108, 1999.

\bibitem{Gaboriau(2001)}
D.~Gaboriau.
\newblock Invariants $l^2$ de relation d'\'equivalence et de groupes.
\newblock Preprint, Lyon, 2001.

\bibitem{Gaboriau(2001b)}
D.~Gaboriau.
\newblock On orbit equivalence of measure preserving actions.
\newblock Preprint, Lyon, 2001.

\bibitem{Grigorchuk-Linnell-Schick-Zuk(2000)}
R.~I. Grigorchuk, P.~A. Linnell, T.~Schick, and A.~{\.Z}uk.
\newblock On a question of {A}tiyah.
\newblock {\em C. R. Acad. Sci. Paris S\'er. I Math.}, 331(9):663--668, 2000.

\bibitem{Grigorchuk-Zuk(2001)}
R.~I. Grigorchuk and A.~{\.Z}uk.
\newblock The {L}amplighter group as a group generated by a 2-state automaton,
  and its spectrum.
\newblock {\em Geom. Dedicata}, 87(1-3):209--244, 2001.

\bibitem{Gromov(1982)}
M.~Gromov.
\newblock Volume and bounded cohomology.
\newblock {\em Inst. Hautes \'Etudes Sci. Publ. Math.}, (56):5--99 (1983),
  1982.

\bibitem{Gromov(1986a)}
M.~Gromov.
\newblock Large {R}iemannian manifolds.
\newblock In {\em Curvature and topology of Riemannian manifolds (Katata,
  1985)}, pages 108--121. Springer-Verlag, Berlin, 1986.

\bibitem{Gromov(1987)}
M.~Gromov.
\newblock Hyperbolic groups.
\newblock In {\em Essays in group theory}, pages 75--263. Springer-Verlag, New
  York, 1987.

\bibitem{Gromov(1991)}
M.~Gromov.
\newblock K\"ahler hyperbolicity and ${L}\sb 2$-{H}odge theory.
\newblock {\em J. Differential Geom.}, 33(1):263--292, 1991.

\bibitem{Gromov(1993)}
M.~Gromov.
\newblock Asymptotic invariants of infinite groups.
\newblock In {\em Geometric group theory, Vol.\ 2 (Sussex, 1991)}, pages
  1--295. Cambridge Univ. Press, Cambridge, 1993.

\bibitem{Gromov(1999a)}
M.~Gromov.
\newblock {\em Metric structures for {R}iemannian and non-{R}iemannian spaces}.
\newblock Birkh\"auser Boston Inc., Boston, MA, 1999.
\newblock Based on the 1981 French original [MR 85e:53051], with appendices by
  M.\ Katz, P.\ Pansu and S.\ Semmes, Translated from the French by
  S.~M.~Bates.

\bibitem{Gromov-Shubin(1991)}
M.~Gromov and M.~A. Shubin.
\newblock Von {N}eumann spectra near zero.
\newblock {\em Geom. Funct. Anal.}, 1(4):375--404, 1991.

\bibitem{Hausmann-Weinberger(1985)}
J.-C. Hausmann and S.~Weinberger.
\newblock Caract\'eristiques d'{E}uler et groupes fondamentaux des vari\'et\'es
  de dimension $4$.
\newblock {\em Comment. Math. Helv.}, 60(1):139--144, 1985.

\bibitem{Hempel(1976)}
J.~Hempel.
\newblock {\em $3$-{M}anifolds}.
\newblock Princeton University Press, Princeton, N. J., 1976.
\newblock Ann. of Math. Studies, No. 86.

\bibitem{Hess-Schick(1998)}
E.~Hess and T.~Schick.
\newblock ${L}\sp 2$-torsion of hyperbolic manifolds.
\newblock {\em Manuscripta Math.}, 97(3):329--334, 1998.

\bibitem{Higson-Roe-Schick(2001)}
N.~Higson, J.~Roe, and T.~Schick.
\newblock Spaces with vanishing $l\sp 2$-homology and their fundamental groups
  (after {F}arber and {W}einberger).
\newblock {\em Geom. Dedicata}, 87(1-3):335--343, 2001.

\bibitem{Hillman(1999)}
J.~A. Hillman.
\newblock Deficiencies of lattices in connected {L}ie groups.
\newblock Preprint, 1999.

\bibitem{Hitchin(1974a)}
N.~Hitchin.
\newblock Compact four-dimensional {E}instein manifolds.
\newblock {\em J. Differential Geometry}, 9:435--441, 1974.

\bibitem{Hog-Lustig-Metzler(1985)}
C.~Hog, M.~Lustig, and W.~Metzler.
\newblock Presentation classes, $3$-manifolds and free products.
\newblock In {\em Geometry and topology (College Park, Md., 1983/84)}, pages
  154--167. Springer-Verlag, Berlin, 1985.

\bibitem{Howie(1982)}
J.~Howie.
\newblock On locally indicable groups.
\newblock {\em Math. Z.}, 180(4):445--461, 1982.

\bibitem{Ivanov(1987)}
N.~Ivanov.
\newblock Foundations of the theory of bounded cohomology.
\newblock {\em J. Soviet Math.}, 37:1090--1114, 1987.

\bibitem{Johnson-Kotschick(1993)}
F.~E.~A. Johnson and D.~Kotschick.
\newblock On the signature and {E}uler characteristic of certain
  four-manifolds.
\newblock {\em Math. Proc. Cambridge Philos. Soc.}, 114(3):431--437, 1993.

\bibitem{Jost-Xin(2000)}
J.~Jost and Y.~L. Xin.
\newblock Vanishing theorems for ${L}\sp 2$-cohomology groups.
\newblock {\em J. Reine Angew. Math.}, 525:95--112, 2000.

\bibitem{Kadison-Ringrose(1986)}
R.~V. Kadison and J.~R. Ringrose.
\newblock {\em Fundamentals of the theory of operator algebras. {V}ol. {I}{I}}.
\newblock Academic Press Inc., Orlando, FL, 1986.
\newblock Advanced theory.

\bibitem{Kaniuth(1969)}
E.~Kaniuth.
\newblock Der {T}yp der regul\"aren {D}arstellung diskreter {G}ruppen.
\newblock {\em Math. Ann.}, 182:334--339, 1969.

\bibitem{Knapp(1986)}
A.~W. Knapp.
\newblock {\em Representation theory of semisimple groups}.
\newblock Princeton University Press, Princeton, NJ, 1986.
\newblock An overview based on examples.

\bibitem{Kotschick(1994)}
D.~Kotschick.
\newblock Four-manifold invariants of finitely presentable groups.
\newblock In {\em Topology, geometry and field theory}, pages 89--99. World
  Sci. Publishing, River Edge, NJ, 1994.

\bibitem{Linnell(1992)}
P.~A. Linnell.
\newblock Zero divisors and ${L}\sp 2({G})$.
\newblock {\em C. R. Acad. Sci. Paris S\'er. I Math.}, 315(1):49--53, 1992.

\bibitem{Linnell(1993)}
P.~A. Linnell.
\newblock Division rings and group von {N}eumann algebras.
\newblock {\em Forum Math.}, 5(6):561--576, 1993.

\bibitem{Lott(1992a)}
J.~Lott.
\newblock Heat kernels on covering spaces and topological invariants.
\newblock {\em J. Differential Geom.}, 35(2):471--510, 1992.

\bibitem{Lott(1996b)}
J.~Lott.
\newblock The zero-in-the-spectrum question.
\newblock {\em Enseign. Math. (2)}, 42(3-4):341--376, 1996.

\bibitem{Lott(1999c)}
J.~Lott.
\newblock Deficiencies of lattice subgroups of {L}ie groups.
\newblock {\em Bull. London Math. Soc.}, 31(2):191--195, 1999.

\bibitem{Lott-Lueck(1995)}
J.~Lott and W.~L{\"u}ck.
\newblock ${L}\sp 2$-topological invariants of $3$-manifolds.
\newblock {\em Invent. Math.}, 120(1):15--60, 1995.

\bibitem{Lubotzky(1983)}
A.~Lubotzky.
\newblock Group presentation, $p$-adic analytic groups and lattices in
  ${SL}\sb{2}(\mathbb{{C}})$.
\newblock {\em Ann. of Math. (2)}, 118(1):115--130, 1983.

\bibitem{Lueck(1989)}
W.~L{\"u}ck.
\newblock {\em Transformation groups and algebraic ${K}$-theory}.
\newblock Springer-Verlag, Berlin, 1989.
\newblock Mathematica Gottingensis.

\bibitem{Lueck(1994c)}
W.~L{\"u}ck.
\newblock Approximating ${L}\sp 2$-invariants by their finite-dimensional
  analogues.
\newblock {\em Geom. Funct. Anal.}, 4(4):455--481, 1994.

\bibitem{Lueck(1994a)}
W.~L{\"u}ck.
\newblock ${L}\sp 2$-torsion and $3$-manifolds.
\newblock In {\em Low-dimensional topology (Knoxville, TN, 1992)}, pages
  75--107. Internat. Press, Cambridge, MA, 1994.

\bibitem{Lueck(1997a)}
W.~L{\"u}ck.
\newblock Hilbert modules and modules over finite von {N}eumann algebras and
  applications to ${L}\sp 2$-invariants.
\newblock {\em Math. Ann.}, 309(2):247--285, 1997.

\bibitem{Lueck(1997c)}
W.~L{\"u}ck.
\newblock ${L}\sp 2$-{I}nvarianten von {M}annigfaltigkeiten und {G}ruppen.
\newblock {\em Jahresber. Deutsch. Math.-Verein.}, 99(3):101--109, 1997.

\bibitem{Lueck(1998a)}
W.~L{\"u}ck.
\newblock Dimension theory of arbitrary modules over finite von {N}eumann
  algebras and ${L}\sp 2$-{B}etti numbers. {I}. {F}oundations.
\newblock {\em J. Reine Angew. Math.}, 495:135--162, 1998.

\bibitem{Lueck(1998b)}
W.~L{\"u}ck.
\newblock Dimension theory of arbitrary modules over finite von {N}eumann
  algebras and ${L}\sp 2$-{B}etti numbers. {I}{I}. {A}pplications to
  {G}rothendieck groups, ${L}\sp 2$-{E}uler characteristics and {B}urnside
  groups.
\newblock {\em J. Reine Angew. Math.}, 496:213--236, 1998.

\bibitem{Lueck(2001c)}
W.~L{\"u}ck.
\newblock ${L}\sp 2$-invariants and their applications to geometry, group
  theory and spectral theory.
\newblock In {\em Mathematics Unlimited---2001 and Beyond}, pages 859--871.
  Springer-Verlag, Berlin, 2001.

\bibitem{Lueck(2002a)}
W.~L{\"u}ck.
\newblock {$L\sp 2$}-invariants of regular coverings of compact manifolds and
  {CW}-complexes.
\newblock In {\em Handbook of geometric topology}, pages 735--817.
  North-Holland, Amsterdam, 2002.

\bibitem{Lueck(2002)}
W.~L{\"u}ck.
\newblock {\em {$L\sp 2$}-invariants: theory and applications to geometry and
  {$K$}-theory}, volume~44 of {\em Ergebnisse der Mathematik und ihrer
  Grenzgebiete}. 3. Folge. 
\newblock Springer-Verlag, Berlin, 2002.

\bibitem{Lueck-Reich(2003)}
W.~L{\"u}ck and H.~Reich.
\newblock The {B}aum-{C}onnes and the {F}arrell-{J}ones conjectures in ${K}$-
  and ${L}$-theory.
\newblock preprint, will be submitted to the handbook of $K$-theory,
  2003.

\bibitem{Lueck-Reich-Schick(1999)}
W.~L{\"u}ck, H.~Reich, and T.~Schick.
\newblock Novikov-{S}hubin invariants for arbitrary group actions and their
  positivity.
\newblock In {\em Tel Aviv Topology Conference: Rothenberg Festschrift (1998)},
  pages 159--176. Amer. Math. Soc., Providence, RI, 1999.

\bibitem{Lueck-Roerdam(1993)}
W.~L{\"u}ck and M.~R{\o}rdam.
\newblock Algebraic ${K}$-theory of von {N}eumann algebras.
\newblock {\em $K$-Theory}, 7(6):517--536, 1993.

\bibitem{Lueck-Schick(1999)}
W.~L{\"u}ck and T.~Schick.
\newblock ${L^2}$-torsion of hyperbolic manifolds of finite volume.
\newblock {\em Geometric and Functional Analysis}, 9:518--567, 1999.

\bibitem{Lyndon-Schupp(1977)}
R.~C. Lyndon and P.~E. Schupp.
\newblock {\em Combinatorial group theory}.
\newblock Springer-Verlag, Berlin, 1977.
\newblock Ergebnisse der Mathematik und ihrer Grenzgebiete, Band 89.

\bibitem{Mathai(1992)}
V.~Mathai.
\newblock $L^2$-analytic torsion.
\newblock {\em J. Funct. Anal.}, 107(2):369--386, 1992.

\bibitem{Mathai(1998)}
V.~Mathai.
\newblock $L^2$-invariants of covering spaces.
\newblock In {\em Geometric analysis and Lie theory in mathematics and
  physics}, pages 209--242. Cambridge Univ. Press, Cambridge, 1998.

\bibitem{Mautner(1950)}
F.~I. Mautner.
\newblock The structure of the regular representation of certain discrete
  groups.
\newblock {\em Duke Math. J.}, 17:437--441, 1950.

\bibitem{Mishchenko-Fomenko(1979)}
A.~S. Mi{\v{s}}{\v{c}}enko and A.~T. Fomenko.
\newblock The index of elliptic operators over ${C}\sp{\ast} $-algebras.
\newblock {\em Izv. Akad. Nauk SSSR Ser. Mat.}, 43(4):831--859, 967, 1979.
\newblock English translation in \emph{Math. USSR-Izv.} 15 (1980), no. 1,
  87--112.

\bibitem{Morgan(1984)}
J.~W. Morgan.
\newblock On {T}hurston's uniformization theorem for three-dimensional
  manifolds.
\newblock In {\em The Smith conjecture (New York, 1979)}, pages 37--125.
  Academic Press, Orlando, FL, 1984.

\bibitem{Novikov-Shubin(1986b)}
S.~P. Novikov and M.~A. Shubin.
\newblock Morse inequalities and von {N}eumann ${II}\sb 1$-factors.
\newblock {\em Dokl. Akad. Nauk SSSR}, 289(2):289--292, 1986.

\bibitem{Novikov-Shubin(1986a)}
S.~P. Novikov and M.~A. Shubin.
\newblock Morse inequalities and von {N}eumann invariants of non-simply
  connected manifolds.
\newblock {\em Uspekhi. Matem. Nauk}, 41(5):222--223, 1986.
\newblock in Russian.

\bibitem{Olbrich(2000)}
M.~Olbrich.
\newblock ${L^2}$-invariants of locally symmetric spaces.
\newblock Preprint, G\"ottingen, 2000.

\bibitem{Ornstein-Weiss(1980)}
D.~S. Ornstein and B.~Weiss.
\newblock Ergodic theory of amenable group actions. {I}. {T}he {R}ohlin lemma.
\newblock {\em Bull. Amer. Math. Soc. (N.S.)}, 2(1):161--164, 1980.

\bibitem{Pansu(1995)}
P.~Pansu.
\newblock Cohomologie ${L^p}$: Invariance sous quasiisometries.
\newblock Preprint, Orsay, 1995.

\bibitem{Pansu(1996a)}
P.~Pansu.
\newblock Introduction to $L^2$-{B}etti numbers.
\newblock In {\em Riemannian geometry (Waterloo, ON, 1993)}, pages 53--86.
  Amer. Math. Soc., Providence, RI, 1996.

\bibitem{Ramachandran(1993)}
M.~Ramachandran.
\newblock Von {N}eumann index theorems for manifolds with boundary.
\newblock {\em J. Differential Geom.}, 38(2):315--349, 1993.

\bibitem{Reich(1999)}
H.~Reich.
\newblock {\em Group von {N}eumann algebras and related algebras}.
\newblock PhD thesis, Universit\"at G\"ottingen, 1999.
\newblock http://www.math.uni-muenster.de/u/lueck/publ/diplome/reich.dvi.

\bibitem{Reich(2001)}
H.~Reich.
\newblock On the {$K$}- and {$L$}-theory of the algebra of operators affiliated
  to a finite von {N}eumann algebra.
\newblock {\em $K$-theory}, 24:303--326, 2001.

\bibitem{Sauer(2002)}
R.~Sauer.
\newblock {\em ${L}^2$-Invariants of groups and discrete measured groupoids}.
\newblock PhD thesis, Universit\"at M\"unster, 2002.

\bibitem{Schick(2000c)}
T.~Schick.
\newblock Integrality of ${L}\sp 2$-{B}etti numbers.
\newblock {\em Math. Ann.}, 317(4):727--750, 2000.

\bibitem{Schick(2001b)}
T.~Schick.
\newblock ${L}\sp 2$-determinant class and approximation of ${L}\sp 2$-{B}etti
  numbers.
\newblock {\em Trans. Amer. Math. Soc.}, 353(8):3247--3265 (electronic), 2001.

\bibitem{Thoma(1964)}
E.~Thoma.
\newblock \"{U}ber unit\"are {D}arstellungen abz\"ahlbarer, diskreter
  {G}ruppen.
\newblock {\em Math. Ann.}, 153:111--138, 1964.

\bibitem{Dieck(1987)}
T.~tom Dieck.
\newblock {\em Transformation groups}.
\newblock Walter de Gruyter \& Co., Berlin, 1987.

\bibitem{Wall(1961)}
C.~T.~C. Wall.
\newblock Rational {E}uler characteristics.
\newblock {\em Proc. Cambridge Philos. Soc.}, 57:182--184, 1961.

\bibitem{Wegner(2000)}
C.~Wegner.
\newblock {\em ${L^2}$-invariants of finite aspherical ${CW}$-complexes with
  fundamental group containing a non-trivial elementary amenable normal
  subgroup}.
\newblock PhD thesis, Westf\"alische Wilhelms-Universit\"at M\"unster, 2000.

\bibitem{Wegner(2001)}
C.~Wegner.
\newblock ${L^2}$-invariants of finite aspherical ${CW}$-complexes.
\newblock Preprintreihe SFB 478 --- Geometrische Strukturen in der Mathematik,
  Heft 152, M\"unster, 2001.

\bibitem{Weibel(1994)}
C.~A. Weibel.
\newblock {\em An introduction to homological algebra}.
\newblock Cambridge University Press, Cambridge, 1994.

\bibitem{Whyte(1999)}
K.~Whyte.
\newblock Amenability, bi-{L}ipschitz equivalence, and the von {N}eumann
  conjecture.
\newblock {\em Duke Math. J.}, 99(1):93--112, 1999.

\bibitem{Yu(1998a)}
G.~Yu.
\newblock The {N}ovikov conjecture for groups with finite asymptotic dimension.
\newblock {\em Ann. of Math. (2)}, 147(2):325--355, 1998.

\bibitem{Zucker(1983)}
S.~Zucker.
\newblock ${L}\sb{2}$-cohomology and intersection homology of locally symmetric
  varieties.
\newblock In {\em Singularities, Part 2 (Arcata, Calif., 1981)}, pages
  675--680. Amer. Math. Soc., Providence, RI, 1983.


\end{thebibliography}
\end{document}